\documentclass[10pt,a4paper]{article}
\usepackage{amsmath}
\usepackage{amssymb}
\usepackage{color}

\usepackage{graphicx}
\usepackage{amsfonts,amsmath, amssymb}

\numberwithin{equation}{section}

\usepackage{graphics}
\usepackage{epsfig}
\newtheorem{claim}{\bf \t}[part]

\newtheorem{theorem}{Theorem}[section]

\newtheorem{lemma}[theorem]{Lemma}
\newtheorem{proposition}[theorem]{Proposition}
\newtheorem{remark}[theorem]{Remark}

\def\v{\varepsilon}

\def\t{\theta}

\def\m{\mu}
\def\a{\alpha}
\def\b{\beta}
\def\g{\gamma}
\def\d{\delta}
\def\l{\lambda}

\def\r{\rho}

\def\Om{\Omega}

\def\i{\infty}
\def\f{\frac}
\def\MZ{\mathcal{Z}}

\begin{document}

\title{ Uniform regularity and vanishing viscosity limit for the compressible Navier-Stokes with general  Navier-slip boundary conditions in  3-dimensional domains}

\author{ {Yong Wang$^{\S}$\footnote{ \indent
Email addresses:  yongwang@amss.ac.cn (Yong Wang), zpxin@ims.cuhk.edu.hk (Zhouping Xin), yongyan@usst.edu.cn (Yan Yong)
},~~~ Zhouping Xin$^{\ddag}$, ~~~Yan Yong$^\dag$
}
\\
\ \\
   {\small \it $^\S$Institute of Applied Mathematics, AMSS, CAS, Beijing 100190, China}\\
\  \\
   {\small \it $^\ddag$The Institute of Mathematical Sciences, The Chinese University of Hong Kong, Shatin, Hong Kong } \\
   \  \\
   {\small \it $^\dag$College of Science, University of Shanghai for Science and Technology, Shanghai 200093, China } \\
}

\date{ }
\maketitle

\begin{abstract}

In this paper, we investigate the uniform regularity for the isentropic compressible Navier-Stokes system with general  Navier-slip boundary conditions \eqref{1.10} and the inviscid limit to the compressible Euler system. It is shown that there exists a unique strong solution of the compressible Navier-Stokes equations with general  Navier-slip boundary conditions in  an interval of time which is uniform in the vanishing viscosity limit. The solution is uniformly bounded in a conormal Sobolev space and is uniform  bounded in $W^{1,\infty}$. It is also shown that the boundary layer for the density is weaker than the one for the velocity field.   In particular, it is  proved that the  velocity will be uniform bounded in $L^\infty(0,T;H^2)$ when the boundary is flat and  the Navier-Stokes system is supplemented with the special boundary condition \eqref{10.1}.   Based on such uniform estimates, we prove the  convergence of the viscous solutions to the inviscid ones in  $L^\infty(0,T;L^2)$, $L^\infty(0,T;H^1)$ and $L^\infty([0,T]\times\Omega)$ with a rate of convergence.

\


 \

AMS: 35Q35, 35B65, 76N10
\end{abstract}

\tableofcontents

\section{Introduction}
In this paper, we consider the isentropic compressible Navier-Stokes equations
\begin{equation}\label{1.1}
\begin{cases}
\r_t^\v+\mbox{div}(\r^\v u^\v)=0,\\
\r^\v u^\v_t+\r^\v u^\v\cdot\nabla u^\v+\nabla p^\v=\mu\v \Delta u^\v+(\mu+\l)\v \nabla div u^\v,
\end{cases}
x\in \Omega,~ t>0
\end{equation}
where $\Omega$ is a bounded smooth domain of $\mathbb{R}^3$,  $\r^\v, u^\v$ represent the density and velocity, respectively, $p^\v=p(\r^\v)$ is the pressure function   given by $\g$-law
\begin{equation*}
p(\r)=\r^{\g},~ \mbox{with}~\g>1.
\end{equation*}
The viscous coefficients  $\mu\v$ and $ \l\v$  satisfy the physical restrictions
\begin{equation}
\mu>0,~~2\m+3\l>0,
\end{equation}
where the parameter $\v>0$ is the inverse of the Reynolds number.

Here, we are interested in the existence of strong solution of \eqref{1.1} with uniform bounds on an interval of time independent of viscosity $\v\in(0,1]$ and the vanishing viscosity limit to the corresponding Euler equations as $\v$ vanishes, i.e,
\begin{equation}\label{1.7}
\begin{cases}
\r_t+\mbox{div}(\r u)=0,\\
(\r u)_t+\mbox{div}(\r u \otimes u)+\nabla p=0.
\end{cases}
\end{equation}

\

There has lots of literature on the uniform bounds and the vanishing viscosity limit when the domain has no boundaries, see for instances \cite{Constantin,Constantin-1,Kato-1,Masmoudi}.
However, in the presence of   physical boundaries, the problems become much more complicated and  challenging due to the possible appearance of  boundary layers. Indeed, in presence of a boundary, one of the most important physical boundary conditions for the Euler equations is the   slip boundary condition, i.e,
\begin{equation}\label{1.8}
u\cdot n=0, ~~\partial\Omega,
\end{equation}
and there exists a unique smooth solution for the initial boundary value problem \eqref{1.7} and \eqref{1.8}  at least locally in time. This boundary condition is  characteristic for the Euler equations \eqref{1.7}. Corresponding to  \eqref{1.8},   there are different choices of boundary conditions for the Navier-Stokes equations, and the no-slip boundary condition
\begin{equation*}\label{1.9}
u^\v=0,~~\mbox{on}~~\partial\Omega,
\end{equation*}
is one of the frequently used one.  Another one is the well-known Navier-slip boundary condition, i.e,
\begin{equation}\label{1.2}
u^\v\cdot n=0, ~~(Su^\v\cdot n)_{\tau}=-\a u^\v_{\tau},~~x\in\partial\Omega,
\end{equation}
where $n$ is the outward  unit normal to $\partial\Omega$, $u_\tau$ represents the tangential part of $u$. $S$ is the strain tensor
\begin{equation*}
Su=\f12(\nabla u+ \nabla u^t).
\end{equation*}
The boundary condition \eqref{1.2}, which was introduced by Navier \cite{Navier},  expresses  that the velocity on the boundary is propositional to the tangential component of the stress. This kind of boundary condition allows the fluid to slip at the boundary, and has important applications for problems with rough boundaries.

The Navier-slip boundary condition \eqref{1.2} can be written to the following generalized one
\begin{equation}\label{1.10}
u^\v\cdot n=0, ~~(Su^\v\cdot n)_{\tau}=-(Au^\v)_{\tau},~~x\in\partial\Omega,
\end{equation}
with $A$ a smooth symmetric matrix ,see \cite{Gie-Kelliher}. For smooth solutions, it is noticed that
\begin{equation*}\label{2.6}
\left(2S(v)n-(\nabla\times v)\times n\right)_\tau=-(2S(n)v)_\tau,
\end{equation*}
see \cite{Xiao-Xin-3} for details. Therefore, as  in \cite{Xiao-Xin-1,Xiao-Xin-2},  the boundary condition \eqref{1.10} can  be rewritten in the form of the vorticity as
\begin{equation}\label{2.7}
u^\v\cdot n=0, ~~n\times\omega^\v=[Bu^\v]_\tau,~~x\in\partial\Omega,
\end{equation}
where $\omega^\v=\nabla\times u^\v$ is the vorticity   and  $B=2(A-S(n))$ is a symmetric matrix. Actually,  it turns out  that the form \eqref{2.7} will be  more convenient than \eqref{1.10} in the energy estimates, see \cite{Xiao-Xin-1}.

For the incompressible fluid, the vanishing viscosity limit of the incompressible Navier-Stokes with no-slip boundary condition to the incompressible Euler flows with boundary condition \eqref{1.8} is one of the major open problems due to the possible appearance of   boundary layers, as illustrated   by Prandtl's theory. In \cite{S-C-1,S-C-2}, the authors  proved the(local in time) convergence of the incompressible Navier-Stokes flows to the Euler flows outside the boundary layer and to the prandtl flows in the boundary layer at the inviscid limit for the analytic initial data. Recently, Y. Maekawa \cite{Maekawa} proved this limit when the initial vorticity is located away from the boundary  in 2-D half plane.

On the other hand, for the incompressible Navier-Stokes system with Navier-slip boundary condition
\eqref{1.2}, considerable progress has been made on this problem. Indeed, the uniform $H^3$ bound and a uniform existence time interval as $\v$ tends to zero are obtained by Xiao-Xin in \cite{Xiao-Xin-1} for flat boundaries, which are generalized to $W^{k,p}$ in \cite{Beiro-1,Beiro-2}. However, such results can not be expected for general curved boundaries since boundary layer may appear due to non-trivial curvature as pointed out in \cite{Iftimie}. In such a case, Iftimie and Sueur have proved the convergence of the viscous solutions to the inviscid Euler solutions in $L^\infty(0,T,; L^2)$-space by a careful construction of boundary layer expansions and energy estimates. However, to identify precisely the asymptotic structure and get the convergence in stronger norms such as $L^\infty(0,T; H^s)(s>0)$, further a priori estimates and analysis are needed. Recently, Masmoudi-Rousset \cite{Masmoudi-R} established conormal uniform estimates for 3-dimensional general smooth domains with the Naiver-slip boundary condition \eqref{1.2}, which, in particular, implies the uniform boundedness of the normal first order derivatives of the velocity field. This allows the authors(\cite{Masmoudi-R}) to obtain the convergence of the viscous solutions to the inviscid ones by a compact argument. Based on the uniform estimates in \cite{Masmoudi-R}, better convergence with rates have been studied in \cite{Gie-Kelliher} and \cite{Xiao-Xin-2}. In particular, Xiao-Xin  \cite{Xiao-Xin-2} has proved the convergence in $L^\infty(0,T; H^1)$ with a rate of convergence.

For the compressible Navier-Stokes equations, however, the study is quite limited. Xin and Yanagisawa \cite{Xin-Y} studied the vanishing viscosity limit of the linearized compressible Navier-Stokes system  with the no-slip boundary condition in the 2-D half plane. Recently, Wang and Williams \cite{Wang-Williams} constructed a boundary layer solution of the compressible Navier-Stokes equations with Navier-slip boundary conditions in 2-D half plane. The layers constructed in \cite{Wang-Williams} are of width $O(\sqrt\v)$ as the Prandtl boundary layer, but are of amplitude $O(\sqrt\v)$ which is similar to the one \cite{Iftimie} for the incompressible case. So, in general, it is impossible to obtain the $H^3$ or $W^{2,p}(p>3)$ estimates for the compressible Navier-Stokes system \eqref{1.1} with the generalized Navier-slip boundary condition \eqref{1.10} or \eqref{2.7}. Recently,  Paddick \cite{Paddick} obtained an uniform estimates for the solutions of the compressible isentropic Navier-Stokes system in the 3-D half-space with a Navier boundary condition.  As expected, the boundary layers for the density must be weaker than the one for the velocity, however, this has not been proved  in \cite{Paddick}.

In the present paper,  we aim to obtain the uniform estimates in some anisotropic conormal Sobolev spaces and a control of the Lipschitz norm for solutions of the compressible Navier-Stokes equations \eqref{1.1} with the Navier-slip boundary condition \eqref{1.10} in general $3$-dimensional domains. As a consequence, our uniform estimates will yield that the boundary layers for the density are weaker than the one for the velocity. Furthermore, we obtain an uniform estimate in $L^\infty(0,T;H^2)$ when the boundary is flat. Finally, we study the vanishing viscosity limit of viscous solutions to the inviscid ones with a rate of  convergence.  Since the divergence free condition plays a key role in the analysis of \cite{Masmoudi-R}, delicate estimates for $\mbox{div}u$ are needed to  complete the analysis for the compressible Navier-Stokes system. Moreover, the compressible Navier-Stokes system is much more complicated to handle than the incompressible one.

\

The bounded domain $\Omega\subset\mathbb{R}^3$ is assumed to have   a covering such that
\begin{equation}\label{2.0}
\Omega\subset\Omega_0\cup_{k=1}^n\Omega_k,
\end{equation}
where $\overline{\Omega}_0\subset\Omega$ and in each $\Omega_k$ there exists a function $\psi_k$ such that
\begin{equation}
\Omega\cap\Omega_k=\{x=(x_1,x_2,x_3)~|~x_3>\psi_k(x_1,x_2)\}\cap\Omega_k
~~\mbox{and}~~\partial\Omega\cap\Omega_k=\{x_3=\psi_k(x_1,x_2)\}\cap\Omega_k.\nonumber
\end{equation}
$\Omega$ is said to be  $\mathcal{C}^m$ if the functions $\psi_k$ are $\mathcal{C}^m$-function.

To define the Sobolev conormal spaces, we consider $(Z_k)_{1\leq k\leq N}$ a finite set of generators of vector fields that are tangent to $\partial\Omega$ and set
\begin{equation*}\label{1.3}
H^m_{co}=\Big\{f\in L^2(\Omega)~|~Z^If\in L^2(\Omega), ~~\mbox{for}~|I|\leq m  \Big\},
\end{equation*}
where $I=(k_1,\cdots, k_m)$. We will  use the following notations
\begin{equation*}\label{1.4}
\|u\|^2_{m}=\|u\|^2_{H^m_{co}}=\sum_{j=1}^3\sum_{|I|\leq m}\|Z^Iu_j\|^2_{L^2},
\end{equation*}
\begin{equation*}\label{1.5}
\|u\|^2_{m,\infty}=\sum_{|I|\leq m}\|Z^Iu\|^2_{L^\infty},
\end{equation*}
and
\begin{equation*}\label{1.6}
\|\nabla Z^m u\|^2=\sum_{|I|=m}\|\nabla Z^Iu\|^2_{L^2}.
\end{equation*}
Noting that by using the covering of $\Om$, one can always assume that each vector field is supported in one of the $\Om_i$, moreover, in $\Om_0$ the norm $\|\cdot\|_m$ yields a control of the standard $H^m$ norm, whereas if $\Om_i\cap \partial\Om\neq \O$, there is no control of the normal derivatives.

Denote by  $C_k$   a positive constant independence of $\v\in(0,1] $ which depends only on the $\mathcal{C}^k$-norm of the functions $\psi_j$. Since $\partial\Omega$ is given locally by $x_3=\psi(x_1,x_2)$(we omit the subscript $j$ for notational convenience), it is convenient to use the coordinates:
\begin{equation*}\label{3.5}
\Psi:~(y,z)\longmapsto (y,\psi(y)+z)=x.
\end{equation*}
A local basis is thus given by the vector fields $(\partial_{y^1},\partial_{y^2},\partial_z)$. On the boundary $\partial_{y^1}$ and $\partial_{y^2}$ are tangent to $\partial\Omega$, and in general, $\partial_z$ is not a normal vector field. By using this parametrization, one can take as suitable vector fields compactly supported in $\Omega_j$ in the definition of the $\|\cdot\|_m$ norms:
\begin{equation*}\label{3.6}
Z_i=\partial_{y^i}=\partial_i+\partial_i\psi\partial_z,~i=1,2, ~~Z_3=\varphi(z)\partial_z,
\end{equation*}
where $\varphi(z)=\f{z}{1+z}$ is smooth, supported in $\mathbb{R}_+$ with the property $\varphi(0)=0$,
$\varphi'(0)>0$, $\varphi(z)>0$ for $z>0$. It is easy to check that
\begin{equation*}\label{3.7}
Z_kZ_j=Z_jZ_k,~~j,~k=1,2,3,
\end{equation*}
and
\begin{equation*}\label{3.8}
\partial_zZ_i=Z_i\partial_z,~i=1,2,~~\mbox{and}~~\partial_zZ_3\neq Z_3\partial_z.
\end{equation*}

In this paper, we shall still denote by $\partial_j,~j=1,2,3$ or $\nabla$ the derivatives in the physical space. The coordinates of a vector field $u$ in the basis  $(\partial_{y^1},\partial_{y^2},\partial_z)$ will be denoted by $u^i$, thus
\begin{equation}\label{3.9}
u=u^1\partial_{y^1}+u^2\partial_{y^2}+u^3\partial_{z}.
\end{equation}
We shall denote by $u_j$ the coordinates in the standard basis of $\mathbb{R}^3$, i.e, $u=u_1\partial_1+u_2\partial_2+u_3\partial_3$. Denote by $n$ the unit outward normal in the physical space which is given locally by
\begin{equation*}\label{3.10}
n(x)\equiv n(\Psi(y,z))=\f{1}{\sqrt{1+|\nabla\psi(y)|^2}}\left(\begin{array}{cccc} &\partial_1\psi(y)\\&\partial_2\psi(y)\\&-1\end{array}\right)\doteq\f{-N(y)}{\sqrt{1+|\nabla\psi(y)|^2}},
\end{equation*}
and by $\Pi$ the orthogonal projection
\begin{equation*}
\Pi(x)\equiv\Pi(\Psi(y,z))u=u-[u\cdot n(\Psi(y,z))]n(\Psi(y,z)).
\end{equation*}
which gives the orthogonal projection onto the tangent space of the boundary. Note that $n$ and $\Pi$ are defined in the whole $\Omega_k$ and do not depend on $z$.

For later use and notational convenience, we set
\begin{equation}\label{2.1}
\mathcal{Z}^\a=\partial_t^{\a_0}Z^{\a_1}=\partial_t^{\a_0}Z_1^{\a_{11}}Z_2^{\a_{12}}Z_3^{\a_{13}}.
\end{equation}
and use the following notations
\begin{equation}\label{2.3}
\|f(t)\|^2_{\mathcal{H}^m}=\sum_{|\a|\leq m}\|\mathcal{Z}^\a f(t)\|^2_{L^2_x},~~
\|f(t)\|_{\mathcal{H}^{k,\infty}}=\sum_{|\a|\leq k}\|\mathcal{Z}^\a f(t)\|^2_{L^\infty_x},
\end{equation}
for smooth space-time function $f(x,t)$. Throughout this paper, the positive
generic constants that are independent of $\v$ are denoted by
$c,C$.   $\|\cdot\|$  denotes the standard
$L^2(\Om;dx)$ norm, and $\|\cdot\|_{H^m}~(m=1,2,3,\cdots)$
denotes the Sobolev $H^m(\Om;dx)$ norm. The notation $|\cdot|_{H^m}$ will be used for the standard Sobolev norm of functions defined on $\partial\Omega$. Note that this norm involves only tangential  derivatives.  $P(\cdot)$ denotes a polynomial function.

\

Since the boundary layer may appear in the presence of physical boundaries, in order to obtain the uniform estimation for  solutions of the compressible Navier-Stokes system with Navier-slip boundary condition, one needs to find a suitable functional space. Here, we define the functional space $X_m^\v(T)$ for a pair of function $(p,u)=(p,u)(x,t)$  as follows:
\begin{eqnarray}
X^\v_m(T)=\Big\{(p,u)\in L^\infty([0,T],L^2);~~
\mbox{esssup}_{0\leq t\leq T}\|(p,u)(t)\|_{X_m^\v}<+\infty\Big\},
\end{eqnarray}
where the norm $\|(\cdot,\cdot)\|_{X_m^\v}$ is given by
\begin{eqnarray}
&&\|(p,u)(t)\|_{X^\v_m}=\|(p,u)(t)\|^2_{\mathcal{H}^m}+\|\nabla u(t)\|^2_{\mathcal{H}^{m-1}}+\sum_{k=0}^{m-2}\|\partial_t^k \nabla p(t)\|^2_{m-1-k}+\|\Delta p(t)\|^2_{\mathcal{H}^1}\nonumber\\
&&~~~~~~~~~~~~~~~~~~~~~+\|\nabla u\|^2_{\mathcal{H}^{1,\infty}}+\v\|\nabla\partial_t^{m-1}p(t)\|^2
+\v\|\Delta p(t)\|^2_{\mathcal{H}^2}.
\end{eqnarray}
In the present paper,  we supplement the compressible Navier-Stokes equations \eqref{1.1} with the initial data
\begin{equation}\label{1.11}
(\r^\v,u^\v)(x,0)=(\r_0^\v,u_0^\v)(x),
\end{equation}
such that
\begin{equation}\label{1.12}
0<\f1{C_0}\leq \r_0^\v\leq C_0<\infty,
\end{equation}
and
\begin{eqnarray}\label{1.13}
&&\sup_{0< \v\leq 1}\|(p_0^\v,u_0^\v)\|_{X^\v_m}=\sup_{0< \v\leq 1}\bigg\{ \|(p_0^\v,u_0^\v)\|^2_{\mathcal{H}^m}+\|\nabla u_0^\v\|^2_{\mathcal{H}^{m-1}}+\sum_{k=0}^{m-2}\|\partial_t^k \nabla p_0^\v\|^2_{m-1-k}\nonumber\\
&&~~~~~~~~~~~~~~+\|\Delta p_0^\v\|^2_{\mathcal{H}^1}+\|\nabla u_0^\v\|^2_{\mathcal{H}^{1,\infty}}+\v\|\nabla\partial_t^{m-1}p_0^\v\|^2
+\v\|\Delta p_0^\v\|^2_{\mathcal{H}^2}\bigg\}\leq \tilde{C}_0,
\end{eqnarray}
where $p_0^\v=p(\r_0^\v)$, $C_0>0$, $\tilde{C}_0>0$ are positive constants independent of $\v\in(0,1]$, and  the time derivatives of initial data  in \eqref{1.13} are defined through the compressible Navier-Stokes system \eqref{1.1}. Thus,   the initial data  $(\r^\v_0,u^\v_0)$  is assumed to have a higher space regularity and compatibilities. Notice that the {\it a priori} estimates in Theorem \ref{thm3.1} below is obtained in the case that the approximate solution is sufficient smooth up to the boundary, therefore, in order to obtain a selfcontained result, one  needs  to assume that the approximate initial data  satisfies the boundary compatibility conditions, i.e. \eqref{1.10}(or equivalent to \eqref{2.7}).  For the initial data $(\r_0^\v,u_0^\v)$  satisfying \eqref{1.13}, it is  not clear if there exists an approximate sequence  $(\r_0^{\v,\d},u_0^{\v,\d})$($\d$ being a regularization parameter), which satisfy the boundary compatibilities  and $\|(p_0^{\v,\d}-p_0^\v,u_0^{\v,\d}-u_0^\v)\|_{X^\v_m}\rightarrow0$ as $\d\rightarrow0$. Therefore, we set
\begin{eqnarray}\label{Initial}
&&X_{NS,ap}^{\v,m}=\Big\{(p,u)\in C^{2m}(\bar\Omega)~\Big| \partial_t^k p,~\partial_t^ku,k=1,\cdots,m~\mbox{are defined through the Navier-Stokes} \nonumber\\
&&~~~~~~~~~~~~~~~~~~~~~~~~~~~~~~~~~~~~~~~ \mbox{equations}~\eqref{1.1}~ \mbox{and}~\partial_t^ku,k=0,\cdots,m-1~ \mbox{satisfy }\nonumber\\
&&~~~~~~~~~~~~~~~~~~~~~~~~~~~~~~~~~~~~~~~~~\mbox{ the boundary compatibility condition}\Big\},
\end{eqnarray}
and
\begin{eqnarray}\label{Initial space}
&&X^{\v,m}_{NS}=\mbox{The closure of}~ X^{\v,m}_{NS,ap}~ \mbox{in the norm }~ \|(\cdot,\cdot)\|_{X^\v_m}.
\end{eqnarray}

Then  our main result in this paper is follows:
\begin{theorem}[Uniform Regularity]\label{thm1.1}
Let $m$ be an integer satisfying $m\geq 6$,  $\Omega$ be a $\mathcal{C}^{m+2}$ domain and $A\in\mathcal{C}^{m+1}(\partial\Omega)$. Consider the initial data $(p^\v_0, u^\v_0)\in X^{\v,m}_{NS}$ given in \eqref{1.11} and satisfying \eqref{1.12}
-\eqref{1.13}. Then there exists a time $T_0>0$ and $\tilde{C}_1>0$  independent of  $\v\in(0,1]$, such that there exists a unique solution $(\r^\v, u^\v)$ of \eqref{1.1}, \eqref{1.10}, \eqref{1.11} which is defined on $[0,T_0]$  and satisfies the estimates:
\begin{eqnarray}\label{1.18}
&&\sup_{0\leq t\leq T_0}\bigg\{ \|( u^\v, p^\v)(t)\|^2_{\mathcal{H}^m}+\|\nabla u^\v(t)\|^2_{\mathcal{H}^{m-1}}+\sum_{k=0}^{m-2}\|\partial_t^k \nabla p^\v(t)\|^2_{m-1-k}+\|\Delta p^\v(t)\|^2_{\mathcal{H}^1}\nonumber\\
&&~~~~~~~~~+\|\nabla u^\v(t)\|^2_{\mathcal{H}^{1,\infty}}+\v\|\nabla\partial_t^{m-1}p^\v(t)\|^2
+\v\|\Delta p^\v(t)\|^2_{\mathcal{H}^2}\bigg\}+\int_{0}^{T_0}\|\nabla\partial_t^{m-1}p^\v(t)\|^2dt\nonumber\\
&&~~~~~~~~~+\int_{0}^{T_0}\|\Delta p^\v(t)\|^2_{\mathcal{H}^2}dt
+\v\int_{0}^{T_0}\|\nabla u^\v(t)\|^2_{\mathcal{H}^m}dt
+\v\sum_{k=0}^{m-2}\int_{0}^{T_0}\|\nabla^2\partial_t^ku^\v(t)\|^2_{m-k-1}dt\nonumber\\
&&~~~~~~~~~
+\v^2\int_{0}^{T_0}\|\nabla^2\partial_t^{m-1}u^\v(t)\|^2dt
\leq \tilde{C}_1<\infty,
\end{eqnarray}
and
\begin{equation}\label{1.19}
\f1{2C_0}\leq \r^\v(t)\leq 2C_0~~\forall t\in [0,T_0],
\end{equation}
where $\tilde{C}_1$ depends only  on $C_0,~\tilde{C}_0$ and $C_{m+2}$.

\end{theorem}

\begin{remark}
Recently, we notice that Paddick \cite{Paddick} obtained a similar uniform estimates for the solutions of the  compressible isentropic Navier-Stokes system in  the 3-D half-space with a Navier boundary condition. However, the details of proof are different, and our regularity is  better than the one in \cite{Paddick}, especially, we show that $\|\Delta p^\v(t)\|^2_{\mathcal{H}^1}$ is uniform bounded which yields immediately  that the boundary layer for the density $\r^\v$ is weaker than the one for velocity $u^\v$ as expected.
\end{remark}

\begin{remark}
It is obvious  that $X_{NS}^{\v,m}\subset\{(p,u)\in L^2(\Omega)~|\partial^k_t(p,u)~\mbox{defined through}~ \eqref{1.1},
~\|(p,u)\|_{X^\v_m}<+\infty, 0\leq k\leq m\}$, yet it is not clear whether $"\subset"$ can be changed to $"="$. And we will not address this problem since our main concern is  the uniform regularity of the solution of Navier-Stokes equations. Here, it should be  pointed out that there are lots of data contained in $X_{NS}^{\v,m}$, for example, let $(\r_0^\v,u_0^\v)$ be sufficiently smooth functions, and in a vicinity of the boundary, $\r_0^\v$ is positive constant and $u_0^\v$ vanishes, then it is obvious that $(p(\r_0^\v),u_0^\v)\in X_{NS}^{\v,m}$.
\end{remark}

\begin{remark}
For $(p_0^\v,u_0^\v)\in X_{NS}^{\v,m}$, it must hold that $u^\v_0\cdot n|_{\partial\Omega}=0$ and $(Su^\v_0\cdot n)_{\tau}|_{\partial\Omega}=-(Au^\v_0)_{\tau}|_{\partial\Omega}$ in the trace sense for every fixed $\v\in(0,1]$.  For the solution $(\r^\v, u^\v)(t)$ of \eqref{1.1}, \eqref{1.10}, \eqref{1.11}, the boundary conditions \eqref{1.10} are satisfied in the trace sense for every fixed $\v\in(0,1]$ and $t\in (0,T_0]$.
\end{remark}

\begin{remark}
When time derivative is applied to the boundary layer, it has the same properties as the tangential derivatives. So,   the time derivative is regarded as a tangential derivative in this sense.
\end{remark}


We now outline the proof of Theorem \ref{thm1.1}. First, we   obtain a conormal energy estimates for $(p^\v,u^\v)$ in  $\mathcal{H}^m$-norm(see  \eqref{2.3} above for the definition of $\mathcal{H}^m$). Second, since the $\mbox{div}u^\v$ is no longer free for the compressible Navier-Stokes equations, one has to get enough estimates for $\mbox{div}u^\v$. Indeed,  we can obtain a control of  $\sum_{j=0}^{m-2}\|\partial_t^j(\mbox{div}u^\v,\nabla p^\v)\|^2_{m-1-j}$ at the cost that  the term $\int_0^t\|\nabla\MZ^{m-2}\mbox{div}u^\v\|^2d\tau$ appears in the right hand side of the inequality. And, in general, it is impossible to obtain the uniform bound of $\int_0^t\| \MZ^{m-2}\partial_{zz}u^\v\|^2d\tau$ due to the possible appearance of boundary layers. However, the situation is different for $\int_0^t\|\nabla\MZ^{m-2}\mbox{div}u^\v\|^2d\tau$, because $\mbox{div}u^\v$ is not expected to have boundary layer structure.  Another difficulty is that, due to the  singular behavior at the boundary, we can only obtain  the uniform estimate of $\v\|\partial_t^{m-1}(\mbox{div}u^\v,\nabla p^\v)\|^2$ which is not enough   to get the uniform estimate for $\|\nabla\partial_t^{m-1}u^\v\|$.  Fortunately, we can obtain the uniform estimates for $\int_0^t\|\partial_t^{m-1}\nabla p^\v\|^2$ and get a  control  of $\|\partial_t^{m-1}\mbox{div}u^\v\|$ in terms of  $\sum_{j=0}^{m-2}\|\partial_t^j(\nabla u^\v,\nabla p^\v)\|^2_{m-1-j}$ and  $\|(p^\v,u^\v)\|_{\mathcal{H}^m}$ which are independent of $\|\partial_t^{m-1}(\mbox{div}u^\v,\nabla p^\v)\|^2$. These key observations play an important role in this paper.  The third step is to estimate the $\|\partial_nu^\v\|_{\mathcal{H}^{m-1}}$. Similar to \cite{Masmoudi-R}, due to the Navier-slip condition \eqref{2.7}, it is convenient to study $\eta=\omega^\v\times n+(Bu^\v)_{\tau}$ with a  homogeneous Dirichlet boundary condition. Indeed, we  get a control of $\|\eta\|_{\mathcal{H}^{m-1}}$ by using energy estimates on the equations solved by $\eta$. The fourth step is to estimate $\|\nabla u^\v\|_{\mathcal{H}^{1,\infty}}$. In fact, it suffices to estimate $\|(\partial_n u^\v)_\tau\|_{\mathcal{H}^{1,\infty}}$ since the other terms can be estimated by the Sobolev imbedding. We choose an equivalent quantity such that it satisfies a homogeneous Dirichlet condition and solves a convection-diffusion equation at the leading order. Before performing the estimates, we generalize some results of \cite{Masmoudi-R} in the Appendix, so that it can be applied  to the compressible Navier-Stokes system.  Moreover, we also need to get some control on $\|\nabla\mbox{div}u^\v\|_{L^\infty}$.   Then, all these preparations will enable us to obtain  a control of  $\|\nabla u^\v\|_{\mathcal{H}^{1,\infty}}$. The last step is to obtain the uniform estimate of $\|\Delta p^\v\|_{\mathcal{H}^{1}}$ which gives a control of $\|\nabla p^\v\|_{\mathcal{H}^{1,\infty}}$  from Proposition \ref{prop3.3}. Then Theorem \ref{thm1.1} can be proved by   the above a priori  estimates and a classical iteration method.

\

In general, it is hard to obtain the uniform estimate of $\|u^\v\|_{L^\infty(0,T;H^2)}$ due to the possible boundary layers. However,  the uniform $H^3$ bound and a uniform existence time interval as $\v$ tends to zero are obtained by Xiao-Xin in \cite{Xiao-Xin-1}(which are generalized to $W^{k,p}$ in \cite{Beiro-1,Beiro-2}) when the boundary is  flat and the Navier-Stokes system is imposed with the following special Navier-slip boundary condition
\begin{align}\label{10.1}
n\cdot u^\v=0,~~~n\times\omega^\v=0, ~x\in\partial\Om.
\end{align}
In   Theorem \ref{thm1.3} below, we prove that $\|u^\v\|_{L^\infty(0,T;H^2)}$ is uniformly bounded for the solution of compressible Navier-Stokes system \eqref{1.1} when the boundary is flat and the special boundary condition \eqref{10.1} is imposed.

In order to avoid the unessential technical difficulties, without loss of generality,  we assume that the domain $\Om$ is given by
\begin{align}\label{10.2}
\Omega=\mathbb{T}^2\times(0,1),
\end{align}
and set
\begin{align}\label{10.3}
\Gamma=\{x=(y_1,y_2,z)~|~0\leq y_1,~y_2\leq 1,~\mbox{and}~~z=0~\mbox{or}~z=1 \}.
\end{align}
Then, the boundary condition \eqref{10.1} will be imposed on $\Gamma$. Hereafter, the {\it flat case} means that $\Om=\mathbb{T}^2\times(0,1)$ and the Navier-Stokes system is supplemented with the special Navier-slip boundary condition \eqref{10.1}. In this domain, we define  the conormal derivatives
as following
\begin{align}\label{10.4}
Z_i=\partial_{y_i},~i=1,2,~\mbox{and}~~Z_3=z(1-z)\partial_z.
\end{align}
Then,  we have  better uniform estimates for  $\|u^\v\|_{H^2}$ as follows:
\begin{theorem}[Flat case]\label{thm1.3}
	Let $m\geq6$ and $\Om=\mathbb{T}^2\times(0,1)$. Consider the initial data $(p^\v_0, u^\v_0)\in X^{\v,m}_{NS}\cap H^2$ given in \eqref{1.11} and satisfying \eqref{1.12}-\eqref{1.13}. Then there exists a time $T_0>0$ and $\tilde{C}_1>0$  independent of  $\v\in(0,1]$, such that there exists a unique solution $(\r^\v, u^\v)$ of \eqref{1.1}, \eqref{1.11}, \eqref{10.1} which is defined on $[0,T_0]$  and satisfies the uniform estimates \eqref{1.18} and \eqref{1.19}.
	Especially, it holds that
	\begin{align}\label{10.5}
	\sup_{0\leq t\leq T_0}\|u^\v(t)\|^2_{H^2}+\v\int_0^{T_0}\|u^\v(\tau)\|^2_{H^3}d\tau
	\leq \exp(\tilde{C}_1)(1+ \|u_0\|^2_{H^2}),
	\end{align}
	where $\tilde{C}_1$ depends only  on $C_0,~\tilde{C}_0$.
	
\end{theorem}
\begin{remark}
	This theorem  implies that $\|(\r^\v,u^\v)\|_{L^\infty(0,T;H^2)}$ is uniform bounded which yields immediately that   the boundary layers for $(\r^\v,u^\v)$  is very weak for the flat case.
\end{remark}


\

Based on the uniform estimates of Theorem \ref{thm1.1}, using similar  arguments as \cite{Masmoudi-R}, one can prove the vanishing viscosity limit of viscous solutions to the inviscid one in $L^\infty$-norm by the strong compactness argument, but without convergence rate. However,  we are interested in the  vanishing viscosity limit with rate of convergence.  In Theorem \ref{thm1.2} below, we prove the vanishing viscosity limit with rates of convergence, which generalizes the corresponding results for the incompressible case in \cite{Xiao-Xin-1, Xiao-Xin-2}.

We supplement the compressible Euler equations \eqref{1.7} and the compressible Navier-Stokes system \eqref{1.1}  with the same  initial data $(\r_0,u_0)$  satisfying
\begin{equation}\label{EI-1}
(p_0,u_0)\in H^3\cap X^{\v,m}_{NS} ~~\mbox{with}~~m\geq6.
\end{equation}
It is well known that there exists a unique smooth solution $(\r,u)(t)\in H^3$ for the  problem
\eqref{1.7}, \eqref{1.8} with initial data  $(\r_0,u_0)$ at least locally in time $[0,T_1]$ where $T_1>0$  depends only on $\|(p_0,u_0)\|_{H^3}$.
On the other hand, it follows from  Theorem \ref{thm1.1} that there exists a time $T_0>0$ and $\tilde{C}_1>0$ independent of $\v\in(0,1]$, such that there exists a unique solution $(\r^\v,u^\v)(t)$ of \eqref{1.1},\eqref{1.10} with initial data $(\r_0,u_0)$ and satisfies $\|(p(\r^\v),u^\v)(t)\|_{X_{m}^{\v}}\leq \tilde{C}_1$.

\

We justify the vanishing viscosity limit as follows:
\begin{theorem}[Inviscid Limit]\label{thm1.2}  Let $(\r,u)(t)\in L^\infty(0,T_1;H^3)$ be the smooth solution to Euler equations \eqref{1.7}, \eqref{1.8} with initial data $(\r_0,u_0)$ satisfying \eqref{EI-1}.\\[1mm]
{\bf Part I(General case):} Let  $(\r^\v,u^\v)(t)$ be the solution to the initial boundary value problem of the compressible Navier-Stokes equations \eqref{1.1},\eqref{1.10}  with  initial data $(\r_0,u_0)$ satisfying \eqref{EI-1}. Then, there exists $T_2=\min\{T_0,T_1\}>0$, which is independent of $\v>0$, such that
\begin{eqnarray}
&&\|(\r^\v-\r,u^\v-u)(t)\|^2_{L^2}+\v\int_0^t\|(u^\v-u)(\tau)\|^2_{H^1}d\tau\leq C\v^{\f32},~t\in[0,T_0],\label{1.8-0}\\
&&\|(\r^\v-\r,u^\v-u)(t)\|^2_{H^1}+\v\int_0^t\|(u^\v-u)(\tau)\|^2_{H^2}d\tau\leq C\v^{\f12},~t\in[0,T_0],\label{1.9-0}	\end{eqnarray}
and
\begin{equation}
\|(\r^\v-\r,u^\v-u)\|_{L^\infty(\Omega\times[0,T_0])}\leq \|(\r^\v-\r,u^\v-u)\|^{\f25}_{L^2}\cdot \|(\r^\v-\r,u^\v-u)\|^{\f35}_{W^{1,\infty}}\leq C\v^{\f3{10}},\label{1.10-2}
\end{equation}
where $C$ depend only on the norm $\|(\r_0,u_0)\|_{H^3}+\|(p(\r_0),u_0)\|_{X^{\v}_{m}}$.\\[2mm]
{\bf Part II(Flat case):} Let $\Om=\mathbb{T}^2\times(0,1)$ and $(\r^\v,u^\v)(t)$ be the solution to the initial boundary value problem of the compressible Navier-Stokes equations \eqref{1.1},\eqref{10.1}  with  initial data $(\r_0,u_0)$ satisfying \eqref{EI-1}. Then, there exists $T_2=\min\{T_0,T_1\}>0$, which is independent of $\v>0$, such that
\begin{eqnarray}
&&\|(\r^\v-\r,u^\v-u)(t)\|^2_{L^2}+\v\int_0^t\|(u^\v-u)(\tau)\|^2_{H^1}d\tau\leq C\v^2,~t\in[0,T_0],\label{1.8-1}\\
&&\|(\r^\v-\r,u^\v-u)(t)\|^2_{H^1}+\v\int_0^t\|(u^\v-u)(\tau)\|^2_{H^2}d\tau\leq C\v^{\f32},~t\in[0,T_0],\label{1.9-1}	\end{eqnarray}
and
\begin{equation}
\|(\r^\v-\r,u^\v-u)\|_{L^\infty(\Omega\times[0,T_0])}\leq \|(\r^\v-\r,u^\v-u)\|^{\f25}_{L^2}\cdot \|(\r^\v-\r,u^\v-u)\|^{\f35}_{W^{1,\infty}}\leq C\v^{\f2{5}},\label{1.10-3}
\end{equation}
where $C$ depend only on the norm $\|(\r_0,u_0)\|_{H^3}+\|(p(\r_0),u_0)\|_{X^{\v}_{m}}$.  Moreover, the solution $(\r,u)$ of the Euler system sasifies the  additional boundary condition, i.e.
\begin{align}\label{16.1-1}
n\times\omega=0,~~\mbox{on}~\Gamma.
\end{align}
	
\end{theorem}



\begin{remark}
In general, it is hard to obtain uniform bound for $\|u^\v\|_{L^\i(0,T;H^2)}$, otherwise, the corresponding Euler solution will satisfy \eqref{16.1-1} as above. However, usually, it is impossible for   the solution of Euler system to satisfy the additional boundary condition \eqref{16.1-1} because the boundary condition \eqref{1.8} is enough for the well-posedness of Euler system \eqref{1.7}.
\end{remark}

\begin{remark}
	The multi-scale analysis implies that the convergence should be of order $\v^{\f12}$ in $L^\infty(\Omega\times[0,T])$, so  the justification of this rate is still an difficult problem.	
\end{remark}

The rest of the paper is organized as follows: In the next section, we collect some  inequalities that will be used later.  In section 3,  we prove  the a priori estimates Theorem \ref{thm3.1}. By using the a priori estimates, we  prove Theorem \ref{thm1.1} in section 4. By careful boundary analysis,  Theorem \ref{thm1.3} is proved in  section 5. Based on the uniform estimate in Theorem \ref{thm1.1}, Theorem \ref{thm1.2} is proved in section 6.  In the Appendix, we generalize the Lemma 14 and Lemma 15 of \cite{Masmoudi-R} so that it can be applied to the case of compressible Navier-Stokes equations.


\section{Preliminaries }
The following lemma \cite{Xiao-Xin-1,Teman} allows one to control the $H^m(\Om)$-norm of a vector valued function $u$ by its $H^{m-1}(\Om)$-norm of $\nabla\times u$ and $\mbox{div} u$, together with the $H^{m-\f12}(\partial\Om)$-norm of $u\cdot n$.

\begin{proposition}\label{prop3.1}
Let $m\in \mathbb{N}_+$ be an integer. Let $u\in H^m$ be a vector-valued function. Then, there exists a constant $C>0$ in dependent $u$, such that
\begin{equation}\label{3.1}
\|u\|_{H^m}\leq C\left(\|\nabla \times u\|_{H^{m-1}}+\|\mbox{div}u\|_{H^{m-1}}
+\|u\|_{H^{m-1}}+|u\cdot n|_{H^{m-\f12}(\partial\Om)}\right).
\end{equation}
and
\begin{equation}\label{3.1-1}
\|u\|_{H^m}\leq C\left(\|\nabla \times u\|_{H^{m-1}}+\|\mbox{div}u\|_{H^{m-1}}
+\|u\|_{H^{m-1}}+|n\times u|_{H^{m-\f12}(\partial\Om)}\right).
\end{equation}

\end{proposition}

\

In this paper, we shall use repeatedly the  Gagliardo-Nirenbirg-Morser type inequality, whose proof can be find in \cite{Gues}. First, define the space
\begin{equation}\label{2.2}
\mathcal{W}^m(\Om\times[0,T])=\{f(x,t)\in L^2(\Omega\times[0,T]) ~|~ \mathcal{Z}^\a f \in L^2(\Omega\times[0,T]) , ~~|\a|\leq m~\}.
\end{equation}
Then, the Gagliardo-Nirenbirg-Morser type inequality is as follows:
\begin{proposition}\label{prop3.2}
For $ u,v\in L^\infty(\Om\times[0,T])\cap \mathcal{W}^m(\Om\times[0,T])$ with  $m\in \mathbb{N}_+$ be an integer. It holds that
\begin{equation}\label{3.2}
\int_0^t\|(\mathcal{Z}^{\b}u\mathcal{Z}^{\g}v)(\tau)\|^2d\tau
\lesssim \|u\|^2_{L^\infty_{t,x}}\int_0^t\|v(\tau)\|^2_{\mathcal{H}^m}d\tau
+\|v\|^2_{L^\infty_{t,x}}\int_0^t\|u(\tau)\|^2_{\mathcal{H}^m}d\tau,~~|\b|+|\g|=m.
\end{equation}

\end{proposition}

\

We also need the following anisotropic Sobolev embedding and trace estimates:
\begin{proposition}\label{prop3.3}
Let  $m_1\geq 0,~m_2\geq 0$ be integers, $f\in H^{m_1}_{co}(\Om)\cap  H^{m_2}_{co}(\Om)$  and   $\nabla f\in H^{m_2}_{co}(\Om)$.\\
1) The following anisotropic Sobolev embedding holds:
\begin{equation}\label{3.3}
\|f\|^2_{L^\infty}\leq C \Big(\|\nabla f\|_{H^{m_2}_{co}}+ \|f\|_{H^{m_2}_{co}}\Big)\cdot\|f\|_{H^{m_1}_{co}},
\end{equation}
provided $m_1+m_2\geq 3$.\\[2mm]
2) The following trace estimate holds:
\begin{equation}\label{3.4}
|f|^2_{H^{s}(\partial\Om)}\leq C\Big(\|\nabla f\|_{H^{m_2}_{co}}+ \|f\|_{H^{m_2}_{co}}\Big)\cdot\|f\|_{H^{m_1}_{co}}.
\end{equation}
provided $m_1+m_2\geq 2s\geq 0$.
\end{proposition}

\noindent\textbf{Proof}. The proof is just a using of  the covering $\Om\subset \Om_0\cup_{k=1}^n\Om_k$ and Proposition 2.2 in \cite{Masmoudi-R-1}, the details are tus omitted here.  $\hfill\Box$




\section{A priori Estimates}

The aim of this section is to prove the following {\it a priori} estimates, which is a crucial step to prove Theorem \ref{thm1.1}. For notational convenience, we drop the superscript $\v$ throughout this section.

\begin{theorem}[A priori Estimates]\label{thm3.1}
Let $m$ be an integer satisfying $m\geq 6$,  $\Omega$ be a $\mathcal{C}^{m+2}$ domain and $A\in\mathcal{C}^{m+1}(\partial\Omega)$. For very sufficiently smooth solution defined on $[0,T]$ of \eqref{1.1} and  \eqref{1.10}(or  \eqref{2.7}), then it holds that
\begin{equation}\label{3.0-3}
|\r(x,0)|\exp(-\int_0^t\|\mbox{div}u(\tau)\|_{L^\infty}d\tau)\leq\r(x,t)\leq |\r(x,0)|\exp(\int_0^t\|\mbox{div}u(\tau)\|_{L^\infty}d\tau), \forall t\in[0,T].
\end{equation}
In addition, if
\begin{equation}\label{3.0-4}
0<c_0\leq \r(t)\leq \f1{c_0}<\infty,~~\forall t\in[0,T],
\end{equation}
where $c_0$ is any given small positive constant, then  the following a priori estimate holds
\begin{eqnarray}\label{3.0-1}
&&\mathcal{N}_m(t)+\int_{0}^{t}\|\nabla\partial_t^{m-1}p(\tau)\|^2
+\|\Delta p(\tau)\|^2_{\mathcal{H}^2}d\tau
+\v\int_{0}^{t}\|\nabla u(\tau)\|^2_{\mathcal{H}^{m}}d\tau\nonumber\\
&&~~~~~~~~+\v\sum_{k=0}^{m-2}\int_{0}^{t}\|\nabla^2\partial_t^ku(\tau)\|^2_{m-k-1}d\tau
+\v^2\int_{0}^{t}\|\nabla^2\partial_t^{m-1}u(\tau)\|^2d\tau\nonumber\\
&&\leq \tilde{C}_2 C_{m+2}\Big\{P(\mathcal{N}_m(0))+P(\mathcal{N}_m(t))
\cdot\int_0^tP(\mathcal{N}_m(\tau))d\tau\Big\}, ~~\forall t\in[0,T],
\end{eqnarray}
where $\tilde{C}_2$ depends only on $\f1{c_0}$, $P(\cdot)$ is a polynomial and
\begin{eqnarray}\label{3.0-2}
&&\mathcal{N}_m(t)\triangleq\mathcal{N}_m(p,u)(t)=\sup_{0\leq\tau\leq t}\Big\{1+\|(p,u)(\tau)\|^2_{\mathcal{H}^m}+\|\nabla u(\tau)\|^2_{\mathcal{H}^{m-1}}+\sum_{k=0}^{m-2}\|\partial_t^k \nabla p(\tau)\|^2_{m-1-k}\nonumber\\
&&~~~~~~~~~~~~~~~~~~~~~~~~~~+\|\Delta p(t)\|^2_{\mathcal{H}^1}
+\|\nabla u(\tau)\|^2_{\mathcal{H}^{1,\infty}}+\v\|\nabla\partial_t^{m-1}p(\tau)\|^2
+\v\|\Delta p(\tau)\|^2_{\mathcal{H}^2}\Big\}.
\end{eqnarray}

\end{theorem}

Throughout this section, we shall work on the interval of time $[0,T]$ such that $c_0\leq \r(t)\leq \f1{c_0}$. And we point out that the generic constant $C$ may depend on $\f1{c_0}$ in this section.  Since the proof of Theorem \ref{thm3.1} is quite lengthy and involved, we divide the proof into the following several subsections.

\subsection{Conormal Energy Estimates}

\renewcommand{\theequation}{\arabic{section}.\arabic{subsection}.\arabic{equation}}

Notice that
\begin{equation}\label{2.4}
\Delta u=\nabla\mbox{div} u-\nabla\times\nabla\times u,
\end{equation}
then $\eqref{1.1}_2$ is rewritten as
\begin{equation}\label{2.5}
\r u_t+\r u\cdot\nabla u+\nabla p=-\mu\v \nabla\times\omega+(2\mu+\l)\v \nabla div u,
\end{equation}
where $\omega=\nabla \times u$ is the vorticity.  Since $\mu>0, 2\mu+\l>0$, we  normalize $\mu$ and $2\m+\l$ to be 1 and 2 respectively for simplicity.

\

In this subsection, we first give the basic a priori $L^2$ energy estimate which holds for  \eqref{1.1} with \eqref{1.10}.
\begin{lemma}\label{lem2.1}
For a smooth solution to \eqref{1.1} and \eqref{2.7}, it holds that for $\v\in(0,1]$
\begin{equation}\label{2.8}
\sup_{0\leq\tau\leq t}\Big(\int \f12\r|u|^2+\f{1}{\g-1}\r^\g dx\Big)+c_1\v \int_0^t\|\nabla u\|^2d\tau
\leq \int \f12\r_0|u_0|^2+\f{1}{\g-1}\r_0^\g dx+ C \int_0^t\|u\|^2d\tau,
\end{equation}
where $c_1>0$ is a positive constant.
\end{lemma}

\noindent\textbf{Proof}. Multiplying \eqref{2.5} by $u$,  using the boundary condition and integrating by parts, we have that
\begin{equation}\label{2.9}
\f{d}{dt}\int \f12\r|u|^2dx+\int \nabla p udx =-\v \int\nabla\times\omega udx+2\v \int\nabla\mbox{div}u udx.
\end{equation}
By using $\eqref{1.1}_1$, we obtain that
\begin{equation}\label{2.10}
\int \nabla p udx=\frac{\g}{\g-1}\int\nabla\r^{\g-1}\cdot \r udx=\frac{\g}{\g-1}\int\r^{\g-1}\r_tdx
=\f{d}{dt}\int\frac{1}{\g-1}\r^\g dx.
\end{equation}
Integrating by parts and using the boundary conditions \eqref{2.7}, one has that
\begin{eqnarray*}\label{2.11}
&&-\v \int\nabla\times\omega udx=-\v \|\omega\|^2-\v \int_{\partial\Om}(n\times\omega)\cdot u d\sigma
 \leq-\v \|\omega\|^2+C\v|u|^2_{L^2{(\partial\Om)}},
\end{eqnarray*}
and
\begin{eqnarray}\label{2.12}
&&\v \int u\nabla\mbox{div}u dx=-\v \|\mbox{div}u\|^2.
\end{eqnarray}
Substituting \eqref{2.10}-\eqref{2.12} into \eqref{2.9} leads to that
\begin{equation}\label{2.13}
\f{d}{dt}\Big(\int \f12\r|u|^2+\f{1}{\g-1}\r^\g dx\Big)+\v \|\omega\|^2+2\v \|\mbox{div}u\|^2
\leq C\v|u|^2_{L^2{(\partial\Om)}}.
\end{equation}
Due to Proposition \ref{prop3.1}, it holds that
\begin{equation*}\label{2.14}
 \|\omega\|^2+\|\mbox{div}u\|^2
\geq 2c_1\|\nabla u\|^2-C\|u\|^2.
\end{equation*}
The trace theorem implies
\begin{equation}\label{2.15}
\v |u|^2_{L^2{(\partial\Om)}}\leq  \v \|u\|_{H^1}\cdot\|u\|\leq \f12c_1\v \|\nabla u\|^2+C\v\|u\|^2.
\end{equation}
Then \eqref{2.13}-\eqref{2.15} give that
\begin{equation}\label{2.13}
\f{d}{dt}\Big(\int \f12\r|u|^2+\f{1}{\g-1}\r^\g dx\Big)+c_1\v \|\nabla u\|^2
\leq C\|u\|^2.
\end{equation}
Integrating the above inequality with respect to $t$ yields  \eqref{2.8} immediately. Thus   the proof of the Lemma \ref{lem2.1} is completed.
$\hfill\Box$

\

Set
\begin{equation}\label{2.15-1}
\Lambda_m(t) \triangleq \|(p,u)(t)\|^2_{\mathcal{H}^{m}}+\|\nabla u(t)\|^2_{\mathcal{H}^{m-1}}
+\sum_{k=0}^{m-2}\|\nabla\partial_t^{k}p(t)\|^2_{m-1-k}+\v\|\nabla\partial_t^{m-1}p(t)\|^2,
\end{equation}
and
\begin{equation}\label{2.15-2}
Q(t)\triangleq \sup_{0\leq\tau\leq t}\Big\{\|(\nabla p, \nabla u)(t)\|^2_{\mathcal{H}^{1,\infty}}
+\|(p,u, p_t,u_t)(t)\|^2_{L_x^\infty}\Big\}.
\end{equation}

\begin{lemma}\label{lem2.2}
For every $m\in \mathbb{N}_+$, it holds that
\begin{eqnarray}\label{2.16}
&&\sup_{0\leq \tau\leq t}\|(u, p)\|^2_{\mathcal{H}^m}+\v\int_0^t\|\nabla u(\tau)\|^2_{\mathcal{H}^m}d\tau\leq C
C_{m+2} \Big\{\|(u_0,p_0)\|^2_{\mathcal{H}^m}+\d \v^2\int_0^t\|\nabla^2 u(\tau)\|^2_{\mathcal{H}^{m-1}}d\tau\nonumber\\
&&~~~~~~~+\d\int_0^t\|\nabla\partial^{m-1}_t p(\tau)\|^2d\tau+C_\d[1+P(Q(t))]\int_0^t\Lambda_m(\tau) d\tau\Big\}.
\end{eqnarray}

\end{lemma}

\noindent\textbf{Proof}.
The case for $m=0$ is already proved in Lemma \ref{lem2.1}. Assume that \eqref{2.16} is proved for $k\leq m-1$. We shall prove that it holds for $ k=m\geq 1$. By applying $\mathcal{Z}^\a$ with $|\a|=m$ to \eqref{2.5}, we obtain
\begin{equation}\label{2.17}
\r\mathcal{Z}^\a u_t+\r u\cdot\nabla\mathcal{Z}^\a u+\mathcal{Z}^\a\nabla p=-\v\mathcal{Z}^\a\nabla\times\omega+2\v\mathcal{Z}^\a\nabla\mbox{div}u+\mathcal{C}_1^\a+\mathcal{C}_2^\a,
\end{equation}
where
\begin{equation}\label{2.18}
\mathcal{C}_1^\a=-[\mathcal{Z}^\a, \r]u_t=-\sum_{|\b|\geq 1,\b+\g=\a}C_{\b,\g}\mathcal{Z}^\b\r\mathcal{Z}^\g u_t,
\end{equation}
and
\begin{equation}\label{2.19}
\mathcal{C}_2^\a=-[\mathcal{Z}^\a, \r u\cdot\nabla]u=-\sum_{|\b|\geq 1,\b+\g=\a}C_{\b,\g}\mathcal{Z}^\b(\r u)\mathcal{Z}^\g\nabla u-\r u\cdot [\mathcal{Z}^\a,\nabla]u.
\end{equation}
Multiplying \eqref{2.17} by $\mathcal{Z}^\a u$, and  integrating by parts, one gets that
\begin{eqnarray}\label{2.20}
&&\f{d}{dt}\int\f12\r|\mathcal{Z}^\a u|^2dx+\int\mathcal{Z}^\a\nabla p\mathcal{Z}^\a udx\nonumber\\
&&=-\v\int\mathcal{Z}^\a\nabla\times\omega\cdot\mathcal{Z}^\a udx+2\v\int\mathcal{Z}^\a\nabla\mbox{div}u\cdot\mathcal{Z}^\a udx
+\int(\mathcal{C}_1^\a+\mathcal{C}_2^\a)\mathcal{Z}^\a udx.
\end{eqnarray}
Notice that
\begin{eqnarray}\label{2.21}
&&-\v\int\mathcal{Z}^\a\nabla\times\omega\cdot\mathcal{Z}^\a udx=-\v\int\nabla\times\mathcal{Z}^\a\omega\cdot\mathcal{Z}^\a udx-\v\int[\mathcal{Z}^\a, \nabla\times]\omega\cdot\mathcal{Z}^\a udx\nonumber\\
&&\leq -\v\int\mathcal{Z}^\a\omega\cdot\nabla\times\mathcal{Z}^\a udx-\v\int_{\partial\Om} n\times\mathcal{Z}^\a\omega\cdot\mathcal{Z}^\a ud\sigma\nonumber\\
&&~~~~~~~~~~~~~~~~~~~~~~~~~+C\v\|\nabla^2u\|_{\mathcal{H}^{m-1}}\|u\|_{\mathcal{H}^{m}}
+C(\|\nabla u\|^2_{\mathcal{H}^{m-1}}+\|u\|^2_{\mathcal{H}^{m}})\nonumber\\
&&\leq -\f{3\v}{4}\|\nabla\times\mathcal{Z}^\a u\|^2-\v\int_{\partial\Om} n\times\mathcal{Z}^\a\omega\cdot\mathcal{Z}^\a ud\sigma\nonumber\\
&&~~~~~~~~~~~~~~~~~~~~~~~~+\d\v^2\|\nabla^2u\|^2_{\mathcal{H}^{m-1}}+C_\d(\|\nabla u\|^2_{\mathcal{H}^{m-1}}+\|u\|^2_{\mathcal{H}^{m}}).
\end{eqnarray}
In order to complete the estimates of \eqref{2.21}, one needs to estimate the boundary terms involving  $\MZ^\a u$ with $\a_{13}=0$(For $\a_{13}\neq 0$,   $\MZ^\a u|_{\partial\Om}=0$ by definition). Due to \eqref{2.7}, one has for $|\a_0|+|\a_1|=m$
\begin{eqnarray}\label{2.22}
&&|n\times\mathcal{Z}^\a\omega|_{L^2(\partial\Om)}\leq C_{m+2} \left(|\partial_t^{\a_0}\omega|_{H^{|\a_1|-1}(\partial\Om)}
+|\partial_t^{\a_0}u|_{H^{|\a_1|}(\partial\Om)}\right)\nonumber\\
&&\leq C_{m+2} \|\nabla\partial_t^{\a_0}\omega\|^{\f12}_{|\a_1|-1}\cdot\|\partial_t^{\a_0}\omega\|^{\f12}_{|\a_1|-1}
+C_{m+2} \|\nabla\partial_t^{\a_0}u\|^{\f12}_{|\a_1|}\cdot\|\partial_t^{\a_0}u\|^{\f12}_{|\a_1|}\nonumber\\
&&~~~~~~~~~~~~~~~~~~~~~+C_{m+2}(\|\nabla u\|_{\mathcal{H}^{m-1}}+\|u\|_{\mathcal{H}^{m}})\nonumber\\
&&\leq C_{m+2} \Big(\|\nabla^2u\|^{\f12}_{\mathcal{H}^{m-1}}\cdot\|\nabla u\|^{\f12}_{\mathcal{H}^{m-1}}
+\|\nabla u\|^{\f12}_{\mathcal{H}^{m}}\cdot \|u\|^{\f12}_{\mathcal{H}^{m}}+\|\nabla u\|_{\mathcal{H}^{m-1}}+\|u\|_{\mathcal{H}^{m}}\Big),
\end{eqnarray}
thus
\begin{eqnarray}\label{2.23}
&&\v\left|\int_{\partial\Om} n\times\mathcal{Z}^\a\omega\cdot\mathcal{Z}^\a ud\sigma\right|
\leq \v|n\times\mathcal{Z}^\a\omega|_{L^2(\partial\Om)}\cdot |\mathcal{Z}^\a u|_{L^2(\partial\Om)}\nonumber\\
&&\leq C_{m+2}\v(\|\nabla u\|^{\f12}_{\mathcal{H}^{m}}+\|u\|^{\f12}_{\mathcal{H}^{m}})\cdot \|u\|^{\f12}_{\mathcal{H}^{m}}
\left(\|\nabla^2u\|_{\mathcal{H}^{m-1}}
+\|\nabla u\|_{\mathcal{H}^{m}}+\|\nabla u\|_{\mathcal{H}^{m-1}}+ \|u\|_{\mathcal{H}^{m}} \right)\nonumber\\
&&\leq \d\v \|\nabla u\|^2_{\mathcal{H}^{m}}+\d\v^2 \|\nabla^2 u\|^2_{\mathcal{H}^{m-1}}
+C_\d C_{m+2} \left( \|u\|^2_{\mathcal{H}^{m}}+ \|\nabla u\|^2_{\mathcal{H}^{m-1}} \right).
\end{eqnarray}
This, together with   \eqref{2.21},   yields that
\begin{eqnarray}\label{2.24}
&&-\v\int\mathcal{Z}^\a\nabla\times\omega\cdot\mathcal{Z}^\a udx
\leq -\f{3\v}{4}\|\nabla\times\mathcal{Z}^\a u\|^2+\d\v \|\nabla u\|^2_{\mathcal{H}^{m}}+\d\v^2\|\nabla^2u\|^2_{\mathcal{H}^{m-1}}\nonumber\\
&&~~~~~~~~~~~~~~~~~~~~~~~~~~~~~~~~~~~~~+C_\d C_{m+2} (\|\nabla u\|^2_{\mathcal{H}^{m-1}}+\|u\|^2_{\mathcal{H}^{m}}).
\end{eqnarray}
Notice that
\begin{eqnarray}\label{2.25}
&&\v\int\mathcal{Z}^\a\nabla\mbox{div}u\cdot\mathcal{Z}^\a udx=\v\int\nabla\mathcal{Z}^\a\mbox{div}u\cdot\mathcal{Z}^\a udx+\v\int[\mathcal{Z}^\a, \nabla]\mbox{div}u\cdot\mathcal{Z}^\a udx\nonumber\\
&&\leq-\v\int\mathcal{Z}^\a\mbox{div}u\cdot\mbox{div}\mathcal{Z}^\a udx+\v\int_{\partial\Om} \mathcal{Z}^\a\mbox{div}u\cdot(\mathcal{Z}^\a u\cdot n)d\sigma\nonumber\\
&&~~~~~~~~~~~~~~~~~~~~~~~~~~+\v\|\nabla^2u\|_{\mathcal{H}^{m-1}}\cdot\|u\|_{\mathcal{H}^{m}}+C(\|\nabla u\|^2_{\mathcal{H}^{m-1}}+\|u\|^2_{\mathcal{H}^{m}})\nonumber\\
&&\leq -\f{3\v}{4}\|\mbox{div}\mathcal{Z}^\a u\|^2+\v\int_{\partial\Om} \mathcal{Z}^\a\mbox{div}u\cdot(\mathcal{Z}^\a u\cdot n)d\sigma\nonumber\\
&&~~~~~~~~~~~~~~~~~~~~~~~~~~+\d\v^2\|\nabla^2u\|^2_{\mathcal{H}^{m-1}}+C_\d(\|\nabla u\|^2_{\mathcal{H}^{m-1}}+\|u\|^2_{\mathcal{H}^{m}}).
\end{eqnarray}
In order to estimate the boundary term in  the above term, one needs to discuss the following two cases. If $|\a_0|=|\a|$, then by \eqref{2.7}, one has that  $\partial_t^{\a_0}u\cdot n|_{\partial\Om}=0$ which implies that
\begin{equation}\label{2.26}
\v\int_{\partial\Om} \mathcal{Z}^\a\mbox{div}u\cdot(\mathcal{Z}^\a u\cdot n)d\sigma=0.
\end{equation}
If $|\a_1|\geq 1$, then by using \eqref{2.7} and integrating by parts along the boundary, one obtains  that
\begin{eqnarray}\label{2.27}
&&\v\int_{\partial\Om} \mathcal{Z}^\a\mbox{div}u\cdot(\mathcal{Z}^\a u\cdot n)d\sigma
=\v\int_{\partial\Om} Z^{\a_1}\partial_t^{\a_0}\mbox{div}u\cdot(Z^{\a_1}\partial_t^{\a_0} u\cdot n)d\sigma\nonumber\\
&&=-\v\int_{\partial\Om} Z^{\a_1-1}\partial_t^{\a_0}\mbox{div}u\cdot Z_y(Z^{\a_1}\partial_t^{\a_0} u\cdot n)d\sigma\leq \v|Z^{\a_1-1}\partial_t^{\a_0}\mbox{div}u|_{L^2} |Z_y(Z^{\a_1}\partial_t^{\a_0} u\cdot n)|_{L^2}\nonumber\\
&&\leq C_{m+2} \v|Z^{\a_1-1}\partial_t^{\a_0}\mbox{div}u|_{L^2}\cdot |\partial_t^{\a_0} u|_{H^{|\a|-|\a_0|}}\nonumber\\
&&\leq C_{m+2} \v(\|\nabla^2u\|^{\f12}_{\mathcal{H}^{m-1}}+\|\nabla u\|^{\f12}_{\mathcal{H}^{m-1}})\cdot\|\nabla u\|^{\f12}_{\mathcal{H}^{m-1}}\cdot(
\|\nabla u\|^{\f12}_{\mathcal{H}^{m}}+\|u\|^{\f12}_{\mathcal{H}^{m}})\|u\|^{\f12}_{\mathcal{H}^{m}}
\nonumber\\
&&\leq \d\v \|\nabla u\|^2_{\mathcal{H}^{m}}+\d\v^2\|\nabla^2u\|^2_{\mathcal{H}^{m-1}}+C_\d C_{m+2} (\|\nabla u\|^2_{\mathcal{H}^{m-1}}+\|u\|^2_{\mathcal{H}^{m}}),
\end{eqnarray}
where $Z_y$ or $\partial_y$ represents the derivatives involves only the tangential parts.
Then \eqref{2.25}-\eqref{2.27} yield that
\begin{eqnarray}\label{2.28}
&&\v\int\mathcal{Z}^\a\nabla\mbox{div}u\cdot\mathcal{Z}^\a udx
\leq -\f{3\v}{4}\|\mbox{div}\mathcal{Z}^\a u\|^2
+\d\v\|\nabla u\|^2_{\mathcal{H}^{m}}+2\d\v^2\|\nabla^2u\|^2_{\mathcal{H}^{m-1}}\nonumber\\
&&~~~~~~~~~~~~~~~~~~~~~~~~~~~~~~~~~~+C_\d C_{m+2} (\|\nabla u\|^2_{\mathcal{H}^{m-1}}
+\|u\|^2_{\mathcal{H}^{m}}).
\end{eqnarray}
On the other hand, it follows from Proposition \ref{prop3.1},  that
\begin{eqnarray}\label{2.29}
&&2c_1\|\nabla\mathcal{Z}^\a u\|^2_{L^2}\leq  \Big(\|\nabla\times\mathcal{Z}^\a u\|^2_{L^2}
+\|\mbox{div}\mathcal{Z}^\a u\|^2_{L^2}+\|\mathcal{Z}^\a u\|^2_{L^2}+|\mathcal{Z}^\a u\cdot n|^2_{H^{\f12}(\partial\Om)}\Big)\nonumber\\
&&~~~~~~~~~~~~~~\leq \Big(\|\nabla\times\mathcal{Z}^\a u\|^2_{L^2}
+\|\mbox{div}\mathcal{Z}^\a u\|^2_{L^2}\Big)+C_{m+2} \Big(\|u\|^2_{\mathcal{H}^{m}}+\|\nabla u\|^2_{\mathcal{H}^{m-1}}\Big)
\end{eqnarray}
where one has used
\begin{eqnarray}\label{4.11}
|Z_y^{m-k}\partial_t^k u\cdot n|_{H^{\f12}}\leq
\begin{cases}
0, ~\mbox{if}~~k=m,\\[3mm]
C_{m+2} |\partial_t^{k}u|_{H^{m-k-\f12}}
\leq \sum_{|\a|\leq m-1}|\MZ^\a u|_{H^\f12}\\
\leq C_{m+2}(\|\nabla u\|_{\mathcal{H}^{m-1}}+\|u\|_{\mathcal{H}^{m}}), ~\mbox{if}~~k\leq m-1,
\end{cases}
\end{eqnarray}
which is a consequence of \eqref{2.7} and \eqref{3.4}.

Integrating the resulting equation \eqref{2.20}, and  substituting  \eqref{2.29}, \eqref{2.28} and \eqref{2.24} into \eqref{2.20}, one gets that
\begin{eqnarray}\label{2.30}
&&\f12\int\r|\mathcal{Z}^\a u|^2dx+\int_0^t\int\mathcal{Z}^\a\nabla p\cdot\mathcal{Z}^\a udxd\tau
+2c_1\v\int_0^t\|\nabla\mathcal{Z}^\a u(\tau)\|^2_{L^2}d\tau\nonumber\\
&&\leq \f12\int\r_0|\mathcal{Z}^\a u_0|^2dx+C\d\v^2\int_0^t\|\nabla^2u(\tau)\|^2_{\mathcal{H}^{m-1}}d\tau
+C\d\v\int_0^t\|\nabla u(\tau)\|^2_{\mathcal{H}^{m}}d\tau\nonumber\\
&&~~~~+C_\d C_{m+2} \int_0^t\|\nabla u(\tau)\|^2_{\mathcal{H}^{m-1}}+\|u(\tau)\|^2_{\mathcal{H}^{m}}d\tau
+\int_0^t\int(\mathcal{C}_1^\a+\mathcal{C}_2^\a)\cdot\mathcal{Z}^\a udxd\tau.
\end{eqnarray}
Now we estimate the pressure term on the left hand side of \eqref{2.30}. Notice that
\begin{eqnarray}\label{2.31}
&&\int_0^t\int\mathcal{Z}^\a\nabla p\cdot\mathcal{Z}^\a udxd\tau
=\int_0^t\int\nabla\mathcal{Z}^\a p\cdot\mathcal{Z}^\a udxd\tau
+\int_0^t\int[\mathcal{Z}^\a, \nabla] p\cdot\mathcal{Z}^\a udxd\tau\nonumber\\
&&\geq -\int_0^t\int\mathcal{Z}^\a p\cdot\mbox{div}\mathcal{Z}^\a udxd\tau+\int_0^t\int_{\partial\Om}\mathcal{Z}^\a p\mathcal{Z}^\a u\cdot n d\sigma d\tau
-C\int_0^t\|u\|_{\mathcal{H}^{m}}\|\nabla p\|_{\mathcal{H}^{m-1}}d\tau\nonumber\\
&&\geq -\int_0^t\int\mathcal{Z}^\a p\cdot\mathcal{Z}^\a \mbox{div}udxd\tau+\int_0^t\int_{\partial\Om}\mathcal{Z}^\a p\mathcal{Z}^\a u\cdot n d\sigma d\tau
-\d\int_0^t\|\nabla p\|^2_{\mathcal{H}^{m-1}}d\tau\nonumber\\
&&~~~~~~~~~~~~~~~-C_\d\int_0^t\|(p,u)\|^2_{\mathcal{H}^{m}}+\|\nabla u\|^2_{\mathcal{H}^{m-1}}d\tau
\end{eqnarray}
First, we treat the  boundary term when $\a_{13}=0$(for $\a_{13}\neq0$,  $\mathcal{Z}^\a u=0$ on the boundary) in the right hand side of \eqref{2.31}. If $|\a_0|=|\a|$,   one has, from \eqref{4.11}, that
\begin{equation}\label{2.32}
\int_0^t\int_{\partial\Om}\mathcal{Z}^\a p\mathcal{Z}^\a u\cdot n d\sigma d\tau=0.
\end{equation}
If $|\a_1|\geq 1$, integrating by parts along the boundary and using \eqref{4.11}, one has that
\begin{eqnarray}\label{2.33}
&&|\int_0^t\int_{\partial\Om}\mathcal{Z}^\a p\mathcal{Z}^\a u\cdot n d\sigma d\tau|
=|\int_0^t\int_{\partial\Om}Z_y^{\a_1}\partial_t^{\a_0} p Z_y^{\a_1}\partial_t^{\a_0}u\cdot n d\sigma d\tau|\nonumber\\
&&=\int_0^t|Z_y^{\a_1-1}\partial_t^{\a_0} p|_{H^{\f12}} \cdot |Z_y^{\a_1}\partial_t^{\a_0}u\cdot n|_{H^{\f12}} d\sigma d\tau\nonumber\\
&&\leq C_{m+2} \int_0^t \left(\|\nabla Z^{\a_1-1}\partial_t^{\a_0} p\|+ \| Z^{\a_1-1}\partial_t^{\a_0} p\|\right)^{\f12}\|Z^{\a_1-1}\partial_t^{\a_0} p\|^{\f12}_{H^1_{co}}(\|\nabla u\|_{\mathcal{H}^{m-1}}+\|u\|_{\mathcal{H}^{m}})d\tau\nonumber\\
&&\leq \d\int_0^t \|\nabla p\|^2_{\mathcal{H}^{m-1}}d\tau+C_\d C_{m+2} \int_0^t(\|\nabla u\|^2_{\mathcal{H}^{m-1}}+\|(p,u)\|^2_{\mathcal{H}^{m}})d\tau.
\end{eqnarray}
Therefore, it follows  from \eqref{2.32} and \eqref{2.33},   that
\begin{equation}\label{2.34}
\left|\int_0^t\int_{\partial\Om}\mathcal{Z}^\a p\mathcal{Z}^\a u\cdot n d\sigma d\tau\right|\leq \d\int_0^t \|\nabla p\|^2_{\mathcal{H}^{m-1}}d\tau+C_\d C_{m+2}\int_0^t(\|\nabla u\|^2_{\mathcal{H}^{m-1}}+\|(p,u)\|^2_{\mathcal{H}^{m}})d\tau.
\end{equation}
In order to estimate the first term on the right hand side of \eqref{2.31}, one rewrites the  equation  $\eqref{1.1}_1$ as
\begin{equation}\label{2.35}
\mbox{div}u=-(\ln\r)_t-u\cdot\nabla\ln\r=-\f{p_t}{\g p}-\f{u}{\g p}\cdot\nabla p.
\end{equation}
Applying $\MZ^\a$ to \eqref{2.35} yields  that
\begin{eqnarray}\label{2.36}
&&\mathcal{Z}^\a\mbox{div}u=-\f{1}{\g p}\mathcal{Z}^{\a}p_t-\f{u}{\g p}\cdot\mathcal{Z}^{\a}\nabla p
-\sum_{|\beta|\geq1, \b+\g=\a}C_{\b,\g}\MZ^\b(\f1{\g p})\cdot\MZ^{\g}p_t\nonumber\\
&&~~~~~~~~~~~~~~~~~~~~-\sum_{|\beta|\geq1, \b+\g=\a}C_{\b,\g}\MZ^\b(\f{u}{\g p})\cdot\MZ^{\g}\nabla p.
\end{eqnarray}
It is easy to get that
\begin{eqnarray}\label{2.37}
&&\int_0^t\int\mathcal{Z}^\a p\cdot\f{1}{\g p}\MZ^{\a}p_tdxd\tau=\int_0^t\int(\f{1}{2\g p}|\MZ^{\a}p|^2)_tdxd\tau-\int_0^t\int(\f{1}{2\g p})_t|\MZ^{\a}p|^2dxd\tau\nonumber\\
&&\geq \int\f{1}{2\g p}|\MZ^{\a}p|^2dx-\int\f{1}{2\g p_0}|\MZ^{\a}p_0|^2dx
-C\|p_t\|_{L^\infty}\int_0^t\|\MZ^\a p\|^2d\tau.
\end{eqnarray}
It follows from integrating by parts and  \eqref{2.7}, that
\begin{eqnarray}\label{2.38}
&&\int_0^t\int\mathcal{Z}^\a p\cdot\f{u}{\g p}\cdot\mathcal{Z}^{\a}\nabla pdxd\tau\nonumber\\
&&\geq\int_0^t\int\f{u}{2\g p}\cdot\nabla(|\mathcal{Z}^{\a}p|^2)dxd\tau
-C\|\f{u}{p}\|_{L^\infty}\int_0^t(\|p\|_{\mathcal{H}^m}+\|\nabla p\|_{\mathcal{H}^{m-1}})\|p\|_{\mathcal{H}^m}d\tau\nonumber\\
&&\geq -\d\int_0^t \|\nabla p\|^2_{\mathcal{H}^{m-1}}d\tau-  C_\d P(Q(t))\int_0^t\|p\|^2_{\mathcal{H}^{m}}d\tau.
\end{eqnarray}
Due to Proposition \ref{prop3.2}, one has that
\begin{eqnarray}\label{2.39}
&&\sum_{|\beta|\geq1, \b+\g=\a}|\int_0^t\int C_{\b,\g}\mathcal{Z}^\a p\MZ^\b(\f1{ p})\cdot\MZ^{\g}p_tdxd\tau|\nonumber\\
&&\lesssim\left(\int_0^t\|p\|^2_{\mathcal{H}^{m}}d\tau\right)^{\f12}
\sum_{|\beta|\geq1, \b+\g=\a}\left(\int_0^t\|\MZ^\b(\f1{ p})\cdot\MZ^{\g}p_t\|^2d\tau\right)^{\f12}\nonumber\\
&&\lesssim \left(\int_0^t\|p\|^2_{\mathcal{H}^{m}}d\tau\right)^{\f12}
\left(\sup_{0\leq\tau\leq t}\|\MZ p\|^2_{L^{\infty}}\int_0^t\|p_t\|^2_{\mathcal{H}^{m-1}}d\tau
+\sup_{0\leq\tau\leq t}\|p_t\|^2_{L^{\infty}}\int_0^t\|\f1p\|^2_{\mathcal{H}^{m}}d\tau\right)^{\f12}\nonumber\\
&&\lesssim [1+P(Q(t))]\int_0^t\|p\|^2_{\mathcal{H}^{m}}d\tau,
\end{eqnarray}
and
\begin{eqnarray}\label{2.40}
&&\sum_{|\beta|\geq1, \b+\g=\a}|\int_0^t\int C_{\b,\g}\mathcal{Z}^\a p\MZ^\b(\f u{ p})\cdot\MZ^{\g}\nabla pdxd\tau|\nonumber\\
&&\lesssim \left(\int_0^t\|p\|^2_{\mathcal{H}^{m}}d\tau\right)^{\f12}
\sum_{|\beta|\geq1, \b+\g=\a}\left(\int_0^t\|\MZ^\b(\f u{ p})\cdot\MZ^{\g}\nabla p\|^2d\tau\right)^{\f12}\nonumber\\
&&\lesssim \left(\int_0^t\|p\|^2_{\mathcal{H}^{m}}d\tau\right)^{\f12}
\left([1+P(Q(t))]\int_0^t\|\nabla p\|^2_{\mathcal{H}^{m-1}}+\|\f up\|^2_{\mathcal{H}^{m}}d\tau
\right)^{\f12}\nonumber\\
&&\leq\d\int_0^t\|\nabla p\|^2_{\mathcal{H}^{m-1}}d\tau+ [1+P(Q(t))]\int_0^t\|(p,u)\|^2_{\mathcal{H}^{m}}d\tau,
\end{eqnarray}
where in the estimates of \eqref{2.39} and \eqref{2.40}, one has used
\begin{eqnarray}\label{2.41}
&&c_2\int_0^t\|p\|^2_{\mathcal{H}^m}d\tau-C[1+P(Q(t))]\int_0^t\|p\|^2_{\mathcal{H}^{m-1}}d\tau\leq \int_0^t\|(p^a, \ln\r)\|^2_{\mathcal{H}^m}d\tau\nonumber\\
&&\leq \f1{c_2}\int_0^t\|p\|^2_{\mathcal{H}^m}d\tau+C[1+P(Q(t))]\int_0^t\|p\|^2_{\mathcal{H}^{m-1}}d\tau,
\end{eqnarray}
here $a\in\mathbb{R}$ is any given constant and $c_2$ is a positive constant depending on $\f1{c_0}, a$.
It follows from \eqref{2.36}-\eqref{2.40},  that
\begin{eqnarray}\label{2.42}
&&-\int_0^t\int\mathcal{Z}^\a p\cdot\mathcal{Z}^\a \mbox{div}udxd\tau\geq \int\f{1}{2\g p}|\MZ^{\a}p|^2dx-\int\f{1}{2\g p_0}|\MZ^{\a}p_0|^2dx\nonumber\\
&&~~~~~~~~~~~~~~~~~-C\d\int_0^t \|\nabla p\|^2_{\mathcal{H}^{m-1}}d\tau
-C_\d [1+P(Q(t))]\int_0^t\|(p,u)\|^2_{\mathcal{H}^{m}}d\tau.
\end{eqnarray}

In order to complete the estimates in \eqref{2.30}, it remains to estimate the terms involving $\mathcal{C}_1^\a$ and $\mathcal{C}_2^\a$.  It follows from
Proposition \ref{prop3.2} and \eqref{2.41} that
\begin{eqnarray}\label{2.43}
&&\int_0^t\|\mathcal{C}_1^\a\|^2 dxd\tau\leq C\sum_{|\beta|\geq1, \b+\g=\a}\int_0^t \|\MZ^\b\r\cdot\MZ^{\g}u_t\|^2d\tau\nonumber\\
&&\leq
C\sup_{0\leq\tau\leq t}\|\MZ p\|^2_{L^{\infty}}\int_0^t\|u_t\|^2_{\mathcal{H}^{m-1}}d\tau
+C\sup_{0\leq\tau\leq t}\|u_t\|^2_{L^{\infty}}\int_0^t\|\r\|^2_{\mathcal{H}^{m}}d\tau\nonumber\\
&&\leq [1+P(Q(t))]\int_0^t\|(p,u)\|^2_{\mathcal{H}^{m}}d\tau,
\end{eqnarray}
and
\begin{eqnarray}\label{2.43-1}
&&\int_0^t\|\mathcal{C}_2^\a\|^2 dxd\tau\leq C\sum_{|\beta|\geq1, \b+\g=\a}\int_0^t \|\MZ^\b(\r u)\cdot\MZ^{\g}\nabla u\|^2d\tau
+\sup_{0\leq\tau t}\|\r u\|_{L^{\infty}}\int_0^t\|\nabla u\|^2_{\mathcal{H}^{m-1}}d\tau\nonumber\\
&&\leq
C[1+P(Q(t))]\int_0^t\|\nabla u\|^2_{\mathcal{H}^{m-1}}+\|\r u\|^2_{\mathcal{H}^{m}}d\tau
\leq
C[1+P(Q(t))]\int_0^tP(\Lambda_m(\tau))d\tau.
\end{eqnarray}
As a consequence of  \eqref{2.43}, \eqref{2.43-1} and the Cauchy inequality, one has that
\begin{eqnarray}\label{2.43-2}
&&\int_0^t\int(\mathcal{C}_1^\a+\mathcal{C}_2^\a)\cdot\mathcal{Z}^\a udxd\tau\leq C[1+P(Q(t))]\int_0^tP(\Lambda_m(\tau))d\tau.
\end{eqnarray}

Therefore, substituting \eqref{2.42} and \eqref{2.43-2} into \eqref{2.30} yields \eqref{2.16}. Thus  the proof of the Lemma \ref{lem2.2} is completed.
 $\hfill\Box$

\subsection{Estimates for $\mbox{div}u$ and  $\nabla{p}$}
\setcounter{equation}{0}

To deal with the compressibility of the system, we need to derive some uniform estimates on $\|\mbox{div}u\|_{\mathcal{H}^{m-1}}$, which will imply the desired  uniform estimates on $\|\nabla u\|_{\mathcal{H}^{m-1}}$.

\begin{lemma}\label{lem4.1}
For every $m\in \mathbb{N}_+$, it holds that
\begin{eqnarray}\label{4.1}
&&\sup_{0\leq \tau\leq t}\Big(\int\f12\r|\mbox{div}u(\tau)|^2+\f{1}{2\g p}|\nabla p(\tau)|^2dx\Big)+\f34\v\int_0^t\|\nabla \mbox{div}u(\tau)\|^2d\tau\nonumber\\
&&\leq \int\f12\r_0|\mbox{div}u_0|^2+\f{1}{2\g p_0}|\nabla p_0|^2dx
+C_3[1+P(Q(t))]\int_0^t\Lambda_m(\tau) d\tau.
\end{eqnarray}

\end{lemma}

\noindent\textbf{Proof}. Multiplying \eqref{2.5} by $\nabla\mbox{div}u$ yields that
\begin{eqnarray}\label{4.2}
&&\int_0^t\int(\r u_t+\r u\cdot\nabla u)\nabla\mbox{div}u dxd\tau
+\int_0^t\int\nabla p\cdot\nabla\mbox{div}u dxd\tau\nonumber\\
&&=-\v\int_0^t\int\nabla\times\omega\cdot\nabla\mbox{div}u dxd\tau
+2\v\int_0^t\|\nabla\mbox{div}u\|^2 d\tau.
\end{eqnarray}
Integrating by parts and using the boundary conditions  \eqref{2.7}, one obtains that
\begin{eqnarray}\label{4.3}
&&\int_0^t\int(\r u_t+\r u\cdot\nabla u)\nabla\mbox{div}u dxd\tau=-\int_0^t\int(\r \mbox{div}u_t+\r u\cdot\nabla \mbox{div}u)\mbox{div}u dxd\tau\nonumber\\
&&~~~-\int_0^t\int(\nabla\r\cdot u_t+\nabla(\r u)^t\cdot\nabla u) \mbox{div}u dxd\tau
+\int_0^t\int_{\partial\Om}\r(u\cdot\nabla)u\cdot n\mbox{div}ud\sigma d\tau\nonumber\\
&&\leq -\f12\int\r|\mbox{div}u|^2dx+\f12\int\r_0|\mbox{div}u_0|^2dx+C[1+P(Q(t))]\int_0^t\|(u_t,\nabla u)\|^2d\tau\nonumber\\
&&~~~~~~~~~~~~~~~~~~~+\Big|\int_0^t\int_{\partial\Om}\r(u\cdot\nabla)n\cdot u\mbox{div}ud\sigma d\tau\Big|\nonumber\\
&&\leq -\f12\int\r|\mbox{div}u|^2dx+\f12\int\r_0|\mbox{div}u_0|^2dx+C[1+P(Q(t))]\int_0^t\|(u_t,\nabla u)\|^2d\tau\nonumber\\
&&~~~~~~~~~~~~~~~~~~~+C_2[1+\|\r \nabla u\|_{L^\infty}]\int_0^t\|\nabla u\|\|u\|d\tau\nonumber\\
&&\leq -\f12\int\r|\mbox{div}u|^2dx+\f12\int\r_0|\mbox{div}u_0|^2dx+C_2[1+P(Q(t))]\int_0^t\Lambda_m(\tau) d\tau.
\end{eqnarray}
Using \eqref{2.35} and integrating by parts lead to
\begin{eqnarray}\label{4.4}
&&\int_0^t\int\nabla{p}\cdot\nabla\mbox{div}u dxd\tau
=-\int_0^t\int\nabla{p}\cdot\nabla(\f{p_t}{\g p}) dxd\tau
-\int_0^t\int\nabla{p}\cdot\nabla(u\cdot\f{\nabla p}{\g p}) dxd\tau\nonumber\\
&&\leq -\int_0^t\int\nabla{p}\cdot\f{\nabla p_t}{\g p} dxd\tau-\int_0^t\int \frac{u}{2\g p}\cdot\nabla(|\nabla p|^2)dxd\tau+C[1+P(Q(t))]\int_0^t\|\nabla p\|^2 d\tau\nonumber\\
&&\leq -\int\f{1}{2\g p}|\nabla{p}|^2dx+\int\f{1}{2\g p_0}|\nabla{p_0}|^2dx
+C[1+P(Q(t))]\int_0^t\|\nabla p\|^2 d\tau,
\end{eqnarray}
and  \eqref{2.7}, together with  integration by parts along the boundary, implies that
\begin{eqnarray}\label{4.5}
&&\v\left|\int_0^t\int\nabla\times\omega\cdot\nabla\mbox{div}u dxd\tau\right|
=\v\left|\int_0^t\int_{\partial\Om} n\times\omega\cdot\nabla\mbox{div}u d\sigma d\tau\right|\nonumber\\
&&=\v\left|\int_0^t\int_{\partial\Om} (Bu)\cdot \Pi(\nabla\mbox{div}u) d\sigma d\tau\right|
=\v\left|\int_0^t\int_{\partial\Om} (Bu)\cdot Z_y\mbox{div}u d\sigma d\tau\right|\nonumber\\
&&\leq  C_3\v\int_0^t |u|_{H^{\f12}}|\mbox{div}u|_{H^{\f12}} d\tau
\leq C_3\v\int_0^t \|u\|_{H^1}\|\mbox{div}u\|_{H^1}d\tau\nonumber\\
&&\leq \f\v4\int_0^t\|\nabla\mbox{div}u\|^2d\tau+C_3\v\int_0^t(\|\nabla u\|^2+\|u\|^2) d\tau.
\end{eqnarray}

Substituting \eqref{4.3}-\eqref{4.5} into \eqref{4.2} proves \eqref{4.1}. Thus the proof of Lemma \ref{lem4.1} is completed.
 $\hfill\Box$

 \

Next, we consider the higher order estimates. One starts with the estimates of  $\MZ^\a\mbox{div}u$ for $|\a_0|\leq m-2$ with $|\a|=m-1$.

\begin{lemma}\label{lem4.2}
For every $m\geq1$ and $|\a|\leq m-1$ with $|\a_0|\leq m-2$, it holds that
\begin{eqnarray}\label{4.6}
&&\sup_{0\leq \tau\leq t}\Big(\int\r|\MZ^\a\mbox{div}u(\tau)|^2+\f{1}{\g p}|\MZ^\a\nabla p(\tau)|^2dx\Big)+\v\int_0^t\|\nabla\MZ^\a\mbox{div}u(\tau)\|^2d\tau\nonumber\\
&&\leq \int\r_0|\MZ^\a\mbox{div}u_0|^2+\f{1}{\g p_0}|\MZ^\a\nabla p_0|^2dx+C C_{m+2}\Big\{\d\int_0^t\|\nabla\partial_t^{m-1}p(\tau)\|^2d\tau\nonumber\\
&&~~~~~+(\d+\v)\int_0^t\|\nabla\MZ^{m-2}\mbox{div}u(\tau)\|^2d\tau
+C_\d [1+P(Q(t))]\int_0^t\Lambda_m(\tau)d\tau\nonumber\\
&&~~~~~~~~~~~~~~~~~~~~+\v\int_0^t\|\nabla^2u(\tau)\|^2_{\mathcal{H}^{m-2}}d\tau\Big\}
\end{eqnarray}
where the last term doesn't appear if $m-2<0$.
\end{lemma}

\noindent\textbf{Proof}. The case for $|\a|=0$ is already proved in Lemma \ref{lem4.1}. Assume that it is proved for $|\a|\leq m-2$, one needs to prove it for $|\a|=m-1\geq 1$ with $|\a_0|\leq m-2$.

Multiplying \eqref{2.17} by $\nabla\MZ^\a\mbox{div}u$ leads to
\begin{eqnarray}\label{4.7}
&&\int_0^t\int(\r\mathcal{Z}^\a u_t+\r u\cdot\nabla\mathcal{Z}^\a u)\nabla\MZ^\a\mbox{div}udxd\tau+\int_0^t\int\mathcal{Z}^\a\nabla p\cdot\nabla\MZ^\a\mbox{div}udxd\tau\nonumber\\
&&=-\v\int_0^t\int\mathcal{Z}^\a\nabla\times\omega\cdot\nabla\MZ^\a\mbox{div}udxd\tau
+2\v\int_0^t\int\mathcal{Z}^\a\nabla\mbox{div}u\cdot\nabla\MZ^\a\mbox{div}udxd\tau\nonumber\\
&&~~~~~~~~~~~~~~~~~+\int_0^t\int(\mathcal{C}_1^\a+\mathcal{C}_2^\a)\nabla\MZ^\a\mbox{div}udxd\tau.
\end{eqnarray}
Integrating by parts gives that
\begin{eqnarray}\label{4.8}
&&\int_0^t\int(\r\mathcal{Z}^\a u_t+\r u\cdot\nabla\mathcal{Z}^\a u)\nabla\MZ^\a\mbox{div}udxd\tau\nonumber\\
&&=-\int_0^t\int(\r\mbox{div}\mathcal{Z}^\a u_t+\r u\cdot\nabla\mbox{div}\mathcal{Z}^\a u)\MZ^\a\mbox{div}udxd\tau\nonumber\\
&&~~~~-\int_0^t\int(\nabla\r\cdot\mathcal{Z}^\a u_t+\nabla(\r u)^t\cdot\nabla\mathcal{Z}^\a u)\MZ^\a\mbox{div}udxd\tau\nonumber\\
&&~~~~+\int_0^t\int_{\partial\Om}(\r\mathcal{Z}^\a u_t\cdot n+\r(u\cdot\nabla)\mathcal{Z}^\a u\cdot n)\MZ^\a\mbox{div}ud\sigma d\tau\triangleq I_1+I_2+I_3.
\end{eqnarray}
For $I_1$ and $I_2$, one   obtains easily that
\begin{eqnarray}\label{4.9}
I_1&&=-\int_0^t\int(\r\mathcal{Z}^\a \mbox{div}u_t+\r u\cdot\nabla\mathcal{Z}^\a \mbox{div}u)\MZ^\a\mbox{div}udxd\tau\nonumber\\
&&~~~~~~~~~-\int_0^t\int\Big(\r[\mbox{div},\mathcal{Z}^\a]u_t+\r(u_1Z_{y^1}+u_2Z_{y^2}+\frac{u\cdot N}{\varphi(z)}Z_3)[\mbox{div},\mathcal{Z}^\a]u\Big)\MZ^\a\mbox{div}udxd\tau\nonumber\\
&&\leq -\int\frac{\r}{2}|\MZ^\a\mbox{div}u(t)|^2dx+\int\frac{\r_0}{2}|\MZ^\a\mbox{div}u_0|^2dx
+C_2[1+P(Q(t))]\int_0^t\Lambda_m(\tau)d\tau,
\end{eqnarray}
and
\begin{eqnarray}\label{4.10}
I_2\leq C[1+P(Q(t))]\int_0^t\Lambda_m(\tau) d\tau.
\end{eqnarray}
Noting that $\MZ^\a$ contains at least one tangential derivative $Z_y$, integrating by parts along the boundary and using \eqref{4.11} , one obtains that
\begin{eqnarray}\label{4.12}
&&I_3=\int_0^t\int_{\partial\Om}[\r\mathcal{Z}^\a u_t\cdot n-\r(u\cdot\nabla)n\cdot\mathcal{Z}^\a u+\r(u\cdot\nabla)(\mathcal{Z}^\a u\cdot n)]\MZ^\a\mbox{div}ud\sigma d\tau\nonumber\\
&&~~~=\int_0^t\int_{\partial\Om}[\r\mathcal{Z}^\a u_t\cdot n-\r(u\cdot\nabla)n\cdot\mathcal{Z}^\a u+\r(u_1Z_{y^1}+u_2Z_{y^2})(\mathcal{Z}^\a u\cdot n)]\MZ^\a\mbox{div}ud\sigma d\tau\nonumber\\
&&~~~\leq  C[1+P(Q(t))]\int_0^t\Big(|\mathcal{Z}^\a u_t\cdot n|_{H^{\f12}}+|\mathcal{Z}^\a u|_{H^{\f12}}+|\mathcal{Z}^\a u\cdot n|_{H^{\f32}}\Big)\cdot|\MZ^{m-2}\mbox{div}u|_{H^{\f12}}d\tau\nonumber\\
&&~~\leq  C_\d C_{m+2}[1+P(Q(t))]\int_0^t\Lambda_m(\tau)d\tau
+\d\int_0^t\|\nabla\MZ^{m-2}\mbox{div}u\|^2d\tau.
\end{eqnarray}
Substituting \eqref{4.9},\eqref{4.10} and \eqref{4.12} into \eqref{4.8}, yields that
\begin{eqnarray}\label{4.13}
&&\int_0^t\int(\r\mathcal{Z}^\a u_t+\r u\cdot\nabla\mathcal{Z}^\a u)\nabla\MZ^\a\mbox{div}udxd\tau\nonumber\\
&&\leq -\int\frac{\r}{2}|\MZ^\a\mbox{div}u(t)|^2dx+\int\frac{\r_0}{2}|\MZ^\a\mbox{div}u_0|^2dx
+\d\int_0^t\|\nabla\MZ^{m-2}\mbox{div}u\|^2d\tau\nonumber\\
&&~~~~~~~~+C_\d[1+P(Q(t))]\int_0^t\Lambda_m(\tau)d\tau.
\end{eqnarray}
It is easy to obtain
\begin{eqnarray}\label{4.14}
&&\v\int_0^t\int\mathcal{Z}^\a\nabla\mbox{div}u\cdot\nabla\MZ^\a\mbox{div}udxd\tau\nonumber\\
&&= \v\int_0^t\|\nabla\MZ^\a\mbox{div}u\|^2d\tau
+\v\int_0^t\int[\mathcal{Z}^\a,\nabla]\mbox{div}u\cdot\nabla\MZ^\a\mbox{div}udxd\tau\nonumber\\
&&\geq \f34\v\int_0^t\|\nabla\MZ^\a\mbox{div}u\|^2d\tau
-C\v\int_0^t\|\nabla\MZ^{m-2}\mbox{div}u\|^2+\Lambda_m(\tau)d\tau.
\end{eqnarray}
It follows from \eqref{2.7} and Proposition \ref{prop3.3},   that for $k\leq m-1$
\begin{eqnarray}\label{4.15}
&&|n\times Z_y^{m-k-1}\partial_t^k\omega|_{H^{\f12}}\leq C_{m+2}|\partial_t^k\omega|_{H^{m-k-\f32}}
+C_{m+2}|\partial_t^k u|_{H^{m-k-\f12}}\nonumber\\
&&\leq C_{m+2}\|\nabla^2u\|^{\f12}_{\mathcal{H}^{m-2}}\cdot\|\nabla u\|^{\f12}_{\mathcal{H}^{m-1}}
+C_{m+2}\Big(\|\nabla u\|_{\mathcal{H}^{m-1}}+\|u\|_{\mathcal{H}^{m}}\Big).
\end{eqnarray}
This, together with integrating by parts, shows that
\begin{align}\label{4.16}
&-\v\int_0^t\int\mathcal{Z}^\a\nabla\times\omega\cdot\nabla\MZ^\a\mbox{div}udxd\tau\nonumber\\
&=-\v\int_0^t\int\nabla\times\mathcal{Z}^\a\omega\cdot\nabla\MZ^\a\mbox{div}udxd\tau
-\v\int_0^t\int[\mathcal{Z}^\a,\nabla\times]\omega\cdot\nabla\MZ^\a\mbox{div}udxd\tau\nonumber\\
&\geq -\f{\v}8\int_0^t\|\nabla\MZ^\a\mbox{div}u\|^2d\tau -\v\int_0^t\int_{\partial\Om}n\times\mathcal{Z}^\a\omega\cdot\Pi(\nabla\MZ^\a\mbox{div}u)d\sigma d\tau\nonumber\\
&~~~~~~~~~~~~~~~~~~~-C\v\int_0^t\|\nabla\MZ^{m-2}\omega\|^2+\Lambda_m(\tau)d\tau\nonumber\\
&\geq -\f{\v}8\int_0^t\|\nabla\MZ^\a\mbox{div}u\|^2d\tau -C\v\int_0^t |n\times\mathcal{Z}^\a\omega|_{H^{\f12}}\cdot|\MZ^\a\mbox{div}u|_{H^{\f12}}d\tau\nonumber\\
&~~~~~~~~~~~~~~~~~~~-C\v\int_0^t\|\nabla\MZ^{m-2}\omega\|^2+\Lambda_m(\tau)d\tau\nonumber\\
&\geq -\f\v4\int_0^t\|\nabla\MZ^\a\mbox{div}u\|^2d\tau -C\v\int_0^t\|\nabla^2u\|^2_{\mathcal{H}^{m-2}}d\tau - C_{m+2} \int_0^tP(\Lambda_m(\tau)) d\tau.
\end{align}

Now we estimate the terms involving the pressure. It follows from  \eqref{2.35} that
\begin{eqnarray}\label{4.17}
&&\int_0^t\int\mathcal{Z}^\a\nabla p\cdot\nabla\MZ^\a\mbox{div}udxd\tau\nonumber\\
&&\leq \int_0^t\int\mathcal{Z}^\a\nabla p\cdot\MZ^\a\nabla\mbox{div}udxd\tau
+\d\int_0^t\|\nabla\MZ^{m-2}\mbox{div}u\|^2d\tau+C_\d\int_0^t\Lambda_m(\tau)d\tau\nonumber\\
&&\leq -\int_0^t\int\mathcal{Z}^\a\nabla p\cdot\MZ^\a\nabla(\f{p_t}{\g p})dxd\tau
-\int_0^t\int\mathcal{Z}^\a\nabla p\cdot\MZ^\a\nabla(\f{u}{\g p}\cdot\nabla p)dxd\tau\nonumber\\
&&~~~~~~~~+\d\int_0^t\|\nabla\MZ^{m-2}\mbox{div}u\|^2d\tau+C_\d\int_0^t\Lambda_m(\tau)d\tau.
\end{eqnarray}
For the first term on the right hand side of \eqref{4.17}, one notices that
\begin{eqnarray}\label{4.18}
&&\MZ^\a\nabla(\f{p_t}{\g p})=\f{1}{\g p}\mathcal{Z}^{\a}\nabla p_t
+\sum_{|\beta|\geq1, \b+\g=\a}C_{\b,\g}\MZ^\b(\f1{\g p})\cdot\MZ^{\g}\nabla p_t\nonumber\\
&&~~~~~~~~~~~~~~~~~~~~+\sum_{\b+\g=\a}C_{\b,\g}\MZ^\b p_t\cdot\MZ^{\g}\nabla(\f1{\g p}).\nonumber
\end{eqnarray}
Therefore,  Proposition \ref{prop3.2} shows that
\begin{eqnarray}\label{4.19}
&&-\int_0^t\int\mathcal{Z}^\a\nabla p\cdot\MZ^\a\nabla(\f{p_t}{\g p})dxd\tau\leq -\int\f{1}{2\g p}|\mathcal{Z}^\a\nabla p|^2 dx+\int\f{1}{2\g p_0}|\mathcal{Z}^\a\nabla p_0|^2dx\nonumber\\
&&~~~~~~~~~~~~~~+\d\int_0^t\|\nabla\partial_t^{m-1}p\|^2d\tau +C_\d(1+P(Q(t)))\int_0^t\Lambda_m(\tau)d\tau.
\end{eqnarray}
To estimate the second term on the right hand side of \eqref{4.17}, we notice that
\begin{eqnarray}\label{4.20}
&&\MZ^\a\nabla(\f{u}{\g p}\cdot\nabla p)=\sum_{i=1,2}\f{u_i}{\g p}\mathcal{Z}^{\a}\nabla\partial_{y^i} p
+\f{u\cdot N}{\g p}\mathcal{Z}^{\a}\partial_z\nabla p\nonumber\\
&&~~~~~+\sum_{i=1,2}\sum_{|\beta|\geq1, \b+\g=\a}C_{\b,\g}\MZ^\b(\f{u_i}{\g p})\cdot\MZ^{\g}\nabla\partial_{y^i} p+\sum_{|\beta|\geq1, \b+\g=\a}C_{\b,\g}\MZ^\b(\f{u\cdot N}{\g p})\cdot\MZ^{\g}\partial_z\nabla p\nonumber\\
&&~~~~~+\sum_{i=1,2}\sum_{\b+\g=\a}C_{\b,\g}\MZ^\b\nabla(\f{u_i}{\g p})\cdot\MZ^{\g}\partial_{y^i} p
+\sum_{\b+\g=\a}C_{\b,\g}\MZ^\b\nabla(\f{u\cdot N}{\g p})\cdot\MZ^{\g}\partial_z p.
\end{eqnarray}
Then integrating by parts gives that
\begin{eqnarray}\label{4.21-1}
&&-\int_0^t\int\mathcal{Z}^\a\nabla p\cdot\Big(\sum_{i=1,2}\f{u_i}{\g p}\mathcal{Z}^{\a}\nabla\partial_{y^i} p+\f{u\cdot N}{\g p}\mathcal{Z}^{\a}\partial_z\nabla p\Big) dxd\tau\nonumber\\
&&=-\int_0^t\int\mathcal{Z}^\a\nabla p\cdot\Big(\sum_{i=1,2}\f{u_i}{\g p}\partial_{y^i} \mathcal{Z}^{\a}\nabla p+\f{u\cdot N}{\g p}\partial_z\mathcal{Z}^{\a}\nabla p\Big) dxd\tau\nonumber\\
&&~~~~~~~~~~~~~~-\int_0^t\int\mathcal{Z}^\a\nabla p\cdot\Big(\sum_{i=1,2}\f{u_i}{\g p}[\mathcal{Z}^{\a}\nabla,\partial_{y^i}]  p+\f{u\cdot N}{\g p\varphi(z)}\varphi(z)[\mathcal{Z}^{\a}, \partial_z]\nabla p\Big) dxd\tau\nonumber\\
&&\leq C_2[1+P(Q(t))]\int_0^t\Lambda_m(\tau)d\tau.
\end{eqnarray}
Using Proposition \ref{prop3.2}, one has that
\begin{eqnarray}\label{4.21-2}
&&-\sum_{i=1,2}\Big(\sum_{\b+\g=\a}C_{\b,\g}\int_0^t\int\mathcal{Z}^\a\nabla p\cdot\MZ^\b\nabla(\f{u_i}{\g p})\cdot\MZ^{\g}\partial_{y^i} pdxd\tau\nonumber\\
&&~~~~~~~~~~~~~~~~~-\sum_{|\beta|\geq1, \b+\g=\a}C_{\b,\g} \int_0^t\int\mathcal{Z}^\a\nabla p\cdot\MZ^\b(\f{u_i}{\g p})\cdot\MZ^{\g}\nabla\partial_{y^i} pdxd\tau\Big)\nonumber\\
&&\leq \d\int_0^t\|\nabla\partial_t^{m-1}p\|^2d\tau +C_\d(1+P(Q(t)))\int_0^t\Lambda_m(\tau)d\tau,
\end{eqnarray}
and
\begin{eqnarray}\label{4.21-4}
&&-\sum_{\b+\g=\a}C_{\b,\g}\int_0^t\int\mathcal{Z}^\a\nabla p\cdot\MZ^\b\nabla(\f{u\cdot N}{\g p})\cdot\MZ^{\g}\partial_z p dxd\tau\nonumber\\
&&\leq \d\int_0^t\|\nabla\partial_t^{m-1}p\|^2d\tau + C_\d C_{m+1}(1+P(Q(t)))\int_0^t\Lambda_m(\tau)d\tau.
\end{eqnarray}
On the other hand, notice that for $ |\beta|\geq1, \b+\g=\a,$ and $|\a|=m-1$
\begin{eqnarray}\label{4.21-3}
&&\MZ^\b(\f{u\cdot N}{p})\cdot\MZ^{\g}\partial_z\nabla p=\f1{\varphi(z)}\MZ^\b(\f{u\cdot N}{ p})\cdot\varphi(z)\MZ^{\g}\partial_z\nabla p\nonumber\\
&&=\sum_{\tilde{\beta}\leq \beta,\tilde{\g}\leq \g}C_{\tilde\beta,\tilde\g}\MZ^{\tilde\b}(\f{u\cdot N}{ p\varphi(z)})\cdot\MZ^{\tilde\g}(Z_3\nabla p),\nonumber
\end{eqnarray}
where $|\tilde\beta|+|\tilde\g|\leq m-1$, $|\tilde\g|\leq m-2$ and $C_{\tilde\beta,\tilde\g} $ is some smooth bounded coefficient.

If $\tilde\beta=0$, and hence $|\tilde\g|\leq m-2$, it holds that
\begin{eqnarray}\label{4.21-5}
&&\int_0^t\|\MZ^{\tilde\b}(\f{u\cdot N}{ p\varphi(z)})\cdot\MZ^{\tilde\g}(Z_3\nabla p)\|^2d\tau\leq C\|\f{u\cdot N}{ p\varphi(z)}\|^2_{L^\infty}\int_0^t\|Z_3\nabla p(\tau)\|^2_{\mathcal{H}^{m-2}}d\tau\nonumber\\
&&\leq C_2(1+\|u\|^2_{W^{1,\infty}})\int_0^t\Lambda_m(\tau)d\tau\leq
C_2(1+P(Q(t)))\int_0^t\Lambda_m(\tau)d\tau,
\end{eqnarray}
where one has used $\|\f{u\cdot N}{\varphi}\|_{L^\infty}\leq C\|u\|_{W^{1,\infty}}$ due to the boundary condition $u\cdot n=0$.
On the other hand, by using Proposition \ref{prop3.2} for $\tilde\b\neq0$, one obtains that
\begin{eqnarray}\label{4.21-6}
&&-\sum_{\tilde{\beta}\leq \beta,\tilde{\g}\leq \g}C_{\tilde\beta,\tilde\g}\int_0^t\int\mathcal{Z}^\a\nabla p\cdot\MZ^{\tilde\b}(\f{u\cdot N}{ p\varphi(z)})\cdot\MZ^{\tilde\g}(Z_3\nabla p) dxd\tau\nonumber\\
&&\leq C\Big(\int_0^t\|\mathcal{Z}^\a\nabla p\|^2d\tau\Big)^{\f12}
\cdot\sup_{0\leq\tau\leq t}\|\MZ(\f{u\cdot N}{ p\varphi(z)},\MZ\nabla p)\|_{L^\infty}
\Big\{\int_0^t\|\MZ(\f{u\cdot N}{ p\varphi(z)})\|^2_{\mathcal{H}^{m-2}}+\|Z_3\nabla p\|^2_{\mathcal{H}^{m-2}}\Big\}^{\f12}\nonumber\\
&&\leq  C_{m+1}[1+P(Q(t))]\int_0^t\Lambda_m(\tau)d\tau,
\end{eqnarray}
where one has used the following hardy inequality in the last step of \eqref{4.21-6}
\begin{equation}\label{4.21-7}
\|\MZ(\f{u\cdot N}{ \varphi(z)})\|_{\mathcal{H}^{m-2}}\leq C_{m+1}(\|\nabla u\|_{\mathcal{H}^{m-1}}+\|u\|_{\mathcal{H}^{m-1}}),
\end{equation}
which is proved in page543  of \cite{Masmoudi-R}. Then \eqref{4.21-5} and \eqref{4.21-6} yield immediately
\begin{eqnarray}\label{4.21-8}
&&-\sum_{|\b|\geq1, \b+\g=\a}C_{\b,\g}\int_0^t\int\mathcal{Z}^\a\nabla p\cdot\MZ^\b(\f{u\cdot N}{p})\cdot\MZ^{\g}\partial_z\nabla p dxd\tau\nonumber\\
&&\leq C_{m+1}(1+P(Q(t)))\int_0^t\Lambda_m(\tau)d\tau.
\end{eqnarray}
Combining \eqref{4.20}-\eqref{4.21-4} with  \eqref{4.21-8}, one obtains that
\begin{eqnarray}\label{4.21}
&&-\int_0^t\int\mathcal{Z}^\a\nabla p\cdot\MZ^\a\nabla(\f{u}{\g p}\cdot\nabla p)dxd\tau\nonumber\\
&&\leq C\d\int_0^t\|\nabla\partial_t^{m-1}p\|^2d\tau+C_\d C_{m+1}(1+P(Q(t)))\int_0^t\Lambda_m(\tau)d\tau.
\end{eqnarray}
Due to \eqref{4.19}, \eqref{4.21} and \eqref{4.17}, it holds that
\begin{eqnarray}\label{4.22}
&&\int_0^t\int\mathcal{Z}^\a\nabla p\cdot\nabla\MZ^\a\mbox{div}udxd\tau\nonumber\\
&&\leq -\int\f{1}{2\g p}|\mathcal{Z}^\a\nabla p|^2 dx+\int\f{1}{2\g p_0}|\mathcal{Z}^\a\nabla p_0|^2dx +C\d\int_0^t\|\nabla\MZ^{m-2}\mbox{div}u\|^2d\tau\nonumber\\
&&~~~~~~~~~~~~~~+C\d\int_0^t\|\nabla\partial_t^{m-1}p\|^2d\tau+C_\d C_{m+1}(1+P(Q(t)))\int_0^t\Lambda_m(\tau)d\tau.
\end{eqnarray}

\

Finally,  using Proposition \ref{prop3.2} and integrating by parts, one can get that
\begin{eqnarray}\label{4.23}
&&|\int_0^t\int(\mathcal{C}_1^\a+\mathcal{C}_2^\a)\nabla\MZ^\a\mbox{div}udxd\tau|\nonumber\\
&&\leq |\int_0^t\int(\mathcal{C}_1^\a+\mathcal{C}_2^\a)Z\nabla\MZ^{\a-1}\mbox{div}udxd\tau|
+C\int_0^t\int(|\mathcal{C}_1^\a|+|\mathcal{C}_2^\a|)|\nabla\MZ^{\a-1}\mbox{div}u|dxd\tau\nonumber\\
&&\leq |\int_0^t\int(|Z\mathcal{C}_1^\a|+|Z\mathcal{C}_2^\a|+|\mathcal{C}_1^\a|+|\mathcal{C}_2^\a|)
\cdot|\nabla\MZ^{\a-1}\mbox{div}u|dxd\tau\nonumber\\
&&\leq \d\int_0^t\|\nabla\MZ^{m-2}\mbox{div}u\|^2d\tau+C\int_0^t\|(\mathcal{C}_1^\a,\mathcal{C}_2^\a)\|^2_{H^1_{co}}d\tau\nonumber\\
&&\leq \d\int_0^t\|\nabla\MZ^{m-2}\mbox{div}u\|^2d\tau+C_\d(1+P(Q(t)))\int_0^t\Lambda_m(\tau)d\tau,
\end{eqnarray}
where the following estimate has been used:
\begin{equation}\label{4.24}
\int_0^t\|(\mathcal{C}_1^\a,\mathcal{C}_2^\a)\|^2_{\mathcal{H}^1_{co}}d\tau
\leq C(1+P(Q(t)))\int_0^t\Lambda_m(\tau)d\tau.
\end{equation}
Substituting \eqref{4.13}, \eqref{4.14}, \eqref{4.16}, \eqref{4.22} and \eqref{4.23} into \eqref{4.7} proves \eqref{4.6}. Therefore,  the proof of Lemma \ref{lem4.2} is completed.     $\hfill\Box$

\begin{remark}
It should be pointed out that in general, it is hard to derive an uniform estimate to the term $\int_0^t\| \MZ^{m-2}\partial_{zz}u\|^2d\tau$ due to the possible appearance of boundary layers.	
However, the term  $\int_0^t\|\nabla\MZ^{m-2}\mbox{div}u\|^2d\tau$ is expected to be controllable since one could expect that there is no strong boundary layer in either $\mbox{div}u$ or pressure $p$.
\end{remark}

Although we have obtained the bound for $\|\MZ^\a(\mbox{div}u,\nabla p)\|$ for $|\a|=m-1$ with $|\a_0|\leq m-2$. Yet, the estimates on  $\|\partial_t^{m-1}(\mbox{div}u,\nabla p)\|$ is weaker:
\begin{lemma}\label{lem4.3}
For every $m\geq1$, it holds that
\begin{eqnarray}\label{4.25}
&&\sup_{0\leq \tau\leq t}\Big(\v\int\r|\partial_t^{m-1}\mbox{div}u(\tau)|^2+\f{1}{\g p}|\partial_t^{m-1}\nabla p(\tau)|^2dx\Big)+\v^2\int_0^t\|\nabla\partial_t^{m-1}\mbox{div}u(\tau)\|^2d\tau\nonumber\\
&&\leq \v\int\r_0|\partial_t^{m-1}\mbox{div}u_0|^2+\f{1}{\g p_0}|\partial_t^{m-1}\nabla p_0|^2dx+C_\d C_{m+1}[1+P(Q(t))]\int_0^t\Lambda_m(\tau) d\tau.
\end{eqnarray}

\end{lemma}

\noindent\textbf{Proof}. Applying $\partial_t^{m-1}$ to \eqref{2.5} shows that
\begin{equation}\label{4.26}
\r\partial_t^{m-1}u_t+\r u\cdot\nabla\partial_t^{m-1}u+\partial_t^{m-1}\nabla p=-\v\nabla\times\partial_t^{m-1}\omega+2\v\nabla\partial_t^{m-1}\mbox{div}u+\mathcal{C}_1^{m-1}+\mathcal{C}_2^{m-1},
\end{equation}
where
\begin{equation}\label{4.27}
\mathcal{C}_1^{m-1}=-[\partial_t^{m-1}, \r]u_t=-\sum_{|\b|\geq 1,\b+\g=m-1}C_{\b,\g}\partial_t^\b\r\partial_t^\g u_t,
\end{equation}
and
\begin{equation}\label{4.28}
\mathcal{C}_2^{m-1}=-[\partial_t^{m-1}, \r u\cdot\nabla]u=-\sum_{|\b|\geq 1,\b+\g=m-1}C_{\b,\g}\partial_t^\b(\r u)\partial_t^\g\nabla u.
\end{equation}
The boundary conditions become
\begin{equation}\label{4.29}
n\cdot\partial_t^{m}u=0,~~n\times\partial_t^{m-1}\omega=[B\partial_t^{m-1}u]_\tau, ~~x\in \partial\Om.
\end{equation}

Multiplying \eqref{4.26} by $\v\nabla\mbox{div}\partial_t^{m-1}u$ yields that
\begin{eqnarray}\label{4.30}
&&\v\int_0^t\int(\r\partial_t^{m-1}u_t+\r u\cdot\nabla\partial_t^{m-1}u)\nabla\mbox{div}\partial_t^{m-1}udxd\tau
+\v\int_0^t\int\partial_t^{m-1}\nabla p\cdot\nabla\mbox{div}\partial_t^{m-1}udxd\tau\nonumber\\
&&=-\v^2\int_0^t\int\nabla\times\partial_t^{m-1}\omega\cdot\nabla\mbox{div}\partial_t^{m-1}udxd\tau
+2\v^2\int_0^t\|\nabla\partial_t^{m-1}\mbox{div}u\|^2d\tau\nonumber\\
&&~~~~~~~~~~~~~~~~~+\v\int_0^t\int(\mathcal{C}_1^{m-1}+\mathcal{C}_2^{m-1})\nabla\mbox{div}\partial_t^{m-1}udxd\tau
\end{eqnarray}
It follows from  \eqref{4.29} and  integrating by parts that
\begin{eqnarray}\label{4.32}
&&\v^2|\int_0^t\int\nabla\times\partial_t^{m-1}\omega\cdot\nabla\mbox{div}\partial_t^{m-1}udxd\tau|
=\v^2|\int_0^t\int_{\partial\Omega} n\times\partial_t^{m-1}\omega\cdot\Pi(\nabla\mbox{div}\partial_t^{m-1}u)d\sigma d\tau|\nonumber\\
&&\lesssim \v^2\int_0^t|n\times\partial_t^{m-1}\omega|_{H^{\f12}}\cdot|\mbox{div}\partial_t^{m-1}u|_{H^{\f12}}d\tau
\lesssim  C_3 \v^2\int_0^t|\partial_t^{m-1}u|_{H^{\f12}}\cdot|\mbox{div}\partial_t^{m-1}u|_{H^{\f12}}d\tau\nonumber\\
&&\lesssim  C_3 \v^2\int_0^t\|\partial_t^{m-1}u\|_{H^{1}}\cdot\|\mbox{div}\partial_t^{m-1}u\|_{H^{1}}d\tau\nonumber\\
&&\leq \d\v^4\int_0^t\|\nabla\partial_t^{m-1}\mbox{div}u\|^2d\tau+C_3C_\d\int_0^t\Lambda_m(\tau)d\tau.
\end{eqnarray}
and
\begin{eqnarray}\label{4.31}
&&\v\int_0^t\int(\r\partial_t^{m-1}u_t+\r u\cdot\nabla\partial_t^{m-1}u)\nabla\mbox{div}\partial_t^{m-1}udxd\tau\nonumber\\
&&=-\v\int_0^t\int(\r\partial_t^{m-1}\mbox{div}u_t+\r u\cdot\nabla\partial_t^{m-1}\mbox{div}u)\mbox{div}\partial_t^{m-1}udxd\tau\nonumber\\
&&~~~~~~~-\v\int_0^t\int(\nabla\r\cdot\partial_t^{m-1}u_t+\nabla(\r u)^t\cdot\nabla\partial_t^{m-1}u)\mbox{div}\partial_t^{m-1}udxd\tau\nonumber\\
&&~~~~~~~-\v\int_0^t\int_{\partial\Om}\r(u\cdot\nabla)n\cdot\partial_t^{m-1}u \partial_t^{m-1}\mbox{div}ud\sigma d\tau\nonumber\\
&&\leq -\v\int\f{\r}{2}|\partial_t^{m-1}\mbox{div}u(t)|^2dx+\v\int\f{\r_0}{2}|\partial_t^{m-1}\mbox{div}u_0|^2dx
+C[1+P(Q(t))] \int_0^t\Lambda_m(\tau) d\tau\nonumber\\
&&~~~~~+C[1+P(Q(t))]\v\int_0^t\|\partial_t^{m-1}u\|^{\f12}\|\partial_t^{m-1}u\|^{\f12}_{H^1}
\cdot\|\partial_t^{m-1}\mbox{div}u\|^{\f12}_{H^1}\|\partial_t^{m-1}\mbox{div}u\|^{\f12}d\tau\nonumber\\
&&\leq -\v\int\r|\partial_t^{m-1}\mbox{div}u(t)|^2dx+\v\int\r_0|\partial_t^{m-1}\mbox{div}u_0|^2dx
\nonumber\\
&&~~~~~~~~~~~~~~+C[1+P(Q(t))] \int_0^t\Lambda_m(\tau) d\tau+\f18 \v^2 \int_0^t\|\nabla\partial_t^{m-1}\mbox{div}u\|^2d\tau.
\end{eqnarray}
By similar arguments as in  the proof of \eqref{4.22}, one gets that
\begin{eqnarray}\label{4.33}
&&\v\int_0^t\int\partial_t^{m-1}\nabla p\cdot\nabla\mbox{div}\partial_t^{m-1}udxd\tau\leq -\v\int\f1{2\g p}|\nabla\partial_t^{m-1}p|^2dx+\v\int\f1{2\g p_0}|\nabla\partial_t^{m-1}p_0|^2dx\nonumber\\
&&~~~~~~~~~~~~~~~~~~~~~~~~~~~~~~~~~~~~~~~~~~~~~~~~~+C_{m+1}[1+P(Q(t))]\int_0^t\Lambda_m(\tau)d\tau,
\end{eqnarray}
and by using \eqref{4.24}, one can obtain that
\begin{eqnarray}\label{4.34}
&&\v|\int_0^t\int(\mathcal{C}_1^{m-1}+\mathcal{C}_2^{m-1})\nabla\mbox{div}\partial_t^{m-1}udxd\tau|
\leq \f{1}{8}\v^2\int_0^t\|\nabla\partial_t^{m-1}\mbox{div}u\|^2d\tau
+C\int_0^t\|\mathcal{C}_1^{m-1}\|^2+\|\mathcal{C}_2^{m-1}\|^2d\tau\nonumber\\
&&\leq \f{1}{8}\v^2\int_0^t\|\nabla\partial_t^{m-1}\mbox{div}u\|^2d\tau
+C[1+P(Q(t))]\int_0^t\Lambda_m(\tau)d\tau.
\end{eqnarray}

Substituting \eqref{4.32}-\eqref{4.34} into \eqref{4.30} proves \eqref{4.25}. Therefore Lemma \ref{lem4.3} is proved.
$\hfill\Box$

Since the estimate in Lemma \ref{lem4.3} is not enough  to obtain the uniform estimate for $\nabla\partial_t^{m-1}u$, so we need some new estimates on $\|\partial_t^{m-1}\mbox{div}u\|$. Fortunately,  we have the following subtle control about $\|\partial_t^{m-1}\mbox{div}u\|$:
\begin{lemma}\label{lem4.7}
Define
\begin{equation}\label{4.53}
\Lambda_{1m}(t)=\|(p, u)(t)\|^2_{\mathcal{H}^m}+\sum_{|\b|\leq m-2}\|\MZ^{\b}\nabla p(t)\|_1^2+\sum_{|\b|\leq m-2}\|\MZ^{\b}\nabla u(t)\|_1^2.
\end{equation}
Then, for every $m\geq3$,  it holds that
\begin{eqnarray}\label{4.48}
\|\partial_t^{m-1}\mbox{div}u(t)\|^2\leq C_2 \Big\{P(\Lambda_{1m}(t))
+P(Q(t))\Big\}.
\end{eqnarray}
\end{lemma}
\begin{remark}
It should be pointed out that it does not contain the terms $\|\nabla\partial_t^{m-1}u\|$ and $\|\nabla\partial_t^{m-1}p\|$ in the right hand side of \eqref{4.48}. This estimate allows one to obtain the uniform estimates for $\|\nabla\partial_t^{m-1}u\|$.
\end{remark}

\noindent\textbf{Proof}.
Applying $ \partial_t^{m-1}$ to \eqref{2.35} yields that
\begin{equation}\label{4.49}
\partial_t^{m-1}\mbox{div}u=-\partial_t^{m}(\ln\r)-\partial_t^{m-1}(u_i\partial_{y^i}\ln\r)
-\partial_t^{m-1}(u\cdot N\partial_z\ln\r).
\end{equation}
Since
\begin{equation*}\label{4.50}
\partial_t^m(\ln\r)=\partial_t^{m-1}(\f{\partial_t p}{p})=\sum_{k=0}^{m-1}\partial_t^k (\f1p)\cdot\partial_t^{m-k}p,
\end{equation*}
it holds that
\begin{equation}\label{4.51}
\|\partial_t^m(\ln\r)\|^2\leq C\|\partial_t^mp\|^2+\sum_{k=1}^{[\f{m}2]}\|\partial_t^k (\f1p)\|^2_{L^\infty}\cdot\|\partial_t^{m-k}p\|^2+\sum_{k=1+[\f{m}2]}^{m-1}\|\partial_t^{m-k}p\|^2_{L^\infty}\cdot\|\partial_t^k (\f1p)\|^2.
\end{equation}
Note that for $|\b|\leq [\f{m}2]$ and $m\geq 3$,  it holds that
\begin{eqnarray}\label{4.52}
&&\sum_{|\b|\leq [\f{m}2]}\|\MZ^{\b} p\|^2_{L^\infty}\leq C\sum_{|\b|\leq [\f{m}2]}\|(\nabla\MZ^{\b} p,\MZ^{\b} p)\|_1\cdot\|\MZ^{\b} p\|_2\nonumber\\
&&~~~~~~~~~~~~~~~~~~~~~~\leq C\Big( \|p\|^2_{\mathcal{H}^m}+\sum_{|\b|\leq m-2}\|\MZ^{\b}\nabla p\|^2_1\Big)\leq CP(\Lambda_{1m}(t)),
\end{eqnarray}
where  \eqref{3.3}  has been used. Therefore, substituting \eqref{4.52} into \eqref{4.51} leads to that
\begin{equation}\label{4.54}
\|\partial_t^m(\ln\r)\|^2\leq P(\Lambda_{1m}(t)).
\end{equation}
Similar arguments as \eqref{4.52} and  \eqref{4.54} show  that
\begin{equation}\label{4.55}
\|\ln\r\|^2_{\mathcal{H}^m}+\sum_{|\b|\leq m-2}\|\MZ^\b p\|^2_{L^\infty} \leq P(\Lambda_{1m}(t)),
\end{equation}
and
\begin{equation}\label{4.56}
\sum_{|\b|\leq m-2}\|\MZ^\b u\|^2_{L^\infty}+\sum_{|\b|\leq [\f{m}2]}\|\MZ^{\b} u\|^2_{L^\infty}\leq P(\Lambda_{1m}(t)).
\end{equation}
It follows from  \eqref{4.55} and \eqref{4.56}  that
\begin{eqnarray}\label{4.57}
&&\|\partial_t^{m-1}(u_i\partial_{y^i}\ln\r)\|^2\leq \sum_{k=0}^{[\f{m}2]}\|\partial_t^k u\|^2_{L^\infty}\cdot\|\partial_t^{m-k-1}\partial_{y^i}\ln\r\|^2\nonumber\\
&&~~~~+\sum_{k=1+[\f{m}2]}^{m-1}\|\partial_t^{m-k-1}\partial_{y^i}\ln\r\|^2_{L^\infty}\cdot\|\partial_t^k u\|^2\leq C P(\Lambda_{1m}(t)).
\end{eqnarray}

Additional efforts are needed to bound $\partial_t^{m-1}(u\cdot N\partial_z\ln\r)$,  since it involves $\partial_z\ln\r$.  First, rewrite this term as
\begin{eqnarray*}\label{4.58}
&&\partial_t^{m-1}(u\cdot N\partial_z\ln\r)
=\f{u\cdot N}{\varphi(z)}Z_3\partial_t^{m-1}\ln\r+\partial_z\ln\r\cdot\partial_t^{m-1}(u\cdot N)\nonumber\\
&&~~~~~~~~~~~~~~+\sum_{k=1}^{m-2}C_{k} \partial_t^{k}(u\cdot N)\cdot\partial_z\partial_t^{m-1-k}\ln\r.
\end{eqnarray*}
Hence, by using \eqref{4.55} and \eqref{4.56}, one has  that
\begin{eqnarray}\label{4.59}
&&\|\partial_t^{m-1}(u\cdot N\partial_z\ln\r)\|^2\leq C\|\f{u\cdot N}{\varphi(z)}\|^2_{L^\infty}\|\ln\r\|^2_{\mathcal{H}^m}
+C\|\nabla p\|^2_{L^\infty}\|\partial_t^{m-1}u\|^2\nonumber\\
&&~~~~~~~~~~~~~~~~~~~~~~~~~~~~+C\sum_{k=1}^{m-2} \|\partial_t^{k}u\|^2_{L^\infty}\cdot\|\partial_t^{m-1-k}(\f{\nabla p}{p})\|^2\nonumber\\
&&\leq C_2\Big(\|(\nabla p,\nabla u)\|^2_{L^\infty}+P(\Lambda_{1m}(t))\Big)
\Big(P(\Lambda_{1m}(t))
+\sum_{k=1}^{m-2}\|\partial_t^{k}(\f{\nabla p}{p})\|^2\Big)\nonumber\\
&&\leq C_2\Big(\|(\nabla p,\nabla u)\|^2_{L^\infty}+P(\Lambda_{1m}(t))\Big)
\Big(P(\Lambda_{1m}(t))
+P(\Lambda_{1m}(t))\sum_{k=1}^{m-2}\|\partial_t^{k}(\f{1}{p})\|^2_{L^\infty}
\Big)\nonumber\\
&&\leq C_2\Big(\|(\nabla p,\nabla u)\|^2_{L^\infty}+P(\Lambda_{1m}(t))\Big)P(\Lambda_{1m}(t)).
\end{eqnarray}

Therefore, \eqref{4.48} follows from \eqref{4.49}, \eqref{4.54}, \eqref{4.57} and \eqref{4.59}, which proves Lemma \ref{lem4.7}.
$\hfill\Box$

\begin{lemma}\label{lem4.8}
For every $m\geq1$, it holds that
\begin{eqnarray}\label{4.60}
\int_0^t\|\nabla\partial_t^{m-1}p(\tau)\|^2d\tau
\leq C\v^2\int_0^t\|\nabla^2\partial_t^{m-1}u(\tau)\|^2d\tau
+C[1+P(Q(t))]\int_0^t\Lambda_m(\tau)d\tau.
\end{eqnarray}

\end{lemma}

\noindent\textbf{Proof}.  Applying $\partial_t^{m-1}$ to \eqref{2.5} yields that
\begin{equation*}\label{2.5-10}
\nabla \partial_t^{m-1}p=-\partial_t^{m-1}(\r u_t-\r u\cdot\nabla u-\v \nabla\times\omega+\v \nabla div u).
\end{equation*}
Using Proposition \ref{prop3.2}, one obtains \eqref{4.60}, which thus proves lemma \ref{lem4.8}.
$\hfill\Box$

\

In general, it is hard to obtain the uniform estimate for the term $\int_0^t\|\nabla\partial_t^{m-2}\mbox{div}u\|^2d\tau$ since it may involve two spacial   derivatives in the normal direction. However, $\mbox{div}u$ can be expressed by the variation of the density which is expected to have no strong boundary layers.

\begin{lemma}\label{lem4.5}
For every $m\geq2$, it holds that
\begin{eqnarray}\label{4.43}
\int_0^t\|\nabla\MZ^{m-2}\mbox{div}u(\tau)\|^2d\tau
\leq C\int_0^t\|\nabla\partial_t^{m-1}p(\tau)\|^2d\tau
+C_m[1+P(Q(t))]\int_0^t\Lambda_m(\tau)d\tau.
\end{eqnarray}

\end{lemma}

\noindent\textbf{Proof}.  Applying $\nabla\MZ^{\a}$ to \eqref{2.35} with $|\a|\leq m-2$ gives that
\begin{equation}\label{4.44}
\nabla\MZ^{\a}\mbox{div}u=-\nabla\MZ^{\a}(\ln\r)_t-\nabla\MZ^{\a}(u_i\partial_{y^i}\ln\r)
-\nabla\MZ^{\a}(u\cdot N\partial_z\ln\r)
\end{equation}
It follows from Proposition \ref{prop3.2} that
\begin{eqnarray}\label{4.45}
\int_0^t\|\nabla\MZ^{\a}(\ln\r)_t\|^2d\tau\leq C\int_0^t\|\nabla\partial_t^{m-1}p\|^2d\tau
+C[1+P(Q(t))]\int_0^t\Lambda_m(\tau)d\tau,
\end{eqnarray}
and
\begin{eqnarray}\label{4.46}
\int_0^t\|\nabla\MZ^{\a}(u_i\partial_{y^i}\ln\r)\|^2d\tau\leq C[1+P(Q(t))]\int_0^t\Lambda_m(\tau)d\tau.
\end{eqnarray}
Finally, by using  Proposition \ref{prop3.2} and the hardy inequality \eqref{4.21-7}, one has that
\begin{eqnarray}\label{4.47}
&&\int_0^t\|\nabla\MZ^{\a}(u\cdot N\partial_z\ln\r)\|^2d\tau\leq
\int_0^t\|\nabla(u\cdot N\partial_z\ln\r)\|^2_{\mathcal{H}^{m-2}}d\tau\nonumber\\
&&\leq \int_0^t\|\nabla(u\cdot N)\cdot\partial_z\ln\r)\|^2_{\mathcal{H}^{m-2}}d\tau
+\int_0^t\|(\frac{u\cdot N}{\varphi(z)})\cdot Z_3\nabla\ln\r)\|^2_{\mathcal{H}^{m-2}}d\tau\nonumber\\
&&\leq  C[1+P(Q(t))+\sup_{0\leq\tau\leq t}\|\frac{u\cdot N}{\varphi(z)}\|^2_{L^\infty}]
\int_0^t\Big(\Lambda_m(\tau)+\|\frac{u\cdot N}{\varphi(z)}\|^2_{\mathcal{H}^{m-2}}\Big)d\tau\nonumber\\
&&\leq  C_m[1+P(Q(t))]\int_0^t\Lambda_m(\tau)d\tau.
\end{eqnarray}
Therefore, collecting all the estimates,  one obtains \eqref{4.43}. Thus  the proof of Lemma \ref{lem4.5} is completed.
$\hfill\Box$

\

It follows from  \eqref{4.43} and \eqref{4.60} that
\begin{eqnarray}\label{3.80-1}
&&\d\int_0^t\|\nabla\MZ^{m-2}\mbox{div}u(\tau)\|^2+\|\nabla\partial_t^{m-1}p(\tau)\|^2d\tau+\v\int_0^t\|\nabla\partial_t^{m-1}p(\tau)\|^2d\tau\nonumber \\
&&\leq C\d\int_0^t\|\nabla\partial_t^{m-1}p(\tau)\|^2d\tau+C\int_0^t\v\|\nabla\partial_t^{m-1}p(\tau)\|^2d\tau
+C_m[1+P(Q(t))]\int_0^t\Lambda_m(\tau)d\tau\nonumber\\
&&\leq C\d\v^2\int_0^t\|\nabla^2\partial_t^{m-1}u(\tau)\|^2d\tau
+C[1+P(Q(t))]\int_0^t\Lambda_m(\tau)d\tau.
\end{eqnarray}
Substituting \eqref{3.80-1} into \eqref{4.6} yields that
\begin{eqnarray}\label{4.61}
&&\sup_{0\leq \tau\leq t}\sum_{k=0}^{m-2}\|(\partial_t^{k}\nabla p,\partial_t^{k}\mbox{div}u)(\tau)\|^2_{m-1-k}
+\v\int_0^t\sum_{k=0}^{m-2}\|\partial_t^{k}\nabla\mbox{div}u(\tau)\|^2_{m-1-k}d\tau\nonumber\\
&&\leq C\int\r_0|\MZ^\a\mbox{div}u_0|^2+\f{1}{\g p_0}|\MZ^\a\nabla p_0|^2dx+ CC_{m+2}\Big\{\d\v^2\int_0^t\|\nabla^2\partial_t^{m-1}u(\tau)\|^2d\tau\nonumber\\
&&~~~~~+\v\int_0^t\|\nabla^2u(\tau)\|^2_{\mathcal{H}^{m-2}}d\tau+C_\d[1+P(Q(t))]\int_0^tP(\Lambda_m(\tau)) d\tau\Big\}.
\end{eqnarray}
Then, combining \eqref{4.61} with \eqref{4.25} shows that
\begin{eqnarray}\label{4.62}
&&\sup_{0\leq \tau\leq t}\left\{\sum_{k=0}^{m-2}\|(\partial_t^{k}\nabla p,\partial_t^{k}\mbox{div}u)(\tau)\|^2_{m-1-k}+\v\|(\partial_t^{m-1}\mbox{div}u,\partial_t^{m-1}\nabla p)(\tau)\|^2\right\}\nonumber\\
&&~~~~~~~~+\v\int_0^t\sum_{k=0}^{m-2}\|\partial_t^{k}\nabla\mbox{div}u(\tau)\|^2_{m-1-k}d\tau
+\v^2\int_0^t\|\partial_t^{m-1}\nabla\mbox{div}u(\tau)\|^2d\tau\nonumber\\
&&\leq CC_{m+2}\Big\{\Lambda_m(0)+C\d\v^2\int_0^t\|\nabla^2\partial_t^{m-1}u(\tau)\|^2d\tau
+\v\int_0^t\|\nabla^2u(\tau)\|^2_{\mathcal{H}^{m-2}}d\tau\nonumber\\
&&~~~~~~~~~~~~~~~~~~+[1+P(Q(t))]\int_0^tP(\Lambda_m(\tau)) d\tau\Big\}.
\end{eqnarray}

\subsection{Normal Derivatives Estimates}
\setcounter{equation}{0}

In order to estimate $\|\nabla u\|_{\mathcal{H}^{m-1}}$, it remains  to estimate $\|\chi\partial_n u\|_{\mathcal{H}^{m-1}}$, where $\chi$ is compactly supported in one of the $\Om_j$ and with value one in a neighborhood of the boundary. Indeed, it follows from the definition of the norm that
$\|\chi\partial_{y^i}u\|_{\mathcal{H}^{m-1}}\leq C\| u\|_{\mathcal{H}^{m}}$ for $i=1,2$.  So it suffices estimate $\|\chi\partial_n u\|_{\mathcal{H}^{m-1}}$.

Note that
\begin{equation}\label{5.1}
\mbox{div}u=\partial_nu\cdot n+(\Pi\partial_{y^1}u)_1+(\Pi\partial_{y^2}u)_2.
\end{equation}
and
\begin{equation}\label{5.2}
\partial_nu=[\partial_nu\cdot n] n+\Pi(\partial_nu)
\end{equation}
Thus it follows from \eqref{5.1} and \eqref{5.2} that
\begin{eqnarray*}\label{5.3}
\|\chi\partial_n u\|_{\mathcal{H}^{m-1}}
&&\leq \|\chi\partial_n u\cdot n\|_{\mathcal{H}^{m-1}}
+\|\chi\Pi(\partial_nu)\|_{\mathcal{H}^{m-1}}\nonumber\\
&&\leq  C_{m}\Big\{\|\chi\mbox{div}u\|_{\mathcal{H}^{m-1}}
+\|\chi\Pi(\partial_nu)\|_{\mathcal{H}^{m-1}}+ \|u\|_{\mathcal{H}^{m}}\Big\}.
\end{eqnarray*}
Thus it suffices to estimate $\|\chi\Pi(\partial_nu)\|_{\mathcal{H}^{m-1}}$, since  $ \|u\|_{\mathcal{H}^{m}} $ and $\|\chi\mbox{div}u\|_{\mathcal{H}^{m-1}} $ have been estimated in subsection 3.1 and subsection 3.2, respectively. We extend the smooth symmetric matrix $A$   in \eqref{1.10} to be
$$A(y,z)=A(y).$$
Define
\begin{equation}\label{5.4}
\eta\triangleq \chi\Big(\omega\times n+\Pi(Bu)\Big)=\chi\Big(\Pi(\omega\times n)+\Pi(Bu)\Big).
\end{equation}
The $\eta$ defined here, which enable us to avoid the term $\nabla^2p$, is slightly different from the one in \cite{Masmoudi-R}. Then in view of the Navier-slip boundary condition \eqref{2.7},  $\eta$ satisfies:
\begin{equation}\label{5.5}
\eta|_{\partial\Om}=0.
\end{equation}
Since $\omega\times n=(\nabla u-(\nabla u)^t)\cdot n$, so $\eta$ can be rewritten as
\begin{equation}\label{5.6}
\eta=\chi\Big\{\Pi(\partial_nu)-\Pi(\nabla(u\cdot n))+\Pi((\nabla n)^t\cdot u)+\Pi(Bu)\Big\},
\end{equation}
which  yields immediately that
\begin{equation}\label{5.7}
\|\chi\Pi(\partial_n u)\|_{\mathcal{H}^{m-1}}\leq C_{m+1}(\|\eta\|_{\mathcal{H}^{m-1}}+\|u\|_{\mathcal{H}^{m}}).
\end{equation}
Hence, it remains  to estimate $\|\eta\|_{\mathcal{H}^{m-1}}$. In fact, one can get the following conormal estimates for $\eta$:

\begin{lemma}\label{lem5.1}
For every $m\geq1$, it holds that
\begin{eqnarray}\label{5.8}
&&\sup_{0\leq \tau\leq t}\|\eta(\tau)\|^2_{\mathcal{H}^{m-1}}+\v\int_0^t\|\nabla\eta(\tau)\|^2_{\mathcal{H}^{m-1}}d\tau
\leq CC_{m+2}(\|u(0)\|^2_{\mathcal{H}^{m-1}}+\|\nabla u(0)\|^2_{\mathcal{H}^{m-1}})\nonumber\\
&&~~~+CC_{m+2}\Big\{\d\v^2\int_0^t\|\nabla^2u(\tau)\|^2_{\mathcal{H}^{m-1}}d\tau +C_\d [1+P(Q(t))]\int_0^tP(\Lambda_m(\tau))d\tau\Big\}.
\end{eqnarray}

\end{lemma}

\noindent\textbf{Proof}.  Notice that
$$\nabla\times((u\cdot\nabla) u)=(u\cdot\nabla)\omega-(\omega\cdot\nabla)u+\mbox{div}u \cdot\omega,$$
so $\omega$ satisfies the following  equations
\begin{eqnarray}\label{5.9}
\r \omega_t+\r(u\cdot\nabla)\omega=\v \Delta\omega+F_1,
\end{eqnarray}
with
\begin{equation}\label{5.10}
F_1\triangleq -\nabla\r\times u_t-\nabla\r\times(u\cdot\nabla)u+\r(\omega\cdot\nabla)u-\r\mbox{div}u \omega.
\end{equation}
Consequently, the system of  $\eta$ is
\begin{eqnarray}\label{5.11}
&&\r\eta_t+\r u_1\partial_{y^1}\eta+\r u_2\partial_{y^2}\eta+\r u\cdot N\partial_z\eta-\v\Delta\eta\nonumber\\
&&~~~~~~=\chi[F_1\times n+\Pi(B F_2)]+\chi F_3+F_4+\v\chi\Delta(\Pi B)\cdot u =:F,
\end{eqnarray}
where
\begin{eqnarray}
&&F_2=\v\nabla\mbox{div}u-\nabla p, \label{5.12}\\
&&F_3=-2\sum_{j=1}^2\v\partial_j\omega\times\partial_jn-\v\omega\times\Delta n+\sum_{i=1}^2\r u_i\omega\times\partial_{y^i}n \nonumber\\
&&~~~~~~+ \sum_{i=1}^2\r u_i\partial_{y^i}(\Pi B) u-2\v\sum_{j=1}^2 \partial_j(\Pi B)\partial_ju,\label{5.13}\\
&&F_4=\sum_{i=1}^2\r u_i\partial_{y^i}\chi\cdot (\omega\times n+\Pi(Bu))+\r u\cdot N\partial_z\chi\cdot (\omega\times n+\Pi(Bu))\nonumber\\
&&~~~~~~~~~~-2\sum_{j=1}^3\v\partial_j\chi\partial_j(\omega\times n+\Pi(Bu))-\v\Delta\chi\cdot (\omega\times n+\Pi(Bu))
\label{5.15}.
\end{eqnarray}

We start with $m=1$. Multiplying \eqref{5.11} by $\eta$ and integrating lead to that
\begin{equation}\label{5.16}
\sup_{0\leq \tau\leq t}\int\r|\eta|^2dx+2\v\int_0^t\|\nabla\eta\|^2d\tau
\leq \int\r_0|\eta_0|^2dx+\int_0^t\int F\eta dxd\tau.
\end{equation}
To handle the right-hand side, one notes that
\begin{equation}\label{5.17}
\int_0^t\|\chi(F_1\times n)\|^2_{\mathcal{H}^{m-1}}d\tau\leq C_m[1+P(Q(t))]
\int_0^t\Big(\|\nabla\partial_t^{m-1}p\|^2+\Lambda_m\Big)d\tau,
\end{equation}
\begin{equation}\label{5.18}
\int_0^t\|\chi\Pi(BF_2)\|^2_{\mathcal{H}^{m-1}}d\tau\leq
C_{m+1}\Big\{\int_0^t\Lambda_m+\|\nabla\partial_t^{m-1}p\|^2d\tau
+\v^2\int_0^t\|\chi\nabla\mbox{div}u\|^2_{\mathcal{H}^{m-1}}d\tau\Big\},
\end{equation}
and
\begin{equation}\label{5.19}
\int_0^t\|\chi F_3\|^2_{\mathcal{H}^{m-1}}d\tau\leq
C_{m+2}\Big\{\v^2\int_0^t\|\chi\nabla^2u\|^2_{\mathcal{H}^{m-1}}d\tau
+[1+P(Q(t))]\int_0^t \Lambda_m(\tau)d\tau\Big\},
\end{equation}
Since all the terms in $F_4$ are supported away from the boundary, one can estimate all the derivatives by the $\|\cdot\|_{\mathcal{H}^{m}}$ norms. Therefore, it is easy to obtain
\begin{equation}\label{5.21}
\int_0^t\| F_4\|^2_{\mathcal{H}^{m-1}}d\tau\leq C_{m+1}\Big\{\v^2\int_0^t\|\chi\nabla^2 u\|^2_{\mathcal{H}^{m-1}}d\tau+ C[1+P(Q(t))]\int_0^t \Lambda_m(\tau)d\tau\Big\}.
\end{equation}
Finally, by integrating by parts, it is easy to obtain that
\begin{equation}\label{5.21-1}
\int_0^t\int\v\MZ^{m-1}(\Delta(\Pi B)\cdot u)\cdot \MZ^{m-1}\eta d\tau\leq \d\v\int_0^t \|\nabla\eta\|^2_{\mathcal{H}^{m-1}}d\tau +C_{m+2}\int_0^t \Lambda_m(\tau)d\tau.
\end{equation}
Consequently, substituting these estimates into \eqref{5.16} and using the Cauchy inequality,  one has that
\begin{eqnarray}\label{5.22}
&&\sup_{0\leq \tau\leq t}\int\r|\eta|^2dx+2\v\int_0^t\|\nabla\eta\|^2d\tau \leq \int\r_0|\eta_0|^2dx+C_{m+2}\Big\{\d\int_0^t\|\nabla\partial_t^{m-1}p(\tau)\|^2d\tau\nonumber\\
&&~~~~~~~~~~~~~~~~ +C_\d[1+P(Q(t))]\int_0^t \Lambda_m(\tau)d\tau+\d\v^2\int_0^t\|\chi \nabla^2 u\|^2_{\mathcal{H}^{m-1}}d\tau\Big\}.\nonumber
\end{eqnarray}
Thus,  \eqref{5.8} is proved for $m=1$ by using Lemma \ref{lem4.8}.

\

To prove the general case, we assume that \eqref{5.8} is proved for $k\leq m-2$.  Applying $\MZ^\a$
to \eqref{5.11} for $|\a|=m-1$ yields  that
\begin{equation}\label{5.23}
\r \MZ^\a\eta_t+\r (u\cdot\nabla )\MZ^\a\eta-\v\MZ^\a\Delta\eta= \MZ^\a F+\mathcal{C}_3^\a+\mathcal{C}_4^\a,
\end{equation}
with
\begin{eqnarray}
&&\mathcal{C}_3^\a=-[\MZ^\a,\r]\eta_t=\sum_{|\b|\geq 1, \b+\g=\a}C_{\b,\g}\MZ^\b\r \MZ^\g\eta_t,\label{5.24}\\
&&\mathcal{C}_4^\a=-\sum_{|\b|\geq 1, \b+\g=\a}\sum_{i=1}^2C_{\b,\g}\MZ^\b(\r u_i)\MZ^\g\partial_{y^i}\eta
-\sum_{|\b|\geq 1, \b+\g=\a}C_{\b,\g}\MZ^\b(\r u\cdot N)\MZ^\g\partial_{z}\eta\nonumber\\
&&~~~~~~~~~~~~~~~~~~~~~-\r (u\cdot N)\sum_{|\b|\leq m-2}C_{\b} \partial_z\MZ^\b\eta, \label{5.25}
\end{eqnarray}
where $C_\b$ and $C_{\b,\g}$ are functions depending only on $z$. Multiplying \eqref{5.23} by $\MZ^\a\eta$ and using \eqref{5.17}-\eqref{5.21-1}, one obtains that
\begin{eqnarray}\label{5.26}
&&\sup_{0\leq \tau\leq t}\int\f12\r|\MZ^\a\eta|^2dx
\leq\v\int_0^t\int\MZ^\a\Delta\eta\MZ^\a\eta dxd\tau +\int\f12\r_0|\MZ^\a\eta_0|^2dx+\int_0^t\int(\mathcal{C}_3^\a+\mathcal{C}_4^\a)\MZ^\a\eta dxd\tau\nonumber\\
&&~~ +C_\d C_{m+2}\Big\{[1+P(Q(t))]\int_0^t \Lambda_md\tau+\d\v^2\int_0^t\|\chi \nabla^2 u\|^2_{\mathcal{H}^{m-1}}d\tau+\d\int_0^t\|\nabla\partial_t^{m-1}p\|^2d\tau\Big\}.
\end{eqnarray}

In the local basis, it holds that
\begin{equation*}\label{5.27}
\partial_j=\b^1_j\partial_{y^1}+\b^2_j\partial_{y^2}+\b^3_j\partial_{z},~~\mbox{for}~~j=1,2,3.
\end{equation*}
Therefore, we have the following commutation expansion
\begin{equation}\label{5.28}
\MZ^\a\Delta\eta=\Delta\MZ^\a\eta+\sum_{|\b|\leq m-2}C_{1\b}\partial_{zz}\MZ^\b\eta+\sum_{|\b|\leq m-1}\left(C_{2\b}\partial_{z}\MZ^\b\eta+C_{3\b} Z_y\MZ^\b\eta\right).
\end{equation}
By using the expansion \eqref{5.28} and the inequalities before \eqref{5.23}, one obtains that
\begin{eqnarray}\label{5.29}
&&\v \int_0^t\int\MZ^\a\Delta\eta\MZ^\a\eta dxd\tau=\v\int_0^t\int\Delta\MZ^\a\eta\cdot\MZ^\a\eta dxd\tau+\sum_{|\b|\leq m-2}\v\int_0^t\int C_{1\b}\partial_{zz}\MZ^\b\eta\cdot\MZ^\a\eta dxd\tau\nonumber\\
&&~~~~~~~~+\sum_{|\b|\leq m-1}\v\int_0^t\int\Big\{ C_{2\b}\partial_{z}\MZ^\b\eta+ C_{3\b} Z_y\MZ^\b\eta\Big\}\MZ^\a\eta dxd\tau\nonumber\\
&&\leq -\f34\v \int_0^t\|\nabla\MZ^\a\eta\|^2d\tau+C\v\int_0^t\|\nabla\eta\|^2_{\mathcal{H}^{m-2}}d\tau
+C_{m+2}\v\int_0^t\Lambda_m(\tau)d\tau.
\end{eqnarray}
Note that there is no boundary term in the integrating by parts since $\MZ^\a\eta$ vanishes on the boundary.

It remains to estimate terms involving $\mathcal{C}_3^\a$ and  $\mathcal{C}_4^\a$. By using Proposition \ref{prop3.2}, it is easy to obtain
\begin{eqnarray}\label{5.30}
&&\int_0^t\|\mathcal{C}_3^\a\|^2d\tau+\sum_{|\b|\geq 1, \b+\g=\a}C_{\b,\g} \int_0^t\|\MZ^\b(\r u_i)\MZ^\g\partial_{y^i}\eta\|^2d\tau\nonumber\\
&&\leq C_{m+2}[1+P(Q(t))]\int_0^t \Lambda_m(\tau) d\tau.
\end{eqnarray}
The remaining terms are more involved, since it is desired to obtain an estimate independent of  $\partial_z\eta$.  First, it is easy to obtain that
\begin{eqnarray}\label{5.31}
&&\sum_{|\b|\leq m-2}\int_0^t\|\r (u\cdot N)C_{\b} \partial_z\MZ^\b\eta\|d\tau
\leq \sum_{|\b|\leq m-2}\int_0^t\|\r (\frac{u\cdot N}{\varphi(z)})C_{\b} Z_3\MZ^\b\eta\|d\tau\nonumber\\
&&\leq C_{m+2}[1+P(Q(t))]\int_0^t \Lambda_m(\tau) d\tau.
\end{eqnarray}
On the other hand, we notice that for $ |\beta|\geq1, \b+\g=\a,$ and $|\a|=m-1$
\begin{eqnarray}\label{5.32}
&& \MZ^\b(\r u\cdot N)\MZ^\g\partial_{z}\eta=\f1{\varphi(z)}\MZ^\b(\r u\cdot N)\cdot\varphi(z)\MZ^\g\partial_{z}\eta \nonumber\\
&&=\sum_{\tilde{\beta}\leq \beta,\tilde{\g}\leq \g}C_{\tilde\beta,\tilde\g}\MZ^{\tilde\b}(\r\f{u\cdot N}{ \varphi(z)})\cdot\MZ^{\tilde\g}(Z_3\eta),\nonumber
\end{eqnarray}
where $|\tilde\beta|+|\tilde\g|\leq m-1$, $|\tilde\g|\leq m-2$ and $C_{\tilde\beta,\tilde\g} $ is some smooth bounded coefficient. Therefore, by  similar arguments as \eqref{4.21-5}, \eqref{4.21-6} and using the hardy inequality \eqref{4.21-7} , one has that
\begin{eqnarray}\label{5.33}
\sum_{|\b|\geq 1, \b+\g=\a}\int_0^t\|C_{\b,\g}\MZ^\b(\r u\cdot N)\MZ^\g\partial_{z}\eta\|^2d\tau
\leq C_{m+1}[1+P(Q(t))]\int_0^t \Lambda_m(\tau) d\tau.
\end{eqnarray}
Then, it follows from  \eqref{5.30}-\eqref{5.33} that
\begin{eqnarray}\label{5.34}
\int_0^t\|(\mathcal{C}_3^\a,\mathcal{C}_4^\a)\|^2d\tau
\leq C_{m+2}[1+P(Q(t))]\int_0^t \Lambda_m(\tau) d\tau.
\end{eqnarray}
Substituting \eqref{5.29} and \eqref{5.34} into \eqref{5.26} and using Lemma \ref{lem4.8}, we obtain that
\begin{eqnarray}\label{5.35}
&&\sup_{0\leq \tau\leq t}\int\f12\r|\MZ^\a\eta|^2dx+ \v \int_0^t\|\nabla\MZ^\a\eta\|^2d\tau
\leq C_{m+2}\Big\{\int\f12\r_0|\MZ^\a\eta_0|^2dx+\v\int_0^t\|\nabla\eta\|^2_{\mathcal{H}^{m-2}}d\tau\nonumber\\
&&~~~~~~~~~~~~~~~~~~~~~~~ +C[1+P(Q(t))]\int_0^t \Lambda_m(\tau)d\tau+\d\v^2\int_0^t\| \nabla^2 u\|^2_{\mathcal{H}^{m-1}}d\tau\Big\}.
\end{eqnarray}
By the induction assumption, one can  eliminate the term $ \v\int_0^t\|\nabla\eta\|^2_{\mathcal{H}^{m-2}}d\tau$. So the proof Lemma \ref{lem5.1} is completed. $\hfill\Box$

\

It follows from \eqref{5.1}-\eqref{5.7} that
\begin{equation}\label{5.36}
\sum_{|\b|\leq m-2}\|\MZ^\b\nabla u\|^2_{H^1_{co}}\leq C_{m+1}\Big(\|u\|^2_{\mathcal{H}^m}+\|\eta\|^2_{\mathcal{H}^{m-1}}
+\sum_{k=0}^{m-2}\|\partial_t^{k}\mbox{div}u\|^2_{m-1-k}\Big),
\end{equation}
\begin{equation}\label{5.37}
\int_0^t\|\nabla^2 u\|^2_{\mathcal{H}^{m-1}}d\tau\leq C_{m+2}\int_0^t\Big(\|\nabla u\|^2_{\mathcal{H}^m}
+\|\nabla\eta\|^2_{\mathcal{H}^{m-1}}+\|\nabla\mbox{div}u\|^2_{\mathcal{H}^{m-1}}
+\Lambda_m \Big)d\tau,
\end{equation}
\begin{eqnarray}\label{5.37-1}
&&\sum_{k=0}^{m-2}\int_0^t\|\nabla^2\partial_t^k u\|^2_{m-1-k}d\tau\leq C_{m+2}\Big\{\int_0^t \|\nabla u\|^2_{\mathcal{H}^m}
+\|\nabla\eta\|^2_{\mathcal{H}^{m-1}}d\tau\nonumber\\
&&~~~~~~~~~~~~~~~~~~~~~~~~~~~~~~~~~~+\sum_{k=0}^{m-2}\int_0^t \|\partial_t^{k}\nabla\mbox{div}u\|^2_{m-1-k}d\tau
+\int_0^t \Lambda_md\tau\Big\},
\end{eqnarray}
and
\begin{equation}\label{5.38}
\v\int_0^t\|\nabla^2\MZ^{m-2}u\|^2d\tau\leq C_{m+1}\v\int_0^t\|\nabla\eta\|^2_{\mathcal{H}^{m-2}}d\tau+ C_{m+1}\int_0^t\Lambda_md\tau,
\end{equation}
where \eqref{4.60}, \eqref{4.43} are used  in the estimate of \eqref{5.38}.
Then, taking $\d$ suitably small and  using \eqref{2.16}, \eqref{4.60},  \eqref{4.60}, \eqref{4.43}, \eqref{4.62}, \eqref{5.8} and \eqref{5.36}-\eqref{5.38}, one has that
\begin{eqnarray}\label{5.39}
&&\sup_{0\leq \tau\leq t}\Big\{\Lambda_{1m}(\tau)+\|\eta(\tau)\|^2_{\mathcal{H}^{m-1}}
+\v\|(\partial_t^{m-1}\mbox{div}u,\partial_t^{m-1}\nabla p)(\tau)\|^2\Big\}
+\v\int_0^t\|\nabla\eta(\tau)\|^2_{\mathcal{H}^{m-1}}d\tau\nonumber\\
&&+\v\int_0^t\|\nabla u(\tau)\|^2_{\mathcal{H}^{m}}d\tau
+\v\int_0^t\sum_{k=0}^{m-2}\|\partial_t^{k}\nabla\mbox{div}u(\tau)\|^2_{m-1-k}d\tau
+\v^2\int_0^t\|\partial_t^{m-1}\nabla\mbox{div}u(\tau)\|^2d\tau\nonumber\\
&&+\v\sum_{k=0}^{m-2}\int_0^t\|\nabla^2\partial_t^k u(\tau)\|^2_{m-1-k}d\tau
+\v^2\int_0^t\|\nabla^2 u(\tau)\|^2_{\mathcal{H}^{m-1}}d\tau+\int_0^t\|\nabla\partial_t^{m-1}p(\tau)\|^2d\tau\nonumber\\
&&\leq CC_{m+2}\Big\{\Lambda_{m}(0) +[1+P(Q(t))]\int_0^tP(\Lambda_m(\tau)) d\tau\Big\}.
\end{eqnarray}

\subsection{$L^\infty$-Estimates}
\setcounter{equation}{0}

To close the estimates, we need to bound the $L^\infty$-norms of $u$ and $p$.  First, one has the following Lemma:

\begin{lemma}\label{lem6.1}
For every $|\a|\geq 0$, it holds that
\begin{eqnarray}
&&\|\MZ^\a(\ln\r, p,u)(t)\|^2_{L^\infty}\leq CP(\Lambda_{1m}(t)),~~\mbox{for}~~m\geq 2+|\a|,\label{6.1}\\[2mm]
&&\|\nabla(\ln\r, p)(t)\|^2_{\mathcal{H}^{1,\infty}}\leq C_3\Big(P(\|\Delta p\|^2_{\mathcal{H}^{1}})+P(\Lambda_{1m}(t))\Big),~\mbox{for}~m\geq 5,\label{6.1-1}\\[2mm]
&&\|\mbox{div}u(t)\|^2_{\mathcal{H}^{1,\infty}}\leq C_3[P(\Lambda_{1m}(t))+P(\|\Delta p\|^2_{\mathcal{H}^{1}})]
,~~\mbox{for}~~m\geq 5,\label{6.4}\\[2mm]
&&\|\nabla\mbox{div}u(t)\|^2_{L^\infty}\leq C_3 P(Q(t)),\label{6.4-1}\\
&&\|\nabla\mbox{div}u(t)\|^2_{\mathcal{H}^{1,\infty}}\leq C_4[1+P(Q(t))]\cdot\Big(C_\d P(\Lambda_{1m}(t))+\d\|\Delta p\|^2_{\mathcal{H}^2}\Big),
~~\mbox{for}~~m\geq6,\label{6.4-2}\\
&&Q(t)\leq C_3\sup_{0\leq \tau\leq t}\Big\{\|\nabla u(\tau)\|^2_{\mathcal{H}^{1,\infty}}+P(\Lambda_{1m}(\tau))+P(\|\Delta p(\tau)\|^2_{\mathcal{H}^1})\Big\},
~\mbox{for}~~m\geq5.\label{6.4-3}
\end{eqnarray}
\end{lemma}

\noindent\textbf{Proof}. The proof of \eqref{6.1} is a consequence of \eqref{3.3} and thus omitted here. In order to prove \eqref{6.1-1}, one notes that
\begin{equation*}\label{6.9}
\partial_{ii}=\partial^2_{y^i}-\partial_{y^i}(\partial_i\psi\partial_z)
-\partial_i\psi\partial_z\partial_{y^i}+(\partial_i\psi)^2\partial^2_z,~~\mbox{for}~~i=1,2,
\end{equation*}
which implies that
\begin{equation}\label{6.10}
\Delta =(1+|\nabla\psi|^2)\partial_{zz}+\sum_{i=1,2}\Big(\partial^2_{y^i}-\partial_{y^i}(\partial_i\psi\partial_z)
-\partial_i\psi\partial_z\partial_{y^i}\Big).
\end{equation}
Since
\begin{eqnarray}\label{7.5-1}
\begin{cases}
\|\Delta\ln\r\|^2_{\mathcal{H}^1}\leq C\Big\{P(\|\Delta p\|^2_{\mathcal{H}^1}) +P(\Lambda_{1m})\Big\}, m\geq3,\\[2mm]
\|\Delta p\|^2_{\mathcal{H}^1}\leq C\Big\{P(\|\Delta\ln\r\|^2_{\mathcal{H}^1})+P(\Lambda_{1m})\Big\},m\geq3,\\[3mm]
\|\Delta\ln\r\|^2_{\mathcal{H}^2}\leq C\Big\{\|\Delta p\|^2_{\mathcal{H}^2}+P(\|\Delta\ln\r\|^2_{\mathcal{H}^1})+P(\Lambda_{1m})\Big\} ,m\geq4\\[2mm]
\|\Delta p\|^2_{\mathcal{H}^2}\leq C\Big\{\|\Delta\ln\r\|^2_{\mathcal{H}^2}+ P(\|\Delta\ln\r\|^2_{\mathcal{H}^1})+P(\Lambda_{1m})\Big\},~~~m\geq4,
\end{cases}
\end{eqnarray}
it then follows from \eqref{7.5-1}, \eqref{3.3} and \eqref{6.10} that for $m\geq 5$
\begin{eqnarray}\label{6.11}
&&\|\nabla(\ln\r, p)\|^2_{L^\infty}
\leq C\left(\|\partial_z\nabla(\ln\r, p)\|_{H^1_{co}}+\|\nabla(\ln\r, p)\|_{H^1_{co}}\right)\|\nabla (\ln\r, p)\|_{H^2_{co}}\nonumber\\
&&\leq C(\|\Delta (\ln\r, p)\|^2_{H^1_{co}}+\|\nabla(\ln\r, p)\|^2_{H^2_{co}})\nonumber\\
&&\leq C[\|\Delta(\ln\r, p)\|^2_{\mathcal{H}^{1}} +\Lambda_{1m}(t)]
\leq C_3[P(\|\Delta p\|^2_{\mathcal{H}^{1}}) +P(\Lambda_{1m}(t))],\nonumber
\end{eqnarray}
and
\begin{eqnarray}\label{6.12}
&&\|\MZ\nabla (\ln\r, p)\|^2_{L^\infty}\leq C\left(\|\partial_z\MZ\nabla(\ln\r, p)\|+\|\MZ\nabla(\ln\r, p)\|\right)\|\MZ\nabla (\ln\r, p)\|_{H^3_{co}}\nonumber\\
&&\leq C[\|\Delta (\ln\r, p)\|^2_{\mathcal{H}^{1}}+P(\Lambda_{1m}(t))]\leq C_3[P(\|\Delta p\|^2_{\mathcal{H}^{1}})+P(\Lambda_{1m}(t))],\nonumber
\end{eqnarray}
which proves \eqref{6.1-1}.

By \eqref{2.35}, \eqref{6.1} and \eqref{6.1-1}, it is easy to obtain, for $m\geq5$, that
\begin{equation*}\label{6.13}
\|\mbox{div}u\|^2_{L^\infty}\leq C (\|\ln\r_t\|^2_{L^\infty}+\|u\|^2_{L^\infty}\cdot\|\nabla\ln\r \|^2_{L^\infty})\leq C_3\Big(P(\Lambda_{1m}(t))+P(\|\Delta p\|^2_{\mathcal{H}^{1}})\Big),
\end{equation*}
and
\begin{equation*}\label{6.13-1}
\|\MZ\mbox{div}u\|^2_{L^\infty}\leq C \Big(\|\MZ(\ln\r)_t\|^2_{L^\infty}
+\|\MZ(u\cdot\nabla\ln\r)\|^2_{L^\infty}\Big)\leq C_3\Big(P(\Lambda_{1m}(t))+P(\|\Delta p\|^2_{\mathcal{H}^{1}})\Big).
\end{equation*}
Thus \eqref{6.4} follows. By using \eqref{2.35}, one obtains \eqref{6.4-1}.

It follows from \eqref{2.35}, \eqref{6.1} and \eqref{3.3} that
\begin{eqnarray}\label{6.13-1}
&&\|\nabla\mbox{div}u\|^2_{\mathcal{H}^{1,\infty}}\leq C [1+P(Q(t))]
\cdot\Big(\|\ln\r\|^2_{\mathcal{H}^{2,\infty}}+\|\nabla p\|^2_{\mathcal{H}^{2,\infty}}\|\frac{1 }{p}\|^2_{\mathcal{H}^{2,\infty}}\Big)\nonumber\\
&&\leq C[1+P(Q(t))]
\cdot\Big(P(\Lambda_{1m}(t))+\|\partial_z\nabla p\|_{\mathcal{H}^{2}}P(\Lambda_{1m}(t))\Big)\nonumber\\
&&\leq C_4[1+P(Q(t))]
\cdot\Big(C_\d P(\Lambda_{1m}(t))+\d\|\Delta p\|^2_{\mathcal{H}^{2}}\Big),~~\mbox{for}~~m\geq 6,
\end{eqnarray}
which gives \eqref{6.4-2}. Finally, \eqref{6.4-3} is an immediately consequence of \eqref{6.1} and \eqref{6.1-1}.  Therefore, Lemma \ref{lem6.1} is proved. $\hfill\Box$

\

The following lemma is devoted to the estimate of $\|\nabla u(t)\|^2_{\mathcal{H}^{1,\infty}}$.
\begin{lemma}\label{lem6.2}
 For $m\geq 6$, it holds that
\begin{eqnarray}\label{6.2}
&&\|\nabla u(t)\|^2_{\mathcal{H}^{1,\infty}} \leq C C_{m+2}\Big\{\|(u_0,\nabla u_0)\|^2_{\mathcal{H}^{1,\infty}}+P(\Lambda_{1m}(t))+P(\|\Delta p(t)\|^2_{\mathcal{H}^{1}})+\v^2t\int_0^t\|\nabla^2u\|^2_{\mathcal{H}^{4}}d\tau \nonumber\\
&&~~~~~~~~~~~~~~~~~~~~~~~~~~~~~~
+ t\int_0^t[1+P(\Lambda_{1m})+P(Q)]\cdot[1+\v^2\|\Delta p\|^2_{\mathcal{H}^{2}}]d\tau\Big\}.
\end{eqnarray}
\end{lemma}

\noindent\textbf{Proof}.
Away from the boundary, it follows easily by the classical isotropic Sobolev embedding theorem that
\begin{equation}\label{6.3}
\|\chi\MZ\nabla u\|^2_{L^\infty}+\|\chi\nabla u\|^2_{L^\infty}\leq C\|u\|^2_{\mathcal{H}^{m}}\leq \Lambda_{1m}(t), ~~\mbox{for}~~m\geq 4,
\end{equation}
where the support of  $\chi$ is away from the boundary. Therefore, by using a partition of unity subordinated to the covering \eqref{2.0}, one needs to estimate only $\|\chi_j\MZ\nabla u\|_{L^\infty}+\|\chi_j\nabla u\|_{L^\infty}$ for $j\geq 1$. For notational convenience, we shall denote $\chi_j$ by $\chi$. Similar to \cite{Masmoudi-R}, we use the local parametrization in the neighborhood of the boundary given by a normal geodesic system in which the Laplacian takes a convenient form. Denote
\begin{equation*}\label{6.5}
\Psi^n(y,z)=\left(\begin{array}{cccc} &y\\&\psi(y)\end{array}\right)
-z n(y)=x,
\end{equation*}
where
\begin{equation*}\label{6.6}
n(y)=\f{1}{\sqrt{1+|\nabla\psi(y)|^2}}\left(\begin{array}{cccc} &\partial_1\psi(y)\\&\partial_1\psi(y)\\&-1\end{array}\right),
\end{equation*}
is the unit outward normal. As before, one can extend $n$ and $\Pi$ in the interior by setting
$$n(\Psi^n(y,z))=n(y),~~\Pi(\Psi^n(y,z))=\Pi(y).$$
Note that $n(y,z)$ and $\Pi(y,z)$ have different definitions from the ones used before. The advantage of this parametrization is that in the associated local basis $(\partial_{y^1},\partial_{y^2}, \partial_z)$ of $\mathbb{R}^3$, it holds that $\partial_z=\partial_n$ and
$$\Big(\partial_{y^i}\Big)\Big|_{\Psi^n(y,z)}\cdot\Big(\partial_{z}\Big)\Big|_{\Psi^n(y,z)}=0.$$
The scalar product on $\mathbb{R}^3$ induces in this coordinate system the Riemannian metric $g$ with
 the form
\begin{equation*}\label{6.7}
g(y,z)=\left(\begin{array}{cccc} &\tilde{g}(y,z)&0\\&0&1 \end{array}\right).
\end{equation*}
Therefore, the Laplacian in this coordinate system reads
\begin{equation}\label{6.8}
\Delta f=\partial_{zz}f+\f12\partial_z(\ln|g|)\partial_zf+\Delta_{\tilde{g}}f
\end{equation}
where $|g|$ denotes the determinant of the matrix $g$, and $\Delta_{\tilde{g}}$  is defined by
$$\Delta_{\tilde{g}}f=\f{1}{\sqrt{|\tilde{g}|}} \sum_{i,j=1,2}\partial_{y^i}(\tilde{g}^{ij}|\tilde{g}|^\f12\partial_{y^j}f),$$
which involves only the tangential derivatives and $\{\tilde{g}^{ij}\}$ is the inverse matrix to $g$.

\

It follows from \eqref{5.1}($n$ and $\Pi$  in the coordinate system we have just defined) and Lemma \ref{lem6.1} that for $m\geq 5$,
\begin{eqnarray}\label{6.17}
&&\|\chi\nabla u\|^2_{L^\infty}+\|\chi\MZ\nabla u\|^2_{L^\infty}\leq C_2 \Big(\|\chi\Pi\partial_nu\|^2_{L^\infty}+\|\MZ(\chi\Pi\partial_nu)\|^2_{L^\infty}+\|\chi\mbox{div}u\|^2_{L^\infty}\nonumber\\
&&~~~~~~~~~~~~~~~~~~~~~~~~~~~~~+\|\MZ\mbox{div}u\|^2_{L^\infty}+\|\MZ Z_yu\|^2_{L^\infty}+\|Z_yu\|^2_{L^\infty}\Big)\nonumber\\
&&\leq C_3\Big\{\|\chi\Pi\partial_nu\|^2_{L^\infty}+\|\MZ(\chi\Pi\partial_nu)\|^2_{L^\infty}
+P(\Lambda_{1m}(t))+P(\|\Delta p\|^2_{\mathcal{H}^{1}})\Big\}.
\end{eqnarray}
Consequently, it suffices to estimate $\|\chi\Pi\partial_nu\|^2_{L^\infty}+\|\MZ(\chi\Pi\partial_nu)\|^2_{L^\infty}$. To this end, it is useful to use the vorticity $\omega$.  Indeed,
\begin{equation*}\label{6.18}
\Pi(\omega\times n)=\Pi((\nabla u-\nabla u^t)\cdot n)=\Pi(\partial_nu-\nabla(u\cdot n)+\nabla n^t\cdot u).
\end{equation*}
Therefore,
\begin{equation}\label{6.19}
\|\chi\Pi\partial_nu\|^2_{L^\infty}+\|\MZ(\chi\Pi\partial_nu)\|^2_{L^\infty}
\leq C_3\Big\{\|\chi\Pi(\omega\times n)\|^2_{L^\infty}+\|\MZ(\chi\Pi(\omega\times n))\|^2_{L^\infty}
+\Lambda_{1m}(t)\Big\},
\end{equation}
which shows that it suffices to estimate $\|\chi\Pi(\omega\times n)\|^2_{L^\infty}$ and $\|\MZ(\chi\Pi(\omega\times n))\|^2_{L^\infty}$.

\

In the support of $\chi$, set
\begin{equation*}\label{6.20}
\tilde{\omega}(y,z)=\omega(\Psi^n(y,z)),~~(\tilde{\r}, \tilde{u})(y,z)=(\r, u)(\Psi^n(y,z)).
\end{equation*}
It follows from \eqref{5.9} and  \eqref{6.8} that
\begin{equation} \label{6.21}
\tilde{\r} \tilde{\omega}_t+\tilde{\r} \tilde{u}^1\partial_{y^1}\tilde{\omega}+\tilde{\r} \tilde{u}^2\partial_{y^2}\tilde{\omega}+\tilde{\r} \tilde{u}\cdot n\partial_{z}\tilde{\omega}=\v(\partial_{zz}\tilde{\omega}+\f12\partial_z(\ln|g|)\partial_z\tilde{\omega}
+\Delta_{\tilde{g}}\tilde{\omega})+\tilde{F}_1,
\end{equation}
and
\begin{equation} \label{6.22}
\tilde{\r} \tilde{u}_t+\tilde{\r} \tilde{u}^1\partial_{y^1}\tilde{u}+\tilde{\r} \tilde{u}^2\partial_{y^2}\tilde{u}+\tilde{\r} \tilde{u}\cdot n\partial_{z}\tilde{u}=\v(\partial_{zz}\tilde{u}+\f12\partial_z(\ln|g|)\partial_z\tilde{u}
+\Delta_{\tilde{g}}\tilde{u})+\tilde{F}_2,
\end{equation}
here
\begin{equation*} \label{6.23}
\tilde{F}_1(y,z)=F_1(\Psi^n(y,z)),~~\tilde{F}_2(y,z)=F_2(\Psi^n(y,z)),
\end{equation*}
where $F_1$ and $F_2$ are defined in \eqref{5.10} and \eqref{5.12}, respectively.  Note that we use the same convention as before for a vector $u$, and $u^j$ denotes the components of $u$ in the local basis $(\partial_{y^1}, \partial_{y^2}, \partial_z)$ just defined in this section, whereas $u_j$ denotes its components in the standard basis of $\mathbb{R}^3$. The vectorial equation of \eqref{6.21} and \eqref{6.22}
have to be understood component by component in the standard basis of $\mathbb{R}^3$.

Similar to \eqref{5.6}, one can define
\begin{equation} \label{6.24}
\tilde{\eta}(y,z)=\chi\Big(\tilde\omega\times n+\Pi(B\tilde{u})\Big),
\end{equation}
where $A$ is extended into the interior domain by $A(y,z)=A(y)$. Thus, \eqref{2.7} implies that
\begin{equation} \label{6.25}
\tilde{\eta}(y,0)=0.
\end{equation}
Due to \eqref{6.21} and \eqref{6.22},  $\tilde\eta$ solves the equations
\begin{eqnarray} \label{6.26}
&&\tilde{\r} \tilde{\eta}_t+\tilde{\r} \tilde{u}^1\partial_{y^1}\tilde{\eta}+\tilde{\r} \tilde{u}^2\partial_{y^2}\tilde{\eta}+\tilde{\r} \tilde{u}\cdot n\partial_{z}\tilde{\eta}\nonumber\\
&& ~~=\v(\partial_{zz}\tilde{\eta}+\f12\partial_z(\ln|g|)\partial_z\tilde{\eta})
+\chi(\tilde{F}_1\times n)+\chi\Pi(B\tilde{F}_2)+F_\chi+\chi F_\kappa,\nonumber
\end{eqnarray}
with
\begin{eqnarray}
&&F_\chi= \Big(\tilde{\r} \tilde{u}^1\partial_{y^1} +\tilde{\r} \tilde{u}^2\partial_{y^2} +\tilde{\r} \tilde{u}\cdot n\partial_{z}\Big)\chi\cdot(\tilde\omega\times n+\Pi(B\tilde{u}))\nonumber\\
&&~~~~~~~~~~-\v\Big(\partial_{zz}\chi+2\partial_z\chi\partial_z+\f12\partial_z(\ln|g|)\cdot\partial_z\chi\Big)
\cdot(\tilde\omega\times n+\Pi(B\tilde{u})), \label{6.27}\nonumber\\
&&F_\kappa= \Big(\tilde{\r} \tilde{u}^1\partial_{y^1}\Pi +\tilde{\r} \tilde{u}^2\partial_{y^2}\Pi\Big) \cdot(B\tilde{u})+ \tilde{\omega}\times\Big(\tilde{\r} \tilde{u}^1\partial_{y^1}n +\tilde{\r} \tilde{u}^2\partial_{y^2}n\Big)\nonumber\\
&&~~~~~~~~~~+\Pi\Big((\tilde{\r} \tilde{u}^1\partial_{y^1} +\tilde{\r} \tilde{u}^2\partial_{y^2} +\tilde{\r} \tilde{u}\cdot n\partial_{z})B\cdot\tilde{u}\Big)+\v\Delta_{\tilde{g}}\tilde{\omega}\times n\nonumber\\
&&~~~~~~~~~~+\v\Pi(B\Delta_{\tilde{g}}\tilde{u}).\label{6.28}\nonumber
\end{eqnarray}
Note that in the derivation of the source terms above, in particular,  $F_\kappa$, which contains all the commutators coming from the fact that $n$ and $\Pi$ are not constant, we have used the fact that in the coordinate system just defined, $n$ and $\Pi$ do not depend on the normal variable. Since $\Delta_{\tilde{g}}$ involves only the tangential derivatives, and the derivatives of $\chi$ are compactly supported away from the boundary, the following estimates hold for $m\geq6$,
\begin{eqnarray}
&&\|\chi(F_1\times n)\|^2_{\mathcal{H}^{1,\infty}}\leq C_2 P(Q(t)), \label{6.29}\\
&&\|F_\chi\|^2_{\mathcal{H}^{1,\infty}}
\leq C_3\Big(\|\r u\|^2_{\mathcal{H}^{1,\infty}}\cdot\|u\|^2_{\mathcal{H}^{2,\infty}}
+\v^2\|u\|^2_{\mathcal{H}^{3,\infty}}\Big) \leq C_3\Big\{ P(Q)+P(\Lambda_{1m}) \Big\},\label{6.30}\\
&&\|\chi F_\kappa\|^2_{\mathcal{H}^{1,\infty}}
\leq C_4\Big\{\|u\|^8_{\mathcal{H}^{1,\infty}}+\|u\|^4_{\mathcal{H}^{1,\infty}}\|\nabla u\|^4_{\mathcal{H}^{1,\infty}}+\|\r\|^4_{\mathcal{H}^{1,\infty}}
+\v^2(\|u\|^2_{\mathcal{H}^{3,\infty}}+\|\nabla u\|^2_{\mathcal{H}^{3,\infty}})\Big\}\nonumber\\
&&~~~~~~~~~~~~~~~\leq C_4\Big\{ P(Q(t))+P(\Lambda_{1m})
+\v^2\|\nabla^2 u\|^2_{\mathcal{H}^{4}}\Big\}, \label{6.31}
\end{eqnarray}
and \eqref{6.4-2} implies that
\begin{eqnarray}\label{6.32}
&&\|\chi\Pi(B\tilde{F}_2)\|^2_{\mathcal{H}^{1,\infty}}\leq C_3\Big\{\v^2\|\nabla\mbox{div}u\|^2_{\mathcal{H}^{1,\infty}}+\|\nabla p\|^2_{\mathcal{H}^{1,\infty}}\Big\}\nonumber\\
&&\leq C_4\Big\{P(Q(t))+P(\Lambda_{1m}(t))+C\v^2[1+P(Q(t))]\cdot\|\Delta p\|^2_{\mathcal{H}^{2}}\Big\}.
\end{eqnarray}
Consequently, it follows from \eqref{6.29}-\eqref{6.32} that for $m\geq 6$
\begin{equation}\label{6.33}
\|\tilde{F}\|^2_{\mathcal{H}^{1,\infty}}\leq C_4\Big\{  P(Q(t))+P(\Lambda_{1m}(t))+\v^2[1+P(Q(t))]\cdot\|\Delta p\|^2_{\mathcal{H}^{2}}+\v^2\|\nabla^2 u\|^2_{\mathcal{H}^{4}}\Big\},
\end{equation}
where $\tilde{F}=\chi(\tilde{F}_1\times n)+\chi\Pi(B\tilde{F}_2)+F_\chi+\chi F_\kappa$.

\

In order to eliminate the term $\f12\partial_z(\ln|g|)\partial_z\tilde{\eta}$, one can define
\begin{equation}\label{6.34}
\tilde\eta=\f{1}{|g|^\f14}\bar{\eta}\doteq\bar{\g} \bar{\eta}.
\end{equation}
Note that
\begin{equation}\label{6.35}
\|\tilde\eta\|_{\mathcal{H}^{1,\infty}} \leq C_3\|\bar{\eta}\|_{\mathcal{H}^{1,\infty}},~~\mbox{and}
~~\|\bar{\eta}\|_{\mathcal{H}^{1,\infty}} \leq C_3\|\tilde\eta\|_{\mathcal{H}^{1,\infty}},
\end{equation}
and $\bar\eta$ solves the equations
\begin{eqnarray} \label{6.36}
&&\tilde{\r} \bar{\eta}_t+\tilde{\r} \tilde{u}^1\partial_{y^1}\bar{\eta}+\tilde{\r} \tilde{u}^2\partial_{y^2}\bar{\eta}+\tilde{\r} \tilde{u}\cdot n\partial_{z}\bar{\eta}-\v\partial_{zz}\bar{\eta}\nonumber\\
&& ~~=\f1{\bar\g}\Big(\tilde{F}+\v\partial_{zz}\bar\g\cdot \bar\eta+\f12\v\partial_z(\ln|g|)\partial_z\g\cdot\bar\eta
-\tilde\r(\tilde{u}\cdot\nabla \bar\g)\bar\eta\Big)=:\mathcal{S}.
\end{eqnarray}

It is difficult to obtain the  explicit solution formula for \eqref{6.36} directly, so one rewrites it as
\begin{eqnarray} \label{6.37}
&&\tilde{\r}(t,y,0)\Big[\bar{\eta}_t+ \tilde{u}^1(t,y,0)\partial_{y^1}\bar{\eta}+ \tilde{u}^2(t,y,0)\partial_{y^2}\bar{\eta}+z\partial_z(\tilde{u}\cdot n)(t,y,0)\partial_{z}\bar{\eta}\Big]-\v\partial_{zz}\bar{\eta}\nonumber\\
&&=\mathcal{S}+[\tilde{\r}(t,y,0)-\tilde{\r}(t,y,z)]\bar\eta_t
+\sum_{i=1,2}[(\tilde{\r}\tilde{u}^i)(t,y,0)-(\tilde{\r}\tilde{u}^i)(t,y,z)]\partial_{y^i}\bar\eta\nonumber\\
&&~~~~~~-\tilde{\r}(t,y,z)[(\tilde{u}\cdot n)(t,y,z)-z\partial_z(\tilde{u}\cdot n)(t,y,0)]\partial_z\bar\eta\nonumber\\[2mm]
&&~~~~~~-[\tilde{\r}(t,y,z)-\tilde{\r}(t,y,0)]\cdot z\partial_z(\tilde{u}\cdot n)(t,y,0)\partial_z\bar\eta=:G
~~\mbox{for}~~z>0,
\end{eqnarray}
with the boundary condition $\bar\eta(t,y,0)=0$. By Lemma \ref{lemA.2} in Appendix, one has that
\begin{eqnarray}\label{6.38}
&&\|\bar\eta\|_{\mathcal{H}^{1,\infty}} \lesssim \|\bar\eta_0\|_{\mathcal{H}^{1,\infty}}
+\int_0^t\|\tilde{\r}^{-1}\|_{L^\infty} \|G\|_{\mathcal{H}^{1,\infty}}d\tau\nonumber\\
&&~~~~~~~~~~~~~~~~~~~~~~~~~~~+\int_0^t(1+\|\tilde{\r}^{-1}\|_{L^\infty})(1+\|(\r,u,\nabla u)\|^2_{\mathcal{H}^{1\infty}}) \|\bar\eta\|_{\mathcal{H}^{1,\infty}}d\tau\nonumber\\
&&\lesssim \|\bar\eta_0\|_{\mathcal{H}^{1,\infty}}
+C\int_0^t\|G\|_{\mathcal{H}^{1,\infty}}d\tau+C\int_0^t(1+P(\Lambda_{1m})+\|\MZ\nabla u\|^2_{L^\infty}) \|\bar\eta\|_{\mathcal{H}^{1,\infty}}d\tau.
\end{eqnarray}
It remains to estimate the right hand side of \eqref{6.38}. First, by \eqref{6.33}, one has that
\begin{eqnarray}\label{6.39}
&&\|\mathcal{S}\|^2_{\mathcal{H}^{1,\infty}}\leq C_4\Big\{ C_\d[P(Q(t))+P(\Lambda_{1m}(t))]+\v^2[1+P(Q(t))]\cdot\|\Delta p\|^2_{\mathcal{H}^{2}}\nonumber\\
&&~~~~~~~+[1+P(\Lambda_{1m}(t))]\cdot\|\bar\eta\|^2_{\mathcal{H}^{1,\infty}}+\d\v^2\|\nabla^2 u\|^2_{\mathcal{H}^{4}}\Big\},~\mbox{for}~m\geq 6.
\end{eqnarray}
Next, by the Taylor formula and the fact that $\bar\eta$ is compactly supported in $z$, one can obtain that
\begin{eqnarray}\label{6.40}
&&\|[\tilde{\r}(t,y,0)-\tilde{\r}(t,y,z)]\bar\eta_t\|^2_{\mathcal{H}^{1,\infty}}\leq C \|\bar\eta\|^2_{\mathcal{H}^{1,\infty}}
+C\|\MZ[\tilde{\r}(t,y,0)-\tilde{\r}(t,y,z)]\|^2_{L^\infty}\cdot\|\bar\eta_t\|^2_{L^\infty}\nonumber\\
&&~~~~~~~~~~~~~~~~~~~~~~~~~~~~~~~~~~~~~~~~~~~
+\|[\tilde{\r}(t,y,0)-\tilde{\r}(t,y,z)]\cdot\MZ\bar\eta_t\|^2_{L^\infty}\nonumber\\
&&\leq C\|\bar\eta\|^2_{\mathcal{H}^{1,\infty}}
+C\|\MZ\r\|^2_{L^\infty}\cdot\|\bar\eta_t\|^2_{L^\infty}
+C\|\nabla \r\|^2_{L^\infty}\cdot\|\varphi(z)\MZ\bar\eta_t\|^2_{L^\infty}.
\end{eqnarray}
By \eqref{3.3}, one has the following inequality, for $|\a|\geq 0$
\begin{eqnarray}\label{6.41}
&&\|\varphi(z)\MZ^\a \bar\eta\|^2_{L^\infty}\leq C\left(\|\nabla(\varphi(z)\MZ^\a \bar\eta)\|_{H^1_{co}}+\|\varphi(z)\MZ^\a \bar\eta\|_{H^1_{co}}\right)\|\varphi(z)\MZ^\a \bar\eta\|_{H^2_{co}}\nonumber\\
&&~~~~~~~~~~~~~~~~~~~\leq C\|\MZ^\a \bar\eta\|^2_{H^2_{co}}.
\end{eqnarray}
Therefore, substituting \eqref{6.41} with $|\a|=2$ into \eqref{6.40} shows that
\begin{eqnarray}\label{6.42}
&&\|[\tilde{\r}(t,y,0)-\tilde{\r}(t,y,z)]\bar\eta_t\|^2_{\mathcal{H}^{1,\infty}}\leq C_4[1+P(Q(t))]\cdot(\|\bar\eta\|^2_{\mathcal{H}^{1,\infty}}+\|\MZ^2 \bar\eta\|^2_{H^2_{co}})\nonumber\\
&&\leq C[1+P(Q(t))]\cdot(\|\bar\eta\|^2_{\mathcal{H}^{1,\infty}}+P(\Lambda_{1m}(t)),~~\mbox{for}~~m\geq 5.
\end{eqnarray}
Similarly, one has that for $m\geq 5$
\begin{eqnarray}\label{6.43}
&&\|[(\tilde{\r}\tilde{u}^1)(t,y,0)-(\tilde{\r}\tilde{u}^1)(t,y,z)]\partial_{y^1}\bar\eta\|^2_{\mathcal{H}^{1,\infty}}
+\|[(\tilde{\r}\tilde{u}^2)(t,y,0)-(\tilde{\r}\tilde{u}^2)(t,y,z)]\partial_{y^2}\bar\eta\|^2_{\mathcal{H}^{1,\infty}}\nonumber\\
&&\leq C_4[1+P(Q(t))]\cdot(\|\bar\eta\|^2_{\mathcal{H}^{1,\infty}}+P(\Lambda_{1m}(t))),
\end{eqnarray}
\begin{eqnarray} \label{6.44}
&&\|[\tilde{\r}(t,y,z)-\tilde{\r}(t,y,0)]\cdot z\partial_z(\tilde{u}\cdot n)(t,y,0)\partial_z\bar\eta\|^2_{\mathcal{H}^{1,\infty}}\nonumber\\
&&\leq C\|\nabla\r\|^2_{\mathcal{H}^{1,\infty}}\|\nabla u\|^2_{\mathcal{H}^{1,\infty}}
\| \varphi(z) Z_3\bar\eta\|^2_{\mathcal{H}^{1,\infty}}\leq C\Big\{P(Q(t))+P(\Lambda_{1m}(t))\Big\},
\end{eqnarray}
and
\begin{eqnarray} \label{6.45}
&&\|\tilde{\r}(t,y,z)[(\tilde{u}\cdot n)(t,y,z)-z\partial_z(\tilde{u}\cdot n)(t,y,0)]\partial_z\bar\eta\|^2_{\mathcal{H}^{1,\infty}}\nonumber\\
&&\leq C\|\r\|^2_{\mathcal{H}^{1,\infty}}\cdot\|[(\tilde{u}\cdot n)(t,y,z)-z\partial_z(\tilde{u}\cdot n)(t,y,0)]\partial_z\bar\eta\|^2_{\mathcal{H}^{1,\infty}}\nonumber\\
&&\leq C_4\|\r\|^2_{\mathcal{H}^{1,\infty}}\Big(\|\nabla u\|^2_{L^\infty}\|\bar\eta\|^2_{\mathcal{H}^{1,\infty}}+\|\MZ[(\tilde{u}\cdot n)(t,y,z)-z\partial_z(\tilde{u}\cdot n)(t,y,0)]\cdot\partial_z\bar\eta\|^2_{L^\infty}\nonumber\\
&&~~~~~~~+\|[(\tilde{u}\cdot n)(t,y,z)-z\partial_z(\tilde{u}\cdot n)(t,y,0)]\MZ\partial_z\bar\eta\|^2_{L^\infty}\Big)\nonumber\\
&&\leq C_4\|\r\|^2_{\mathcal{H}^{1,\infty}}\Big(\|\nabla u\|^2_{L^\infty}\|\bar\eta\|^2_{\mathcal{H}^{1,\infty}}
+\|\nabla u\|^2_{\mathcal{H}^{1,\infty}}\|Z_3\bar\eta\|^2_{L^\infty}
+\|\partial_{zz}(\tilde{u}\cdot n)\|^2_{L^\infty}\|\varphi^2(z)\MZ\partial_z\bar\eta\|^2_{L^\infty}\Big)\nonumber\\
&&\leq C_4[1+P(Q(t))+\|\partial_{zz}(\tilde{u}\cdot n)\|^2_{L^\infty}]
\cdot\Big(\|\bar\eta\|^2_{\mathcal{H}^{1,\infty}}+P(\Lambda_{1m}(t)) \Big).
\end{eqnarray}
To complete the proof, one needs to dealt with the term $\|\partial_{zz}(\tilde{u}\cdot n)\|^2_{L^\infty}$ on the right hand side of \eqref{6.45}.  Since  that $n$ is independent of $z$ and $\partial_z=\partial_n$, it follows from \eqref{5.1} that
\begin{equation}\label{6.46}
\partial_{zz}(\tilde{u}\cdot n)=\partial_z(\partial_n\tilde{u}\cdot n)
=\partial_z\mbox{div}u-\partial_z(\Pi\partial_{y^1}u)_1-\partial_z(\Pi\partial_{y^2}u)_2.
\end{equation}
Then, this, together with \eqref{6.46} and \eqref{6.4-1}, shows that
\begin{equation}\label{6.47}
\|\partial_{zz}(\tilde{u}\cdot n)\|_{L^\infty}
\leq \|\partial_z\mbox{div}u\|_{L^\infty}+\|\partial_z(\Pi\partial_{y^1}u)_1\|_{L^\infty}
+\|\partial_z(\Pi\partial_{y^2}u)_2\|_{L^\infty}\leq C_3P(Q(t)).
\end{equation}
Substituting \eqref{6.47} into \eqref{6.45} yields that
\begin{eqnarray} \label{6.48}
&&\|\tilde{\r}(t,y,z)[(\tilde{u}\cdot n)(t,y,z)-z\partial_z(\tilde{u}\cdot n)(t,y,0)]\partial_z\bar\eta\|^2_{\mathcal{H}^{1,\infty}}\nonumber\\
&&\leq C_4[1+P(Q(t))]
\cdot\Big(\|\bar\eta\|^2_{\mathcal{H}^{1,\infty}}+P(\Lambda_{1m}(t)) \Big),~~\mbox{for}~~m\geq 5.
\end{eqnarray}
Combining \eqref{6.39}, \eqref{6.42}-\eqref{6.44} with \eqref{6.48} leads to that for $m\geq 6$,
\begin{equation} \label{6.49}
\|G\|^2_{\mathcal{H}^{1,\infty}}\leq C_4\Big\{ P(Q(t))+P(\Lambda_{1m}(t))+\v^2[1+P(Q(t))]\cdot\|\Delta p\|^2_{\mathcal{H}^{2}}+\v^2\|\nabla^2 u\|^2_{\mathcal{H}^{4}}\Big\}.
\end{equation}
Then, substituting \eqref{6.49} into \eqref{6.38} gives that for $m\geq 6$,
\begin{eqnarray}\label{6.50}
\|\bar\eta\|^2_{\mathcal{H}^{1,\infty}}
&&\lesssim \|\bar\eta_0\|^2_{\mathcal{H}^{1,\infty}}
+ C_4t\int_0^t [1+P(\Lambda_{1m})+P(Q)]d\tau\nonumber\\
&&~~~+C_4t\v^2\int_0^t\Big([1+P(Q(t))]\cdot\|\Delta p\|^2_{\mathcal{H}^{2}}
+ \|\nabla^2 u\|^2_{\mathcal{H}^{4}}\Big)d\tau.
\end{eqnarray}
Then, \eqref{6.2} follows from  \eqref{6.50}, \eqref{6.35}, \eqref{6.19}, \eqref{6.17} and \eqref{6.3}. This completes the proof of Lemma \ref{lem6.2}.
$\hfill\Box$

\subsection{Uniform Estimate for $\Delta p$}
\setcounter{equation}{0}

In order to complete the a priori estimates, we still need to estimate $\Delta p$. Due to \eqref{7.5-1}, it suffices to estimate $\Delta\ln\r$.  Applying $\mbox{div}$ to  \eqref{2.5} yields that
\begin{equation*}\label{7.1}
-2\v\Delta\mbox{div}u+\Delta p=-\mbox{div}(\r \dot{u}), ~~\mbox{with}~~\dot{u}=u_t+(u\cdot\nabla)u.
\end{equation*}
Substituting \eqref{2.35} into the above equations leads to that
\begin{equation}\label{7.2}
2\v\Delta(\ln\r)_t+2\v u\cdot\nabla\Delta\ln\r+\Delta p
=-2\v \Delta u\cdot \nabla\ln\r-4\v\sum_{k=1}^3\partial_ku\cdot\nabla\partial_k\ln\r-\mbox{div}(\r \dot{u}).
\end{equation}

\begin{lemma}\label{lem7.1}
For $m \geq 6$, it holds that
\begin{eqnarray}\label{7.28}
&&\sup_{0\leq\tau\leq t}\Big(\|\Delta p(\tau)\|^2_{\mathcal{H}^1}+\v\|\Delta p(\tau)\|^2_{\mathcal{H}^2}\Big)
+\int_0^t\|\Delta p(\tau)\|^2_{\mathcal{H}^2}d\tau\nonumber\\
&&\leq  CC_{m+2}\Big\{P(\mathcal{N}_m(0))+[1+P(Q(t))]\int_0^tP(\mathcal{N}_m(\tau))d\tau\Big\},~m\geq 6.
\end{eqnarray}

\end{lemma}

\noindent\textbf{Proof}.  Applying $\MZ^\a(|\a|\leq 2)$ to \eqref{7.2} gives that
\begin{eqnarray} \label{7.14}
2\v\MZ^\a\Delta(\ln\r)_t+\MZ^\a\Delta p&&=-2\v \MZ^\a(\Delta u\cdot \nabla\ln\r)-4\v\sum_{k=1}^3\MZ^\a(\partial_ku\cdot\nabla\partial_k\ln\r)-\MZ^\a\mbox{div}(\r \dot{u})\nonumber\\
&&~~~~-2\v \MZ^\a\Big(\sum_{i=1,2}u_i\cdot\partial_{y^i}\Delta\ln\r
+u\cdot N\cdot\partial_{z}\Delta\ln\r\Big).
\end{eqnarray}
Multiplying \eqref{7.14} by $\MZ^2\Delta\ln\r$, one can obtain that
\begin{eqnarray}\label{7.15}
&& \v\|\MZ^\a\Delta\ln\r\|^2+\int_0^t\int\MZ^\a\Delta p\cdot\MZ^\a\Delta\ln\r dxd\tau\leq \v\|\MZ^\a\Delta\ln\r_0\|^2_{L^2}\nonumber\\
&&~~~-2\v\int_0^t\int\MZ^\a(\Delta u\nabla\ln\r) \MZ^\a\Delta\ln\r dxd\tau-4\v\sum_{k=1}^3\int_0^t\int\MZ^\a(\partial_ku\cdot\nabla\partial_k\ln\r) \MZ^\a\Delta\ln\r dxd\tau\nonumber\\
&&~~~-2\v\int_0^t\int\MZ^\a\Big(\sum_{i=1,2}u_i\cdot\partial_{y^i}\Delta\ln\r
+u\cdot N\cdot\partial_{z}\Delta\ln\r\Big)\MZ^\a\Delta\ln\r dxd\tau\nonumber\\
&&~~~~~~~~~~~~~~~~~~~~~~~~~~~-\int_0^t\int \MZ^\a\mbox{div}(\r \dot{u})\MZ^\a\Delta\ln\r dxd\tau.
\end{eqnarray}

Note that
\begin{equation}\label{7.5}
\Delta p=\g p \Delta\ln\r+\g \nabla p\cdot\nabla\ln\r.
\end{equation}
This implies that
\begin{eqnarray}\label{7.16}
&& \int_0^t\int\MZ^\a\Delta p\cdot\MZ^\a\Delta\ln\r dxd\tau\nonumber\\
&&\geq \f\g2p(c_0)\int_0^t\|\MZ^\a\Delta\ln\r\|^2d\tau
-C[1+P(Q(t))]\int_0^tP(\Lambda_m)+\|\Delta\ln\r\|^2_{\mathcal{H}^1}d\tau.
\end{eqnarray}
The terms on  the right hand side \eqref{7.15} can be estimated separately. First,
it follows from $\eqref{1.1}_2$, \eqref{4.43} and $\eqref{2.7}$ that for $m\geq 5$,
\begin{eqnarray}\label{7.13}
&& \v^2\int_0^t\|\Delta u\|^2_{\mathcal{H}^2}d\tau\leq
C[1+Q(t)]\int_0^t P(\Lambda_m(\tau))d\tau.
\end{eqnarray}
Thanks to \eqref{3.2}, \eqref{7.13} and Cauchy inequality, one has that
\begin{eqnarray}\label{7.17}
&&-2\v\int_0^t\int\MZ^\a(\Delta u\nabla\ln\r) \MZ^\a\Delta\ln\r dxd\tau
\leq \d\int_0^t\|\MZ^\a\Delta\ln\r\|^2d\tau+C_\d\v^2\int_0^t\|\MZ^\a(\Delta u\nabla\ln\r)\|^2d\tau\nonumber\\
&&\leq \d\int_0^t\|\MZ^\a\Delta\ln\r\|^2d\tau
+C_\d\v^2\|\Delta u\|^2_{L^\infty}\int_0^t\|\nabla\ln\r\|^2_{\mathcal{H}^2}d\tau
+C_\d\v^2P(Q(t))\int_0^t\|\Delta u\|^2_{\mathcal{H}^2}d\tau\nonumber\\
&&\leq\d\int_0^t\|\MZ^\a\Delta\ln\r\|^2d\tau
+C_\d[1+P(Q(t))]\int_0^tP(\Lambda_m)d\tau,
\end{eqnarray}
where one has used the following estimate,
\begin{eqnarray}\label{7.18}
\v^2\|\Delta u\|^2_{L^\infty}&&\leq \v^2\|\nabla\mbox{div}u\|^2_{L^\infty}+\|\nabla p\|^2_{L^\infty}
+\|\r u_t\|^2_{L^\infty}+\|\r u\cdot\nabla u\|^2_{L^\infty}\nonumber\\[1mm]
&&\leq C[1+P(Q(t))],\nonumber
\end{eqnarray}
here \eqref{6.4-1} has been used.
It follows from \eqref{3.3} and the Cauchy inequality that for $m\geq 5$,
\begin{eqnarray}\label{7.19}
&&-4\v\sum_{k=1}^3\int_0^t\int\MZ^\a(\partial_ku\cdot\nabla\partial_k\ln\r) \MZ^\a\Delta\ln\r dxd\tau\nonumber\\
&&\leq \d\int_0^t\|\MZ^\a\Delta\ln\r\|^2d\tau+C_\d\v^2\sum_{k=1}^3\int_0^t\|\MZ^\a(\partial_ku\cdot\nabla\partial_k\ln\r)\|^2\nonumber\\
&&\leq \d\int_0^t\|\MZ^\a\Delta\ln\r\|^2d\tau+C_\d\v^2[1+P(Q(t))]\int_0^tP(\Lambda_m(t))
+\|\Delta\ln\r\|^2_{\mathcal{H}^2}d\tau\nonumber\\
&&~~~~~~+C_\d\v^2\int_0^t\|\MZ^\a\nabla u\|^2_{L^\infty}\cdot[P(\Lambda_m(t))
+\|\Delta\ln\r\|^2]d\tau\nonumber\\
&& \leq \d\int_0^t\|\MZ^\a\Delta\ln\r\|^2d\tau+C_\d\v^2[1+P(Q(t))]\int_0^tP(\Lambda_m)
+\|\Delta\ln\r\|^2_{\mathcal{H}^2}d\tau\nonumber\\
&&~~~~~~~~~+ \v^2\int_0^t\|\nabla^2u\|^2_{\mathcal{H}^3}d\tau+C_\d\v^2\int_0^t\|\nabla u\|^2_{\mathcal{H}^4}\cdot[\|\Delta\ln\r\|^4+P(\Lambda_m)]d\tau.
\end{eqnarray}
Note that
\begin{eqnarray}\label{7.9}
\mbox{div}(\r\dot{u})&&=\r\mbox{div}u_t+\r(u\cdot\nabla)\mbox{div}u
+\nabla\r\cdot u_t+\nabla(\r u)^t\cdot\nabla u\nonumber\\
&&=\r\mbox{div}u_t+\nabla\r\cdot u_t+\nabla(\r u)^t\cdot\nabla u
+\sum_{i=1,2}\r u_i\partial_{y^i}\mbox{div}u+\r \frac{u\cdot N}{\varphi(z)}Z_3\mbox{div}u.\nonumber
\end{eqnarray}
This, together with \eqref{3.2} and the hardy inequality \eqref{4.21-7}, implies that
\begin{eqnarray}\label{7.10}
\int_0^t\|\mbox{div}(\r\dot{u})\|^2_{\mathcal{H}^2}d\tau\leq C_{m+2}[1+P(Q(t))]\int_0^tP(\Lambda_m)d\tau
~~\mbox{for}~~m\geq 4,
\end{eqnarray}
which yields immediately that
\begin{eqnarray}\label{7.20}
&&-\int_0^t\int \MZ^\a\mbox{div}(\r \dot{u})\MZ^\a\Delta\ln\r dxd\tau\nonumber\\
&&\leq \d\int_0^t\|\MZ^\a\Delta\ln\r\|^2d\tau+C_\d[1+P(Q)]\int_0^t
P(\Lambda_m(t))d\tau.
\end{eqnarray}
Finally, by integrating by parts, one has that for $m\geq 5$
\begin{eqnarray}\label{7.21}
&&-2\v\int_0^t\MZ^\a\Big(\sum_{i=1,2}u_i\cdot\partial_{y^i}\Delta\ln\r
+u\cdot N\cdot\partial_{z}\Delta\ln\r\Big)\MZ^\a\Delta\ln\r dxd\tau\leq \v^2\int_0^t\|\nabla^2u\|^2_{\mathcal{H}^3}d\tau\nonumber\\
&&+ \d\int_0^t\|\MZ^\a\Delta\ln\r\|^2d\tau+C_\d C_2[1+P(Q(t))]\int_0^t
\v^2\|\Delta\ln\r\|^2_{\mathcal{H}^2}+P(\Lambda_m)d\tau.\nonumber\\
\end{eqnarray}

Substituting \eqref{7.16}, \eqref{7.17}, \eqref{7.19}, \eqref{7.20}, \eqref{7.21} into \eqref{7.15} and choosing $\d$ suitably small, one obtains that
\begin{eqnarray}\label{7.3}
&& \v\|\Delta\ln\r(t)\|^2_{\mathcal{H}^2}+ \int_0^t\|\Delta\ln\r\|^2_{\mathcal{H}^2}d\tau
\leq C\v\|\Delta\ln\r_0\|^2_{\mathcal{H}^2}+C\v^2\int_0^t\|\nabla^2u\|^2_{\mathcal{H}^3}d\tau\nonumber\\
&&+C_\d C_{m+2}[1+P(Q(t))]\int_0^t\Big(\v\|\Delta\ln\r\|^2_{\mathcal{H}^2}+
\|\Delta\ln\r\|^4_{\mathcal{H}^1}+P(\Lambda_m)\Big)d\tau.
\end{eqnarray}
On the other hand, it is easy to get that
\begin{align}\label{7.22}
\|\Delta\ln\r(t)\|^2_{\mathcal{H}^1}&\leq \|\Delta\ln\r(0)\|^2_{\mathcal{H}^1} +\int_0^t\|\partial_t\Delta\ln\r(\tau)\|^2_{\mathcal{H}^1}d\tau\nonumber\\
&\leq \|\Delta\ln\r(0)\|^2_{\mathcal{H}^1} +\int_0^t\|\Delta\ln\r(\tau)\|^2_{\mathcal{H}^2}d\tau.
\end{align}
Combining \eqref{7.3}, \eqref{7.22} with \eqref{5.39}, one obtains that
\begin{eqnarray}\label{7.28-1}
&&\sup_{0\leq\tau\leq t}\Big(\|\Delta\ln\r(\tau)\|^2_{\mathcal{H}^1}+\v\|\Delta\ln\r(\tau)\|^2_{\mathcal{H}^2}\Big)
+\int_0^t\|\Delta\ln\r(\tau)\|^2_{\mathcal{H}^2}d\tau\nonumber\\
&&\leq CC_{m+2}\Big\{[1+P(Q(t))]\int_0^t\Big(1+\v\|\Delta\ln\r\|^2_{\mathcal{H}^2}+
\|\Delta\ln\r\|^4_{\mathcal{H}^1}+P(\Lambda_m)\Big)d\tau\nonumber\\
&&~~~~~~~~~~~~~~~+C(\|\Delta\ln\r_0\|^2_{\mathcal{H}^1}+\v\|\Delta\ln\r_0\|^2_{\mathcal{H}^2}+\Lambda_m(0))\Big\},~m\geq 6.
\end{eqnarray}
Therefore, \eqref{7.28} follows from \eqref{7.5-1},\eqref{5.39} and \eqref{7.28-1}, which completes the proof of Lemma \ref{lem7.1}.

\

\subsection{Proof of Theorem \ref{thm3.1}}
\setcounter{equation}{0}

Noting the definition of $\mathcal{N}_m(t)$ in \eqref{3.0-2}, and using  \eqref{6.4-3}, \eqref{6.2}, \eqref{5.39} and \eqref{7.28}, one has, for $m\geq6$, that
\begin{eqnarray}\label{7.30}
Q(t) &&\lesssim \sup_{0\leq \tau\leq t}\Big\{\|\nabla u(\tau)\|^2_{\mathcal{H}^{1,\infty}}+P(\Lambda_{1m}(\tau))+P(\|\Delta p(\tau)\|^2_{\mathcal{H}^1})\Big\}\nonumber\\
&&\lesssim C_{m+2}\Big\{P(\mathcal{N}_m(0))+P(\mathcal{N}_m(t))\cdot\int_0^tP(\mathcal{N}_m(\tau))d\tau\Big\}.
\end{eqnarray}
In order to close the a priori estimate, one still needs to get the uniform estimates for $\|\nabla\partial_t^{m-1}u\|$. Due to \eqref{7.30}, \eqref{5.39} and  Lemma \ref{lem4.7}, it holds, for $m\geq6$, that
\begin{eqnarray}\label{7.31}
\|\nabla\partial_t^{m-1}u(t)\|^2&&\lesssim C_{m+1}\Big\{ \|u(t)\|^2_{\mathcal{H}^m}+\|\eta(t)\|^2_{\mathcal{H}^{m-1}}+\|\mbox{div}\partial_t^{m-1}u(t)\|^2\Big\}\nonumber\\
&&\lesssim C_{m+2}\Big\{P(\Lambda_{1m}(t))+\|\eta(t)\|^2_{\mathcal{H}^{m-1}}+P(Q(t))\Big\}\nonumber\\
&&\lesssim  C_{m+2}\Big\{ P(\mathcal{N}_m(0))+P(\Lambda_{1m}(t))+\|\eta(t)\|^2_{\mathcal{H}^{m-1}}+P(\mathcal{N}_m(t))\int_0^tP(\mathcal{N}_m(\tau))d\tau\Big\}\nonumber\\
&&\lesssim  C_{m+2}\Big\{P(\mathcal{N}_m(0))+P(\mathcal{N}_m(t))\int_0^tP(\mathcal{N}_m(\tau))d\tau\Big\}.
\end{eqnarray}
Then, one obtains from  \eqref{5.39} and \eqref{7.28}, \eqref{7.30} and \eqref{7.31} that
\begin{eqnarray}\label{7.32}
&&\mathcal{N}_m(t)+\v\int_0^t\|\nabla u(\tau)\|^2_{\mathcal{H}^{m-1}}d\tau
+\v\sum_{k=0}^{m-2}\int_0^t\|\nabla^2\partial_t^ku(\tau)\|^2_{m-1-k}d\tau
+\v^2\int_0^t\|\nabla^2\partial_t^{m-1}u(\tau)\|^2d\tau\nonumber\\
&&~~~~~~~~~~+\int_0^t\|\nabla\partial_t^{m-1}p(\tau)\|^2+\|\Delta p(\tau)\|^2_{\mathcal{H}^2} d\tau\nonumber\\
&&\lesssim C_{m+2}\Big\{P(\mathcal{N}_m(0))+P(\mathcal{N}_m(t))\int_0^tP(\mathcal{N}_m(\tau))d\tau\Big\}.
\end{eqnarray}
Therefore, \eqref{3.0-1} is proved. Furthermore, \eqref{1.1} implies that
\begin{eqnarray}\label{7.33}
|\r(x,0)|\exp(-\int_0^t\|\mbox{div}u(\tau)\|_{L^\infty}d\tau)\leq\r(x,t)\leq |\r(x,0)|\exp(\int_0^t\|\mbox{div}u(\tau)\|_{L^\infty}d\tau),\nonumber
\end{eqnarray}
which proves \eqref{3.0-3}. Thus we complete the proof of Theorem \ref{thm3.1}.


\section{Proof of Theorem \ref{thm1.1}: Uniform Regularity}
\renewcommand{\theequation}{\arabic{section}.\arabic{equation}}

\noindent\textbf{Proof of Theorem \ref{thm1.1}}: In this section, we shall indicate how to combine  our a priori estimates to obtain the uniform existence result. Fix $m\geq6$, we consider initial data $(p^\v_0, u^\v_0)\in X^{\v,m}_{NS}$ such that
\begin{eqnarray}\label{8.1}
&&\mathcal{I}_m(0)=\sup_{0< \v\leq 1}\bigg\{ \|(u_0^\v, p_0^\v)\|^2_{\mathcal{H}^m}+\|\nabla u_0^\v\|^2_{\mathcal{H}^{m-1}}+\sum_{k=0}^{m-2}\|\partial_t^k \nabla p_0^\v\|^2_{m-1-k}+\|\Delta p_0^\v\|^2_{\mathcal{H}^1}\nonumber\\
&&~~~~~~~~~~~~~~~~~~~~~~~~~~~~+\|\nabla u_0^\v\|^2_{\mathcal{H}^{1,\infty}}+\v\|\nabla\partial_t^{m-1}p_0^\v\|^2
+\v\|\Delta p_0^\v\|^2_{\mathcal{H}^2}\bigg\}\leq \tilde{C}_0,
\end{eqnarray}
and
\begin{equation}\label{8.3}
0<\f1{C_0}\leq \r_0^{\v}\leq C_0,
\end{equation}
For such initial data, since we are not aware of a local existence result for \eqref{1.1} and \eqref{1.10}(or \eqref{2.7}), we shall prove this by using the energy estimates obtained in previous sections and a classical iteration scheme. Noting $(p^\v_0, u^\v_0)\in X^{\v,m}_{NS}$, we can find a sequence of smooth approximate  initial data $(p_0^{\v,\delta}, u_0^{\v, \delta})\in X^{\v,m}_{NS,ap}$($\delta$ being a regularization parameter), which have enough space regularity so that the time derivatives at the initial time can be defined by Navier-Stokes equations and the boundary  compatibility conditions are satisfied. We construct  approximate solutions, inductively, as follows\\
(1) Define $u^0=u_0^{\v,\d}$, and \\
(2) Assume that $u^{k-1}$ has been defined for $k\geq1$. Let $(\r^k,u^k)$ be the unique solution to the following linearized  initial boundary value problem:
\begin{eqnarray}\label{8.5}
\begin{cases}
\r^k_t+\mbox{div}(\r^k u^{k-1})=0 ~~\mbox{in}~(0,T)\times\Omega,\\
\r^k u^k_t+\r^k u^{k-1}\cdot\nabla{u}^k+\nabla{p}^k=\v\Delta u^k+\v\nabla\mbox{div}u^k~~\mbox{in}~(0,T)\times\Omega,\\
(\r^k, u^k)|_{t=0}=(\r_0^{\v,\delta}, u_0^{\v, \delta}),~~\mbox{with}~\f2{3C_0}\leq \r_0^{\v,\d}\leq \f32C_0,\\
\mbox{with boundary conditions }~~ \eqref{1.10}~ \mbox{or}~ \eqref{2.7}.
\end{cases}
\end{eqnarray}
Since $\r^k$ and $u^k$ are decoupled,  the existence of global unique  smooth solution $(\r^k,u^k)$ of \eqref{8.5} with $0<\r^k<\infty$ can be obtained by using classical methods, for example,  the same argument of Lemma 3 in Cho, Choe and Kim \cite{Kim} and the standard elliptic regularity results as in Agmon-Douglis-Nirenberg \cite{Nirinberg}. An alternative method is to use a similar arguments (modified slightly) in \cite{Hoff}(or \cite{Zaja}) to get the existence of smooth solutions to \eqref{8.5}.

\

Applying the a priori estimates given in Theorem \ref{thm3.1} and by an induction argument, we obtain a uniform time $\tilde{T}_1>0$ and constant $\tilde{C}_3$(independent of $\v$ and  $\d$),  such that it holds for $(\r^k, u^k),~k\geq 1$ that
\begin{eqnarray}\label{8.6}
&&\mathcal{N}_m(p^k,u^k)(t)+\int_{0}^{t}\|\nabla\partial_t^{m-1}p^k\|^2
+\|\Delta p^k\|^2_{\mathcal{H}^2}d\tau
+\v\int_{0}^{t}\|\nabla u^k\|^2_{\mathcal{H}^{m}}d\tau\nonumber\\
&&+\v\sum_{j=0}^{m-2}\int_{0}^{t}\|\nabla^2\partial_t^j u^k\|^2_{m-j-1}d\tau
+\v^2\int_{0}^{t}\|\nabla^2\partial_t^{m-1}u^k\|^2d\tau\leq \tilde{C}_3,~\forall t\in[0,\tilde{T}_1],
\end{eqnarray}
and
\begin{equation}\label{8.7}
\f1{2C_0}\leq\r^k(x,t)\leq 2C_0, ~~\forall t\in[0,\tilde{T}_1],
\end{equation}
where $\tilde{T}_1$ and $\tilde{C}_3$ depend only on $C_0$ and $\mathcal{I}_m(0)$. Based on the above uniform estimates for $(\r^k, u^k)$,  by the same arguments as section 3 in \cite{Kim}, there exists a uniform time $\tilde{T}_2(\leq \tilde{T}_1)$(independent of $\v$ and  $\d$) such that $(\r^k, u^k)$ converges to a limit $(\r^{\v,\d}, u^{\v,\d})$ as $k\rightarrow +\infty$ in the following strong sense:
$$(\r^k, u^k)\rightarrow (\r^{\v,\d}, u^{\v,\d})~~\mbox{in}~~L^\infty(0,\tilde{T}_2; L^2),~~\mbox{and}~~\nabla u^k\rightarrow \nabla u^{\v,\d}~\mbox{in}~ L^2(0,\tilde{T}_2, L^2).$$
It is easy to check $(\r^{\v,\d}, u^{\v,\d})$ is a weak solution to the problem \eqref{1.1}, \eqref{1.10} with initial data $(\r_0^{\v,\delta}, u_0^{\v, \delta})$. Then, by virtue of the lower semi-continuity of norms, we deduce from the uniform bounds \eqref{8.6} and \eqref{8.7} that $(\r^{\v,\d}, u^{\v,\d})$ satisfies the following regularity estimates
\begin{eqnarray}\label{8.8}
&&\mathcal{N}_m(p^{\v,\d}, u^{\v,\d})(t)+\int_{0}^{t}\|\nabla\partial_t^{m-1}p^{\v,\d}\|^2
+\|\Delta p^{\v,\d}\|^2_{\mathcal{H}^2}d\tau
+\v\int_{0}^{t}\|\nabla u^{\v,\d}\|^2_{\mathcal{H}^{m}}d\tau\nonumber\\
&&+\v\sum_{j=0}^{m-2}\int_{0}^{t}\|\nabla^2\partial_t^j u^{\v,\d}\|^2_{m-j-1}d\tau
+\v^2\int_{0}^{t}\|\nabla^2\partial_t^{m-1}u^{\v,\d}\|^2d\tau\leq \tilde{C}_3,~\forall t\in[0,\tilde{T}_2],
\end{eqnarray}
and
\begin{equation}\label{8.9}
\f1{2C_0}\leq\r^{\v,\d}(x,t)\leq 2C_0, ~~\forall t\in[0,\tilde{T}_2],
\end{equation}
Based on the uniform estimates \eqref{8.8} and \eqref{8.9} for $(\r^{\v,\d}, u^{\v,\d})$, we pass the limit $\d\rightarrow0$ to get a strong solution $(\r^\v,u^\v)$ of \eqref{1.1}, \eqref{1.10} with initial data $(\r_0^\v, u_0^\v)$ satisfying  \eqref{8.1} by using a strong compactness arguments. It follows from \eqref{8.8} that $(p^{\v,\d},u^{\v,\d})$ is bounded uniformly in $L^\infty([0,\tilde{T}_2];H^m_{co})$, while $\nabla(p^{\v,\d},u^{\v,\d})$ is bounded uniformly in $L^\infty([0,\tilde{T}_2];H^{m-1}_{co})$, and $\partial_t(p^{\v,\d},u^{\v,\d})$ is bounded uniformly in $L^\infty([0,\tilde{T}_2];H^{m-1}_{co})$. Then, one can obtain by the strong compactness argument(see \cite{Simon}) that $(p^{\v,\d},u^{\v,\d})$ is compact in $\mathcal{C}([0,\tilde{T}_2];H^{m-1}_{co})$. In particular, there exists a sequence $\d_n\rightarrow 0+$ and $(p^\v,u^\v)\in\mathcal{C}([0,\tilde{T}_2];H^{m-1}_{co})$ such that
\begin{equation*}
(p^{\v,\d_n},u^{\v,\d_n})\rightarrow (p^\v,u^\v)~~\mbox{in}~~\mathcal{C}([0,\tilde{T}_2];H^{m-1}_{co})~~\mbox{as}~~\d_n\rightarrow 0+,
\end{equation*}
or equivalently
\begin{equation}\label{9.3}
(\r^{\v,\d_n},u^{\v,\d_n})\rightarrow (\r^\v,u^\v)~~\mbox{in}~~\mathcal{C}([0,\tilde{T}_2];H^{m-1}_{co})~~\mbox{as}~~\d_n\rightarrow 0+.
\end{equation}
Moreover, applying the lower semi-continuity of norms to the bounds \eqref{8.8} and \eqref{8.9}, one obtains the bounds of
\eqref{1.18} and \eqref{1.19} for $(\r^\v,u^\v)$. It follows from  \eqref{1.18}, \eqref{1.19}, \eqref{9.3} and  the anisotropic Sobolev inequality \eqref{3.3} that
\begin{eqnarray}\label{9.4}
&&\sup_{t\in[0,T]}\|(\r^{\v,\d_n}-\r^\v,u^{\v,\d_n}-u^\v)\|^2_{L^\infty}\nonumber\\
&&\leq C \sup_{t\in[0,T]}\Big( \|\nabla(\r^{\v,\d_n}-\r^\v,u^{\v,\d_n}-u^\v)\|_{H^1_{co}}\cdot\|(\r^{\v,\d_n}-\r^\v,u^{\v,\d_n}-u^\v)\|_{H^2_{co}}\Big)\rightarrow 0.
\end{eqnarray}
Then, it is easy to check that $(\r^\v, u^\v)$ is a weak solution of the Navier-Stokes system.  The uniqueness of the solution $(\r^\v,u^\v)$ is easy since we work on functions with Lipschitz regularity. Therefore, the whole family $(\r^{\v,\d},u^{\v,\d})$   converges to $(\r^\v,u^\v)$. Taking $T_0=\tilde{T}_2$ and $\tilde{C}_1=\tilde{C}_3$,  one completes the proof of Theorem \ref{thm1.1}. $\hfill\Box$


\section{Proof of Theorem \ref{thm1.3}:  Flat Boundary Case}

  Due to Theorem \ref{thm1.1},   there exists a time $T_0>0$ and $\tilde{C}_1>0$  independent of  $\v\in(0,1]$, such that there exists a unique solution $(\r^\v, u^\v)$ of \eqref{1.1},  \eqref{1.11}, \eqref{10.1} which is defined on $[0,T_0]$  and satisfies the uniform estimates \eqref{1.18} and \eqref{1.19}. Therefore, it suffices to prove \eqref{10.5}.

In the case of flat boundary, it follows from  the boundary conditions \eqref{10.1}  that
\begin{align}\label{11.1}
u^\v_3=0,~~\omega^\v_1=0,~\omega^\v_2=0,~~ \mbox{on}~~\Gamma.
\end{align}
\begin{lemma}\label{lem11.1}
	Let \eqref{11.1} holds, then the vectors $(u^\v\cdot\nabla) \omega^\v$ and $(\omega^\v\cdot\nabla) u^\v$ are normal to $\Gamma$.
\end{lemma}
\noindent\textbf{Proof.}  It is easy to check that  $\partial_3u_1^\v=\omega_2^\v+\partial_1u_3^\v=0$ and  $\partial_3u_2^\v=\omega_1^\v+\partial_2u_3^\v=0$.  Then direct calculations yield that $(u^\v\cdot\nabla) \omega^\v$ and $(\omega^\v\cdot\nabla) u^\v$ are normal to $\Gamma$. Thus, the proof of this lemma is completed. $\hfill\Box$

\

The following formula plays a important role in the proof of uniform bounds for $\|u^\v\|_{L^\i(0,T,H^2)}$.
\begin{lemma}\label{lem11.2}
	Let $(\r^\v,u^\v)$ be smooth solution of \eqref{1.1} and \eqref{10.1}, then it holds that
	\begin{align}\label{11.2}
	\mu n\times(\nabla\times\nabla\times\omega^\v)=-n\times\left(\f{\nabla\r^\v}{\r^\v}\times(\mu\Delta u^\v+(\m+\l)\nabla\mbox{div}u^\v)\right), ~~\mbox{on}~\Gamma.
	\end{align}
\end{lemma}
\noindent\textbf{Proof.} Applying $\nabla\times$ to $\eqref{1.1}_2$ gives that
\begin{align}\label{11.3}
\r^\v\omega^\v_t+\r^\v u^\v\cdot\nabla\omega^\v-\r^\v\omega^\v\cdot\nabla u^\v +\r^\v\mbox{div}u^\v \omega^\v=-\mu\v\nabla\times\nabla\times\omega^\v-\nabla\r^\v\times(u^\v_t+u^\v\cdot\nabla u^\v).
\end{align}
Since $\omega^\v\times n=0$ on the boundary, so
\begin{align}\label{11.4}
\omega^\v_t\times n=0,~~\mbox{on}~\Gamma.
\end{align}
It follows from Lemma \ref{lem11.1} that
\begin{align}\label{11.5}
n\times((u^\v\cdot\nabla) \omega^\v)=n\times((\omega^\v\cdot\nabla) u^\v)=0,~~\mbox{on}~\Gamma.
\end{align}
Combining \eqref{10.1} and \eqref{11.3}-\eqref{11.5} gives that
\begin{align}\label{11.6}
\mu\v n\times(\nabla\times\nabla\times\omega^\v)=-n\times(\nabla\r^\v\times(u^\v_t+u^\v\cdot\nabla u^\v)),~\mbox{on}~\Gamma.
\end{align}
It follows from $\eqref{1.1}_2$ that
\begin{align}\label{11.7}
&\nabla\r^\v\times(u^\v_t+u^\v\cdot\nabla u^\v)=\f{\nabla\r^\v}{\r^\v}\times(\r^\v u^\v_t+\r^\v u^\v\cdot\nabla u^\v)\nonumber\\
&=\f{\nabla\r^\v}{\r^\v}\times (-p'(\r^\v)\nabla\r^\v+\mu\v\Delta u^\v+(\mu+\l)\v\nabla\mbox{div}u^\v)=\v\f{\nabla\r^\v}{\r^\v}\times(\mu\Delta u^\v+(\mu+\l)\nabla\mbox{div}u^\v).
\end{align}
Then, substituting \eqref{11.7} into \eqref{11.6} yields \eqref{11.2}. Therefore, the proof of this lemma is completed.
$\hfill\Box$

\

\begin{lemma}\label{lem11.3}
	It holds that
	\begin{align}\label{11.8}
	&\|\nabla\times\omega^\v(t)\|^2_{L^2}+\v\int_0^t\|\Delta\omega^\v\|^2_{L^2}d\tau\nonumber\\
	&\leq C+\|u^\v_0\|^2_{H^2}+C_\d\int_0^t\|\nabla^2u^\v\|^2d\tau+\d\v^2\int_0^t\|\nabla^3u^\v\|^2d\tau,
	\end{align}
	where $C=P(\tilde{C}_1)>0$.
\end{lemma}
\noindent\textbf{Proof.} Applying $\nabla\times$ to \eqref{11.3} shows that
\begin{align}\label{11.9}
&\r^\v(\nabla\times\omega^\v)_t+\r^\v (u^\v\cdot\nabla)(\nabla\times\omega^\v)\nonumber\\
&=-\mu\v\nabla\times\nabla\times\nabla\times\omega^\v-\nabla\r^\v\times\omega^\v_t
-[\nabla\times,\r^\v u^\v\cdot\nabla]\omega^\v\nonumber\\
&~~-\nabla\times\Big(\nabla\r^\v\times(u^\v_t+u^\v\cdot\nabla u^\v)
-\r^\v\omega^\v\cdot\nabla u^\v +\r^\v\mbox{div}u^\v \omega^\v\Big).
\end{align}
Multiplying  \eqref{11.9} by $\nabla\times\omega^\v$, one obtains that
\begin{align}\label{11.10}
&\f12\int\r^\v|\nabla\times\omega^\v|^2dx+\mu\v\int_0^t\|\nabla\times\nabla\times\omega^\v\|^2d\tau\nonumber\\
&=-\mu\v\int_0^t\int_{\Gamma}(n\times(\nabla\times\nabla\times\omega^\v))\cdot(\nabla\times\omega^\v)dyd\tau\nonumber\\
&~~~~~~~~~~~~-\int_0^t\int\Big(\nabla\r^\v\times\omega^\v_t+[\nabla\times,\r^\v u^\v\cdot\nabla]\omega^\v\Big)\cdot\nabla\times\omega^\v dxd\tau\nonumber\\
&~~~~-\int_0^t\int\Big\{\nabla\times\Big(\nabla\r^\v\times(u^\v_t+u^\v\cdot\nabla u^\v)-\r^\v\omega^\v\cdot\nabla u^\v +\r^\v\mbox{div}u^\v \omega^\v\Big)\Big\}\cdot\nabla\times\omega^\v dxd\tau.
\end{align}
First, it follows from uniform estimates of \eqref{1.18}, \eqref{1.19} and Cauchy inequality that
\begin{align}\label{11.11}
&\Big|\int_0^t\int\Big\{\nabla\times\Big(\nabla\r^\v\times(u^\v_t+u^\v\cdot\nabla u^\v)-\r^\v\omega^\v\cdot\nabla u^\v +\r^\v\mbox{div}u^\v \omega^\v\Big)\Big\}\cdot\nabla\times\omega^\v dxd\tau\Big|\nonumber\\
&+\Big|\int_0^t\int\Big(\nabla\r^\v\times\omega^\v_t+[\nabla\times,\r^\v u^\v\cdot\nabla]\omega^\v\Big)\cdot\nabla\times\omega^\v dxd\tau\Big|\nonumber\\
&\leq CP(1+\|(\r^\v,u^\v,\r^\v_t,u^\v_t)\|^2_{W^{1,\infty}})\int_0^t\|u^\v\|^2_{H^2}+\|\Delta p^\v\|^2+\|\nabla p^\v\|^2_{H^1_{co}}+\|p^\v\|^2_{H^2_{co}}d\tau\nonumber\\
&\leq C+C\int_0^t\|u^\v\|^2_{H^2}d\tau.
\end{align}
For the boundary term, it follows from Lemma \ref{11.2} that
\begin{align}\label{11.12}
&\Big|\mu\v\int_0^t\int_{\Gamma}(n\times(\nabla\times\nabla\times\omega^\v))\cdot(\nabla\times\omega^\v)dyd\tau\Big|\nonumber\\
&=\Big|\v\int_0^t\int_{\Gamma}\Big\{n\times\left(\f{\nabla\r^\v}{\r^\v}\times(\mu\Delta u^\v+(\mu+\l)\nabla\mbox{div}u^\v)\right)\Big\}\cdot(\nabla\times\omega^\v)dyd\tau\Big|\nonumber\\
&\leq C\v\int_0^t\int_{\Gamma} |\nabla^2 u^\v|^2dyd\tau\leq C\v\int_0^t\|\nabla^3u^\v\|\|\nabla^2u^\v\|+\|\nabla^2u^\v\|^2d\tau\nonumber\\
&\leq \d\v^2\int_0^t\|\nabla^3u^\v\|^2d\tau+C_\d C\int_0^t\|\nabla^2u^\v\|^2d\tau.
\end{align}
Substituting \eqref{11.12} and \eqref{11.11} into \eqref{11.10} and noting that $\Delta\omega^\v=-\nabla\times\nabla\times\omega^\v$, one proves \eqref{11.8}. Therefore, the proof of this lemma is completed. $\hfill\Box$

\

\begin{lemma}
	It holds that
	\begin{align}\label{11.13}
	&\|\nabla \mbox{div}u^\v\|^2\leq  C< +\infty,
	\end{align}
	and
	\begin{align}\label{11.14}
	&\|\nabla^2 \mbox{div}u^\v\|^2\leq  C(1+\|u^\v\|^2_{H^2}).
	\end{align}
	where $C=P(\tilde{C}_1)>0$.
\end{lemma}
\noindent\textbf{Proof.} Since
\begin{equation}\label{11.15}
\mbox{div}u^\v=-\f{p^\v_t}{\g p^\v}-\f{u^\v}{\g p^\v}\cdot\nabla p^\v,
\end{equation}
then, it follows from the uniform estimates \eqref{1.18} and \eqref{1.19} that
\begin{align}\label{11.16}
&\|\nabla \mbox{div}u^\v\|^2\leq CP(\|(p^\v,u^\v,p^\v_t,\nabla u^\v,\nabla p^\v)\|^2_{L^\infty}+\|\nabla p^\v\|^2_{H^1_{co}}+\|u^\v_3\partial_z^2p^\v\|^2_{L^2})\nonumber\\
&\leq  P(\|(p^\v,u^\v,p^\v_t,\nabla u^\v,\nabla p^\v)\|^2_{L^\infty}+\|\nabla u^\v_3\|^2_{L^\infty}\|\nabla p^\v\|^2_{H^1_{co}})\leq P(\tilde{C}_1),
\end{align}
and
\begin{align}\label{11.17}
&\|\nabla^2 \mbox{div}u^\v\|^2\leq P(1+\|(p^\v,u^\v,p^\v_t,\nabla u^\v,\nabla p^\v)\|^2_{L^\infty})\cdot(\|\Delta p^\v\|^2_{\mathcal{H}^1}+\| p^\v\|^2_{\mathcal{H}^3}+\|u^\v\|^2_{H^2}+\|u^\v_3\partial_z^3 p^\v\|^2)\nonumber\\
&\leq P(\tilde{C}_1)(1+\|u^\v\|^2_{H^2}+\|\nabla u^\v_3\|^2_{L^\infty}\|\Delta p^\v\|^2_{H^1_{co}})\leq P(\tilde{C}_1)(1+\|u^\v\|^2_{H^2})
\end{align}
where one has used $u_3|_{\Gamma}=0$ in \eqref{11.16} and \eqref{11.17}.  Thus, \eqref{11.13} and \eqref{11.14} are proved. Therefore, the proof of this lemma is completed. $\hfill\Box$

\

{\bf Proof of \eqref{10.5}:~}  It follows from \eqref{1.18} that  $\|\nabla u^\v\|_{H^1_{co}}$ and $\|u^\v\|_{H^2_{co}}$ are bounded uniformly.
Note that $\Delta u^\v=-\nabla\times\omega^\v+\nabla\mbox{div}u^\v$, one obtains that
\begin{align}\label{11.18}
\|u^\v\|^2_{H^2}&\leq C(\|\nabla\times\omega^\v\|^2+\|\nabla\mbox{div}u^\v\|^2+\|u^\v\|^2_{H^2_{co}}+\|\nabla u^\v\|^2_{H^1_{co}})\nonumber\\
&\leq C(1+\|\nabla\times\omega^\v\|^2),
\end{align}
where one has used \eqref{11.13} and the uniform estimate \eqref{1.18} in the last inequality.
On the other hand, it follows from \eqref{3.1} that
\begin{align}\label{11.19}
\|u^\v\|^2_{H^3}&\leq C\left(\|\omega^\v\|^2_{H^2}+\|\mbox{div}u^\v\|^2_{H^2}+\|u^\v\|^2_{L^2}+|u^\v\cdot n|^2_{H^{\f52}} \right)\nonumber\\
&\leq C\left(\|\partial^2_z\omega^\v\|^2_{L^2}+\|\partial^2_z\mbox{div}u^\v\|^2_{L^2}+\|\nabla^2u^\v\|^2_{H^1_{co}}+\|\nabla u^\v\|^2_{H^2_{co}}+\|u^\v\|^2_{H^3_{co}} \right)\nonumber\\
&\leq C\left(\|\Delta\omega^\v\|^2_{L^2}+\|\partial^2_z\mbox{div}u^\v\|^2_{L^2}+\|\nabla^2u^\v\|^2_{H^1_{co}}+\|\nabla u^\v\|^2_{H^2_{co}}+\|u^\v\|^2_{H^3_{co}} \right)\nonumber\\
&\leq C\left(1+\|\Delta\omega^\v\|^2_{L^2}+\|\nabla^2u^\v\|^2_{H^1_{co}} \right),
\end{align}
where \eqref{11.14} and the uniform estimate \eqref{1.18} have been used in the last inequality.

Therefore, combining \eqref{11.8}, \eqref{11.18} and \eqref{11.19}, one obtains that
\begin{align*}
&\|u^\v\|^2_{H^2}+\v\int_0^t\|u^\v\|^2_{H^3}d\tau\nonumber\\
&\leq C+\|u^\v_0\|^2_{H^2}+\d\v\int_0^t\|u^\v\|^2_{H^3}d\tau+\v\int_0^t\|\nabla^2u^\v\|^2_{H^1_{co}}d\tau
+C_\d\int_0^t\|\nabla^2u^\v\|^2d\tau.
\end{align*}
Taking $\d$ suitably small and using the uniform estimate \eqref{1.18}, one has that
\begin{align*}
\|u^\v\|^2_{H^2}+\v\int_0^t\|u^\v\|^2_{H^3}d\tau\leq C+\|u^\v_0\|^2_{H^2}
+C\int_0^t\|\nabla^2u^\v\|^2d\tau,
\end{align*}
Then, it follows from Gronwall inequality that
\begin{align}\label{11.20}
\|u^\v\|^2_{H^2}+\v\int_0^t\|u^\v\|^2_{H^3}d\tau\leq\exp(\tilde{C}_1)(1+ \|u^\v_0\|^2_{H^2}),
\end{align}
which proves \eqref{10.5}. Therefore, the proof of Theorem \ref{thm1.3} is completed.  $\hfill\Box$


\section{Proof of Theorem \ref{thm1.2}: Inviscid Limit}

In this section, we  study  the vanishing viscosity limit of viscous solutions to the inviscid one with a rate of convergence. It is well known that the solution $(\r,u)(t)\in H^3$ of Euler system  \eqref{1.7}, \eqref{1.8} with initial data  $(\r_0,u_0)$ satisfies
\begin{equation}\label{15.1}
\sum_{k=0}^3\|(\r,u)\|_{C^k(0,T_1;H^{3-k})}\leq \tilde{C}_4,~~~\f{1}{2C_0}\leq \r^\v\leq 2C_0,
\end{equation}
where $\tilde{C}_4$ depends only on  $\|(\r_0,u_0)\|_{H^3}$. On the other hand,  it follows from Theorem \ref{thm1.1} that the solution $(\r^\v,u^\v)(t)$ of  Navier-Stokes system \eqref{1.1},\eqref{1.10} with initial data $(\r_0,u_0)$ satisfies
\begin{equation}\label{15.2}
\|(p(\r^\v),u^\v)(t)\|_{X_{m}^{\v}}\leq \tilde{C}_1, ~~~\f{1}{2C_0}\leq \r^{\v}(t)\leq 2C_0,~~\forall~t\in[0,T_0],
\end{equation}
where $T_0,~C_0$, and $\tilde{C}_1$  are defined in Theorem \ref{thm1.1}. In particular,  this uniform regularity implies  the following bound
\begin{equation}\label{15.3}
\|(\r^\v,u^\v)\|_{W^{1,\infty}}+\|\partial_t(\r^\v,u^\v)\|_{L^\infty}\leq \tilde{C}_1,
\end{equation}
which plays an important role in the proof of Theorem \ref{thm1.2}. Based on these uniform estimates, using the strong compactness argument as \cite{Masmoudi-R},  one can prove that
\begin{align}\label{17.1}
(\r^\v,u^\v)\rightarrow (\r,u),~\mbox{as}~\v\rightarrow0,~\mbox{in}~L^\infty.
\end{align}
For the flat case, Theorem \ref{thm1.3} and lower semi-continuity of norm imply that
\begin{align}\label{17.2}
\|u\|^2_{H^2}\leq \exp(\tilde{C}_1)(1+\|u_0\|^2_{H^2}),
\end{align}
which yields immediately that, for the flat case,  the solution $(\r,u)$ of Euler system  satisfies  an additional boundary condition, i.e.
\begin{align}\label{16.1}
n\times\omega=0,~~\mbox{on}~\Gamma.
\end{align}
In general, it is impossible for  the solution of the Euler system to satisfy \eqref{16.1}. The observation of \eqref{16.1} will enable us to obtain better convergence rate for the flat case than the general case.   For later use, we extend smoothly the normal  $n$ to $\Omega$. In particular,  for the flat case, we extend the normal $n$ such that it is constant vector in the vicinity of $\Gamma$.

\

Define
\begin{eqnarray}\label{14.1}
\phi^\v=\r^\v-\r,~~\psi^\v=u^\v-u.
\end{eqnarray}
It then follows from \eqref{1.1} and \eqref{1.7} that
\begin{equation}\label{14.2}
\begin{cases}
\phi^\v_t+\r\mbox{div}\psi^\v+u\cdot\nabla\phi^\v=R_1^\v,\\[2mm]
\r\psi^\v_t+\r u\cdot\nabla\psi^\v+\nabla(p^\v-p)+\Phi^\v=-\mu\v\nabla\times(\nabla\times\psi^\v)
+(2\mu+\l)\v\nabla\mbox{div}\psi^\v+R^\v_2,
\end{cases}
\end{equation}
where
\begin{equation}\label{14.3}
\begin{cases}
R_1^\v=-\phi^\v\mbox{div}\psi^\v-\psi^\v\nabla\phi^\v-\phi^\v\mbox{div}u-\nabla\r \psi^\v,\\
R^\v_2=-\phi^\v\psi^\v_t-\phi^\v u_t+\mu\v\Delta{u}+(\mu+\l)\v\nabla\mbox{div}u,
\end{cases}
\end{equation}
and
\begin{align}\label{14.4}
\Phi^\v= (\r^\v u^\v-\r u)\cdot\nabla u^\v=(\r^\v u^\v-\r u)\cdot\nabla\psi^\v
+(\r^\v u^\v-\r u)\cdot\nabla u.
\end{align}
The boundary conditions to \eqref{14.2} are
\begin{align}\label{16.2}
\psi^\v\cdot n=0,~~n\times(\nabla\times\psi^\v)=[B\psi^\v]_\tau+[Bu]_\tau-n\times\omega,~~\mbox{on}~\partial\Om.
\end{align}
In particular, in the flat case, it follows from \eqref{10.1} and \eqref{16.1} that
\begin{align}\label{16.3}
\psi^\v\cdot n=0,~~n\times(\nabla\times\psi^\v)=0,~~\mbox{on}~\Gamma.
\end{align}

\begin{lemma}\label{lem14.1}
	It holds that
	\begin{eqnarray}\label{14.6}
	&&\|(\phi^\v, \psi^\v)(t)\|^2
	+\v\int_0^t\|\psi^\v\|^2_{H^1}d\tau \leq
	\begin{cases}
	C\v^{\f32},~~\mbox{general case},\\[1mm]
	C\v^2,~~\mbox{flat case},
	\end{cases}
	t\in [0,T_2],
	\end{eqnarray}
	where $T_2=\min\{T_0,T_1\}$,  $C>0$ depend only on $C_0, \tilde{C}_1$ and $\tilde{C}_4$.
\end{lemma}

\noindent\textbf{Proof}: Multiplying $\eqref{14.2}_2$ by $\psi^\v$, one obtains that
\begin{eqnarray}\label{14.7}
&&\frac{d}{dt} \int_{\Omega}\f12\r|\psi^\v|^2dx+\int_{\Omega}\Phi^\v\cdot \psi^\v dx+\int_{\Omega}\nabla(p^\v-p)\cdot \psi^\v dx\nonumber\\
&&=-\mu\v\int_{\Omega}\nabla\times(\nabla\times\psi^\v)\cdot \psi^\v dx+(2\mu+\l)\v\int_{\Omega}\nabla\mbox{div}\psi^\v\cdot\psi^\v dx+\int_{\Omega}R_2^\v\cdot \psi^\v dx.
\end{eqnarray}
It is easy to check that
\begin{equation}\label{14.8}
\left|\int_{\Omega}\Phi^\v\cdot\psi^\v dx\right|=\left|\int_{\Omega}((\r^\v u^\v-\r u)\cdot\nabla) u^\v \cdot\psi^\v dx\right|
\leq C(1+\|(\r^\v,u^\v,\nabla u^\v)\|_{L^\infty})\|(\phi^\v,\psi^\v)\|_{L^2}^2,
\end{equation}
and
\begin{eqnarray}\label{14.9}
&&\int_{\Omega}\nabla(p^\v-p)\cdot \psi^\v dx
=-\int_{\Omega} (p^\v-p)\mbox{div}\psi^\v dx\nonumber\\
&&\geq\int_{\Omega}\f1\r p'(\r) \phi^\v[\phi^\v_t+u\cdot\nabla\phi^\v-R_1^\v]dx-C(1+\|\nabla u^\v\|_{L^\infty})\|\phi^\v\|^2\nonumber\\
&&\geq \frac{d}{dt}\int_{\Omega}\f{p'(\r)}{2\r}|\phi^\v|^2dx
-C(1+\|(\r, u, \r^\v, u^\v)\|_{W^{1\infty}})\|\phi^\v\|^2\nonumber\\
&&\geq \frac{d}{dt}\int_{\Omega}\f{p'(\r)}{2\r}|\phi^\v|^2dx
-C\|\phi^\v\|^2.
\end{eqnarray}
Next, \eqref{16.2} implies that
\begin{eqnarray}\label{14.10}
&&-\mu\v\int_{\Omega}\nabla\times(\nabla\times\psi^\v)\cdot \psi^\v dx=-\mu\v\int_{\Omega}|\nabla\times\psi^\v|^2dx
-\mu\v\int_{\partial\Omega}n\times(\nabla\times\psi^\v)\cdot \psi^\v dx\nonumber\\
&&\leq -\mu\v\|\nabla\times\psi^\v\|^2
+C\v\left|\int_{\partial\Omega}[B\psi^\v+Bu-n\times\omega]\cdot \psi^\v dx\right|\nonumber\\
&&\leq
-\mu\v\|\nabla\times\psi^\v\|^2+C\v\left(|\psi^\v|^2_{L^2(\partial\Omega)}
+|\psi^\v|_{L^2(\partial\Omega)}\right).
\end{eqnarray}
For the flat case, it follows from \eqref{16.3} that
\begin{eqnarray}\label{14.10-1}
&&-\mu\v\int_{\Omega}\nabla\times(\nabla\times\psi^\v)\cdot \psi^\v dx\nonumber\\
&&=-\mu\v\int_{\Omega}|\nabla\times\psi^\v|^2dx
-\mu\v\int_{\Gamma}n\times(\nabla\times\psi^\v)\cdot \psi^\v dx=
-\mu\v\|\nabla\times\psi^\v\|^2.
\end{eqnarray}
It is easy to obtain that
\begin{eqnarray}\label{14.11}
&&\v\int_{\Omega}\nabla\mbox{div}\psi^\v\cdot\psi^\v dx
=-\v\|\mbox{div}\psi^\v\|^2,
\end{eqnarray}
and
\begin{eqnarray}\label{14.12}
&&\left|\int_{\Omega}R_2^\v\cdot \psi^\v dx \right|\leq C\|(\phi^\v,\psi^\v)\|_{L^2}^2+C\v^2.
\end{eqnarray}
Collecting all the above estimates, one gets that
\begin{eqnarray}\label{14.6-1}
&&\frac{d}{dt} \left(\f12\int_{\Omega}\r|\psi^\v|^2+\f{p'(\r)}{\r}|\phi^\v|^2dx\right)
+\mu\v\|\nabla\times\psi^\v\|^2+(2\mu+\l)\v\|\mbox{div}\psi^\v\|^2\nonumber\\
&&\leq
\begin{cases}
C\|(\phi^\v,\psi^\v)\|_{L^2}^2+C\v^2+C\v\left(|\psi^\v|^2_{L^2}+|\psi^\v|_{L^2}\right),~\mbox{general case},\\[2mm]
C\|(\phi^\v,\psi^\v)\|_{L^2}^2+C\v^2,~\mbox{flat case}.
\end{cases}
\end{eqnarray}
It follows from  \eqref{3.1}  that
\begin{equation}\label{14.6-2}
\|\psi^\v\|^2_{H^1}\leq C_1 \left(\|\nabla\times\psi^\v\|^2+\|\mbox{div}\psi^\v\|^2+\|\psi^\v\|^2\right).
\end{equation}
The trace theorem yields that
\begin{equation}\label{14.6-3}
|\psi^\v|^2_{L^2}\leq \d\|\nabla\psi^\v\|^2+C_\d\|\psi^\v\|^2,
\end{equation}
and
\begin{eqnarray}\label{14.6-4}
\v|\psi^\v|_{L^2}\leq \d\v\|\nabla\psi^\v\|^2+C_\d\v\|\psi^\v\|^{\f23}\leq\d\v\|\nabla\psi^\v\|^2+ \|\psi^\v\|^2+C_\d\v^{\f32}.
\end{eqnarray}
Substituting \eqref{14.6-2}-\eqref{14.6-4} into \eqref{14.6-1} and taking $\d$ suitably small, one gets that
\begin{align}\label{4.6-5}
\frac{d}{dt} \left(\int_{\Omega}\r|\psi^\v|^2+\f{p'(\r)}{\r}|\phi^\v|^2dx\right)
+c_0\v\|\psi^\v\|^2_{H^1}
\leq
\begin{cases}
C\|(\phi^\v,\psi^\v)\|_{L^2}^2+C\v^{\f32},~\mbox{general case}.\\[2mm]
C\|(\phi^\v,\psi^\v)\|_{L^2}^2+C\v^2,~\mbox{flat case}.
\end{cases}
\end{align}
Then, \eqref{14.6} follows from the Gronwall inequality.  Therefore,  the proof of Lemma \ref{lem14.1} is completed. $\hfill\Box$

\begin{lemma}\label{lem14.2}
	It holds that
	\begin{align}\label{14.13}
	&\|(\mbox{div}\psi^\v,\nabla(p^\v-p))(t)\|^2_{L^2}
	+(2\mu+\l)\v\int_0^t \|\nabla\mbox{div}\psi^\v(\tau)\|^2_{L^2}d\tau\nonumber\\
	&\leq
	\begin{cases}
	\d\int_0^t\|\psi^\v_t\|^2d\tau+C_\d\int_0^t \|(\phi^\v,\psi^\v)\|^2_{H^1} d\tau+C_\d\v^{\f12},~\mbox{general case},\\[2mm]
	\d\int_0^t\|\psi^\v_t\|^2d\tau+C_\d\int_0^t \|(\phi^\v,\psi^\v)\|^2_{H^1} d\tau+C_\d\v^{\f32},~\mbox{flat case},
	\end{cases}
	~t\in[0,T_2],
	\end{align}
where $\d>0$ will be chosen  later.	
\end{lemma}

\noindent\text{\bf Proof:} Multiplying $\eqref{14.2}_2$ by $\nabla\mbox{div}\psi^\v$, one obtains that
\begin{align}\label{14.14}
&\int \left(\r\psi^\v_t+\r u\cdot\nabla\psi^\v\right)\cdot\nabla\mbox{div}\psi^\v dx+\int\nabla(p^\v-p)\cdot\nabla\mbox{div}\psi^\v dx\\
&=-\mu\v\int\nabla\times(\nabla\times\psi^\v)\cdot \nabla\mbox{div}\psi^\v dx
+(2\mu+\l)\v\|\nabla\mbox{div}\psi^\v\|^2+\int R^\v_2\cdot\nabla\mbox{div}\psi^\v dx-\int\Phi^\v\cdot\nabla\mbox{div}\psi^\v dx.\nonumber
\end{align}
First, it follows from \eqref{6.4-1} and \eqref{1.18} that
\begin{eqnarray}\label{14.21-2}
\|\nabla\mbox{div}u^\v\|_{L^\infty}+\|\nabla\mbox{div}u^\v\|_{L^2}\leq C< \infty,
\end{eqnarray}
where $C>0$ depends only on $\tilde{C}_1$.  Integrating by parts and using Holder inequality, one has that
\begin{eqnarray}\label{14.15}
&&\int \left(\r\psi^\v_t+\r u\cdot\nabla\psi^\v\right)\cdot\nabla\mbox{div}\psi^\v dx\leq -\int \left(\r\mbox{div}\psi^\v_t+\r u\cdot\nabla\mbox{div}\psi^\v\right)\mbox{div}\psi^\v dx
\nonumber\\
&&~~~~~~~+|\int \left(\nabla\r \psi^\v_t+\nabla(\r u)^t\nabla\psi^\v\right)\mbox{div}\psi^\v dx|+|\int_{\partial\Omega}\r (u\cdot\nabla)\psi^\v\cdot n \mbox{div}\phi^\v d\sigma|\nonumber\\
&&\leq -\f{d}{dt}\int\f\r2|\mbox{div}\psi^\v|^2dx+\d\|\psi^\v_t\|^2+C_\d \|\nabla\psi^\v\|^2+\left|\int_{\partial\Omega}\r (u\cdot\nabla)n \psi^\v \mbox{div}\psi^\v dx\right|\nonumber\\
&&\leq \begin{cases}
-\f{d}{dt}\int\f\r2|\mbox{div}\psi^\v|^2dx+\d\|\psi^\v_t\|^2+C_\d \|\nabla\psi^\v\|^2+C|\psi^\v|_{L^2(\partial\Omega)},~\mbox{general case},\\
-\f{d}{dt}\int\f\r2|\mbox{div}\psi^\v|^2dx+\d\|\psi^\v_t\|^2+C_\d \|\nabla\psi^\v\|^2,~\mbox{flat case},
\end{cases}
\end{eqnarray}
and
\begin{eqnarray}\label{14.16}
&&\left|\int\Phi^\v\cdot\nabla\mbox{div}\psi^\v dx\right|=\left|\int[(\r^\v u^\v-\r u)\cdot\nabla\psi^\v+(\r^\v u^\v-\r u)\cdot\nabla u]\nabla\mbox{div}\psi^\v dx\right|\nonumber\\
&&\leq C(1+\|(\r^\v,u^\v)\|_{W^{1,\infty}})\|(\phi^\v,\psi^\v)\|^2_{H^1}+\left|\int_{\partial\Om}((\r^\v u^\v-\r u)\cdot\nabla u)\cdot n  \mbox{div}\psi^\v d\sigma\right|\nonumber\\
&&\leq C(1+\|(\r^\v,u^\v)\|_{W^{1,\infty}})\|(\phi^\v,\psi^\v)\|^2_{H^1}+\left|\int_{\partial\Om}((\r^\v u^\v-\r u)\cdot\nabla n)\cdot u  \mbox{div}\psi^\v d\sigma\right|\nonumber\\
&&\leq
\begin{cases}
C(1+\|(\r^\v,u^\v)\|_{W^{1,\infty}})\left(\|(\phi^\v,\psi^\v)\|^2_{H^1}
+|(\phi^\v,\psi^\v)|_{L^2(\partial\Omega)}\right),~\mbox{general case},\\
C(1+\|(\r^\v,u^\v)\|_{W^{1,\infty}})\|(\phi^\v,\psi^\v)\|^2_{H^1},~\mbox{flat case},\\
\end{cases}
\end{eqnarray}
where one has used that  $\nabla n$ is a zero matrix in the vicinity of $\Gamma$ for the flat case.

\

Rewrite  $\eqref{14.2}_1$ as
\begin{equation}\label{14.17}
(p^\v-p)_t+u\cdot\nabla(p^\v-p)+\gamma p(\r^\v)\mbox{div}\psi^\v
=-\psi^\v\cdot \nabla p-\gamma (p^\v-p)\mbox{div}u,
\end{equation}
which implies  immediately
\begin{eqnarray}\label{14.17-1}
&&\nabla\mbox{div}\psi^\v=-\f1{\gamma p^\v}\Big[\nabla(p^\v-p)_t+(u^\v\cdot\nabla)\nabla(p^\v-p)+\gamma\nabla p^\v\mbox{div}\psi^\v\nonumber\\
&&~~~~~~~~~~~~~~~~~~+\nabla u^\v\nabla(p^\v-p)+\nabla(\psi^\v\cdot\nabla p)+\gamma\nabla((p^\v-p)\mbox{div}u)\Big].
\end{eqnarray}
Using \eqref{14.17-1}, one obtains that
\begin{eqnarray}\label{14.18}
&&\int\nabla(p^\v-p)\cdot\nabla\mbox{div}\psi^\v dx \nonumber\\
&&\leq -\int\f{1}{\gamma p^\v}\nabla(p^\v-p)\Big[\nabla(p^\v-p)_t+(u^\v\cdot\nabla)\nabla(p^\v-p)\Big]dx\nonumber\\
&&~~~~~~~~~~+C(1+\|( u^\v, p^\v)\|_{W^{1\infty}})\|(p^\v-p,\psi^\v)\|^2_{H^1}\nonumber\\
&&\leq -\f{d}{dt}\int\f{1}{2\gamma p^\v}|\nabla(p^\v-p)|^2dx+C(1+\|( u^\v, p^\v)\|_{W^{1\infty}})\|(p^\v-p,\psi^\v)\|^2_{H^1}.
\end{eqnarray}
It follows from the integrating by parts  that
\begin{eqnarray}\label{14.20}
&&\v\left|\int\nabla\times(\nabla\times\psi^\v)\cdot \nabla\mbox{div}\psi^\v dx\right|=\v\left|\int_{\partial\Omega}n\times(\nabla\times\psi^\v)\cdot \nabla\mbox{div}\psi^\v d\sigma \right|\nonumber\\
&&=\v\left|\int_{\partial\Omega}\left(B\psi^\v+Bu-n\times\omega\right)\cdot \Pi(\nabla\mbox{div}\psi^\v) d\sigma \right|\leq C\v\left(1+|\psi^\v|_{H^{\f12}}\right)|\mbox{div}\psi^\v|_{H^{\f12}}\nonumber\\
&&\leq C\v \|\mbox{div}\psi^\v\|_{H^1} (1+\|\psi^\v\|_{H^1})\leq \d\v\|\nabla\mbox{div}\psi^\v\|^2+C_\d \v(1+\|\psi^\v\|_{H^1}^2).
\end{eqnarray}
For the flat case, it follows from \eqref{16.3} that
\begin{eqnarray}\label{14.20-1}
&&\v\left|\int\nabla\times(\nabla\times\psi^\v)\cdot \nabla\mbox{div}\psi^\v dx\right|=\v\left|\int_{\Gamma}n\times(\nabla\times\psi^\v)\cdot \nabla\mbox{div}\psi^\v dy \right|=0.
\end{eqnarray}

\

For the term involving $R^\v_2$.  It follows from \eqref{14.21-2} and  integrating by parts  that
\begin{align}\label{14.21}
&\left|\int R^\v_2\nabla\mbox{div}\psi^\v dx\right|
\leq C(1+\|\nabla\mbox{div}u^\v\|_{L^\infty})[\|\phi^\v\|\|\psi^\v_t\|
+\|(\phi^\v,\psi^\v)\|^2_{H^1}]\nonumber\\
&~~~~+\v\|u\|_{H^3}\|\psi^\v\|_{H^1}+\left|\int_{\partial\Om}(\mu\v\Delta{u}+(\mu+\l)\v\nabla\mbox{div}u)\cdot n \mbox{div}\psi^\v d\sigma\right|\nonumber\\
&\leq \d \|\psi^\v_t\|^2+ C_\d[\|(\phi^\v,\psi^\v)\|^2_{H^1} +\v^2]+C\v\|u\|_{H^3}\|\mbox{div}\psi^\v\|^{\f12}_{H^1}\|\mbox{div}\psi^\v\|^{\f12}_{L^2}\nonumber\\
&\leq \f{1}{8}(2\m+\l)\v\|\nabla\mbox{div}\psi^\v\|^2+\d \|\psi^\v_t\|^2+ C_\d[\|(\phi^\v,\psi^\v)\|^2_{H^1} +\v^{\f32}].
\end{align}
It follows from the trace theorem that
\begin{align}\label{16.4}
|(\phi^\v,\psi^\v)|_{L^2}\lesssim \|(\phi^\v,\psi^\v)\|^2_{H^1}+ \|(\phi^\v,\psi^\v)\|^{\f23}_{L^2}\lesssim \|(\phi^\v,\psi^\v)\|^2_{H^1}+ \v^{\f12}.
\end{align}
Collecting  all the above estimates, we obtain \eqref{14.13}. Thus,   the proof of Lemma \ref{lem14.2} is completed.    $\hfill\Box$

\begin{lemma}\label{lem14.3}
	It holds that
	\begin{align}\label{14.22}
	 &\|\nabla\times\psi^\v\|^2+\v\int_0^t\|(\nabla\times\psi^\v)(\tau)\|^2_{H^1}d\tau\leq \d\|\nabla(\phi^\v,\psi^\v)\|^2_{L^2}\nonumber\\
	&
	 +C\d\int_0^t\|\psi^\v_t\|^2+\v\|\nabla^2\psi^\v\|^2d\tau+C_\d\int_0^t \|(\phi^\v,\psi^\v)\|^2_{H^1}d\tau+C_\d\v^\f{1}{2},~\mbox{general case},
	\end{align}
	and
		\begin{align}\label{14.22-1}
		&\|\nabla\times\psi^\v\|^2+\v\int_0^t\|(\nabla\times\psi^\v)(\tau)\|^2_{H^1}d\tau\nonumber\\
		&
		\leq C\d\int_0^t\|\psi^\v_t\|^2+\v\|\nabla^2\psi^\v\|^2d\tau+C_\d\int_0^t \|(\phi^\v,\psi^\v)\|^2_{H^1}d\tau+C_\d\v^\f{3}{2},~\mbox{flat case},
		\end{align}
	where $\d>0$ will be chosen later.
\end{lemma}

\noindent\text{\bf Proof:} Multiplying $\eqref{14.2}_2$ by $\nabla\times(\nabla\times\psi^\v)$, one obtains that
\begin{eqnarray}\label{14.23}
&&\int \r^\v\psi^\v_t \nabla\times(\nabla\times\psi^\v) dx+\int\nabla(p^\v-p)\cdot\nabla\times(\nabla\times\psi^\v) dx\nonumber\\
&&=-\mu\v\|\nabla\times(\nabla\times\psi^\v)\|^2
+(2\mu+\l)\v\int\nabla\times(\nabla\times\psi^\v)\cdot \nabla\mbox{div}\psi^\v dx\nonumber\\
&&~~~~~~~+\int \tilde\Phi^\v \nabla\times(\nabla\times\psi^\v) dx+\int \tilde{R}^\v_2\nabla\times(\nabla\times\psi^\v) dx,
\end{eqnarray}
where one has rewritten $\eqref{14.2}_2$ and
\begin{align}
\tilde{R}_2^\v=-\phi^\v u_t+\mu\v\Delta u+(\mu+\l)\v\nabla\mbox{div}u~~ \mbox{and}~~ \tilde\Phi^\v =\r^\v u^\v\cdot\nabla\psi^\v
+(\r^\v u^\v-\r u)\cdot\nabla u.\nonumber
\end{align}
Note that the second term on the right hand side of \eqref{14.23} has been estimated in \eqref{14.20} and \eqref{14.20-1}. It remains to estimate the other terms of \eqref{14.23}.  First, it follows from integrating by parts  that
\begin{eqnarray}\label{14.24}
&&\int \r^\v\psi^\v_t \nabla\times(\nabla\times\psi^\v) dx\nonumber\\
&&=\int \Big[\r^\v(\nabla\times\psi^\v)_t+\nabla\r^\v\times\psi^\v_t \Big]\cdot(\nabla\times\psi^\v) dx+\int_{\partial\Omega}\r^\v\psi^\v_t\cdot (n\times(\nabla\times\psi^\v))d\sigma\nonumber\\
&&\geq \f{d}{dt}\int\f12\r^\v|(\nabla\times\psi^\v)|^2dx+\int_{\partial\Omega}\r^\v\psi^\v_t[B\psi^\v+Bu-n\times\omega]d\sigma -\d\|\psi^\v_t\|^2-C_\d\|\psi^\v\|^2_{H^1}\nonumber\\
&&\geq \f{d}{dt}\left(\int\f12\r^\v|(\nabla\times\psi^\v)|^2dx+\int\f12\r^\v\psi^\v B\psi^\v +\r^\v\psi^\v(Bu-n\times\omega)d\sigma  \right)\nonumber\\
&&~~~~~~~~~~~~-\d\|\psi^\v_t\|^2-C_\d\left(\|\psi^\v\|^2_{H^1}+|\psi^\v|^2_{L^2}
+|\psi^\v|_{L^2}\right)\nonumber\\
&&\geq \f{d}{dt}\left(\int\f12\r^\v|\nabla\times\psi^\v|^2dx+\int\f12\r^\v\psi^\v B\psi^\v +\r^\v\psi^\v(Bu-n\times\omega)d\sigma  \right)\nonumber\\
&&~~~~~~~~~~~~-\d\|\psi^\v_t\|^2-C_\d\left(\|\psi^\v\|^2_{H^1}
+|\psi^\v|_{L^2}\right).
\end{eqnarray}
For the flat case,  it follows from \eqref{16.3} that
\begin{eqnarray}\label{14.24-1}
&&\int \r^\v\psi^\v_t \nabla\times(\nabla\times\psi^\v) dx\nonumber\\
&&=\int \Big[\r^\v(\nabla\times\psi^\v)_t+\nabla\r^\v\times\psi^\v_t \Big]\cdot(\nabla\times\psi^\v) dx+\int_{\partial\Omega}\r^\v\psi^\v_t\cdot (n\times(\nabla\times\psi^\v))d\sigma\nonumber\\
&&\geq \f{d}{dt}\left(\int\f12\r^\v|\nabla\times\psi^\v|^2dx \right)-\d\|\psi^\v_t\|^2-C_\d\|\psi^\v\|^2_{H^1}.
\end{eqnarray}
Integrating along the boundary, one has that
\begin{eqnarray}\label{14.25}
&&\left|\int\nabla(p^\v-p)\cdot\nabla\times(\nabla\times\psi^\v) dx\right|
=\left|\int_{\partial\Omega}\nabla(p^\v-p)\cdot (n\times(\nabla\times\psi^\v)) d\sigma\right|\nonumber\\
&&=\left|\int_{\partial\Omega}\Pi(\nabla(p^\v-p))\cdot[B\psi^\v+Bu-n\times\omega] dx\right|\leq C\left[|p^\v-p|_{H^{\f12}}|\psi^\v|_{H^{\f12}}+|p^\v-p|_{L^2}\right]\nonumber\\
&&\leq C\left[\|p^\v-p\|^2_{H^1}+\|\psi^\v\|^2_{H^1}+|p^\v-p|_{L^2}\right].
\end{eqnarray}
For the flat case,  it follows from \eqref{16.3} that
\begin{eqnarray}\label{14.25-1}
&&\left|\int\nabla(p^\v-p)\cdot\nabla\times(\nabla\times\psi^\v) dx\right|
=\left|\int_{\partial\Omega}\nabla(p^\v-p)\cdot n\times(\nabla\times\psi^\v) d\sigma\right|=0.
\end{eqnarray}
For the term involving $\tilde{R}^\v_2$,   integrating by parts leads to that
\begin{eqnarray}\label{14.26}
&&\left|\int \tilde{R}^\v_2\nabla\times(\nabla\times\psi^\v) dx\right|\nonumber\\
&&\leq C\|(\phi^\v,\psi^\v)\|^2_{H^1}+\left|\int_{\partial\Omega}\phi^\v u_t\cdot (n\times(\nabla\times\psi^\v)) \right|+C\v\|u\|_{H^3}\|\psi^\v\|_{H^1}+C\v\|u\|_{H^3}|\nabla\times\psi^\v|_{L^2}
\nonumber\\
&&\leq C\|(\phi^\v,\psi^\v)\|^2_{H^1}+\d\v\|\nabla(\nabla\times\psi^\v)\|^2
+\left|\int_{\partial\Omega}\phi^\v u_t\cdot [B\psi^\v+Bu-n\times\omega] \right|+C_\d(\v^{\f32}+\|\psi^\v\|^2_{H^1})\nonumber\\
&&\leq \d\v\|\nabla\times(\nabla\times\psi^\v)\|^2
+C_\d\left(\|(\phi^\v,\psi^\v)\|^2_{H^1}+|(\phi^\v,\psi^\v)|_{L^2}
+\v^{\f32} \right).
\end{eqnarray}
For the flat case, it follows from \eqref{16.3} that
\begin{eqnarray}\label{14.26-1}
&&\left|\int \tilde{R}^\v_2\nabla\times(\nabla\times\psi^\v) dx\right|\nonumber\\
&&\leq C\|(\phi^\v,\psi^\v)\|^2_{H^1}+\left|\int_{\partial\Omega}\phi^\v u_t\cdot (n\times(\nabla\times\psi^\v)) \right|+C\v\|u\|_{H^3}\|\psi^\v\|_{H^1}+C\v\|u\|_{H^3}|\nabla\times\psi^\v|_{L^2}
\nonumber\\
&&\leq \d\v\|\nabla\times(\nabla\times\psi^\v)\|^2
+C_\d\left(\|(\phi^\v,\psi^\v)\|^2_{H^1}+\v^{\f32} \right).
\end{eqnarray}

\

For the term involving $\tilde\Phi^\v$, we will follow the ideas in [26]. Indeed, it follows from the integrating by parts  that
\begin{eqnarray}\label{14.27}
&&\left|\int \tilde\Phi^\v \nabla\times(\nabla\times\psi^\v) dx\right|\nonumber\\
&&\leq \left|\int \nabla\times\tilde\Phi^\v\cdot (\nabla\times\psi^\v) dx\right|+\left|\int_{\partial\Omega} \tilde\Phi^\v (B\psi^\v)_\tau d\sigma\right|+\left|\int_{\partial\Omega} \tilde\Phi^\v[(Bu)_\tau-n\times\omega] d\sigma\right|\nonumber\\
&&\doteq N+BN+BNL.
\end{eqnarray}
Noting
\begin{eqnarray}\label{14.28}
\nabla\times((a\cdot\nabla)b)&=&(a\cdot\nabla)(\nabla\times b)+
\left(\begin{array}{cc} \partial_2 a\cdot \nabla b_3-\partial_3a\cdot\nabla b_2\\  \partial_3 a\cdot \nabla b_1-\partial_1a\cdot\nabla b_3\nonumber\\
\partial_1 a\cdot \nabla b_2-\partial_2a\cdot\nabla b_1
\end{array}
\right)\\
&=&(a\cdot\nabla)(\nabla\times b)+(\nabla a)^{\perp}\cdot\nabla b,
\end{eqnarray}
which implies that
\begin{eqnarray}\label{14.29}
&&\nabla\times\tilde\Phi^\v=\r^\v u^\v\cdot\nabla(\nabla\times\psi^\v)+(\nabla (\r^\v u^\v))^{\perp}\cdot\nabla\psi^\v+(\r^\v u^\v-\r u)\cdot\nabla\omega\nonumber\\
&&~~~~~~~~~~~~~~+(\nabla(\r^\v u^\v-\r u))^{\perp}\cdot\nabla u.
\end{eqnarray}
Integrating by parts and using Cauchy inequality, one has  that
\begin{eqnarray}\label{14.30}
N\leq C(1+\|(\r^\v,u^\v)\|^2_{W^{1,\infty}})\|(\phi^\v,\psi^\v)\|^2_{H^1}.
\end{eqnarray}

\

For the term $BN$, it holds that
\begin{eqnarray}\label{14.31}
&&BN=\left|\int_{\partial\Omega} \tilde\Phi^\v (B\psi^\v)_\tau d\sigma\right|
=\left|\int_{\partial\Omega} n\times\tilde\Phi^\v \cdot [n\times(B\psi^\v)_\tau] d\sigma\right|\nonumber\\
&&~~~=\left|\int_{\partial\Omega} n\times\tilde\Phi^\v \cdot [n\times(B\psi^\v)] d\sigma\right|\nonumber\\
&&~~~=\left|\int\nabla\times\tilde\Phi^\v\cdot[n\times(B\psi^\v)]dx-\int\tilde\Phi^\v\cdot\nabla\times[n\times(B\psi^\v)]dx\right|.
\end{eqnarray}
It follows from  the Cauchy inequality, integrating by parts and \eqref{14.29} that
\begin{eqnarray}\label{14.32}
&&\left|\int\nabla\times\tilde\Phi^\v\cdot[n\times(B\psi^\v)]dx\right|\nonumber\\
&&\leq \left|\int\r^\v (u^\v\cdot\nabla)(\nabla\times\psi^\v)\cdot[n\times(B\psi^\v)]dx\right|
+C(1+\|(\r^\v,u^\v)\|^2_{W^{1,\infty}})\|(\phi^\v,\psi^\v)\|^2_{H^1}\nonumber\\
&& \leq C(1+\|(\r^\v,u^\v)\|^2_{W^{1,\infty}})\|(\phi^\v,\psi^\v)\|^2_{H^1},
\end{eqnarray}
and
\begin{eqnarray}\label{14.33}
&&\left|\int\tilde\Phi^\v\cdot\nabla\times[n\times(B\psi^\v)]dx\right|\leq C(1+\|(\r^\v,u^\v)\|^2_{W^{1,\infty}})\|(\phi^\v,\psi^\v)\|^2_{H^1}.
\end{eqnarray}
Substituting \eqref{14.32} and \eqref{14.33} into \eqref{14.31}, one obtains that
\begin{eqnarray}\label{14.34}
&&BN= C(1+\|(\r^\v,u^\v)\|^2_{W^{1,\infty}})\|(\phi^\v,\psi^\v)\|^2_{H^1}.
\end{eqnarray}

Finally, we estimate the leading order term   $BNL$ on the boundary. Since the boundary layer may appear, in general, the term $(Bu)_\tau-n\times\omega$ is not zero.
\begin{align}\label{14.35}
&BNL=\left|\int_{\partial\Omega} \tilde\Phi^\v[(Bu)_\tau-n\times\omega] d\sigma\right|\nonumber\\
&\leq C(1+\|(\r^\v,u^\v)\|^2_{W^{1,\infty}})|(\phi^\v,\psi^\v)|_{L^2}
+ \left|\int_{\partial\Omega} \r u\cdot\nabla\psi^\v[(Bu)_\tau-n\times\omega] d\sigma\right|.
\end{align}
In order to estimate the last term in \eqref{14.35},  we note that
\begin{eqnarray}\label{14.36}
u\cdot\nabla\psi^\v=u_1\partial_{y^1}\psi^\v+u_2\partial_{y^2}\psi^\v, ~~x\in \partial\Omega.
\end{eqnarray}
It follows from  \eqref{14.36} and integrating by parts along the boundary  that
\begin{eqnarray}\label{14.37}
\left|\int_{\partial\Omega} \r u\cdot\nabla\psi^\v[(Bu)_\tau-n\times\omega] d\sigma\right|&\leq& C |\psi^\v|_{L^2}\cdot|\r u((Bu)_\tau-n\times\omega)|_{H^{1}}\leq  C |\psi^\v|_{L^2}.
\end{eqnarray}
Therefore, substituting \eqref{14.37} into  \eqref{14.35}, one gets that
\begin{eqnarray}\label{14.38}
&&BNL\leq  C(1+\|(\r^\v,u^\v)\|^2_{W^{1,\infty}})|(\phi^\v,\psi^\v)|_{L^2}.
\end{eqnarray}
Substituting \eqref{14.30}, \eqref{14.34}, \eqref{14.38} into \eqref{14.27}, one obtains that
\begin{eqnarray}\label{14.44}
&&\left|\int \tilde\Phi^\v\cdot\nabla\times(\nabla\times\psi^\v) dx\right|\leq C[\|(\phi^\v,\psi^\v)\|^2_{H^1}+|(\phi^\v,\psi^\v)|_{L^2}].
\end{eqnarray}
For the flat case, it follows from \eqref{16.3} that
\begin{align}\label{14.27-1}
&\left|\int \tilde\Phi^\v \nabla\times(\nabla\times\psi^\v) dx\right|= \left|\int \nabla\times\tilde\Phi^\v\cdot (\nabla\times\psi^\v) dx\right|+\left|\int_{\Gamma} \tilde\Phi^\v\cdot(n\times\nabla\times\psi^\v) d\sigma\right|\nonumber\\
&\leq \left|\int \nabla\times\tilde\Phi^\v\cdot (\nabla\times\psi^\v) dx\right|\leq C\|(\phi^\v,\psi^\v)\|^2_{H^1},
\end{align}
where the last inequality follows from \eqref{14.30}.

\

Combining \eqref{14.23}-\eqref{14.26-1}, \eqref{14.44}, \eqref{14.27-1} and \eqref{14.20}-\eqref{14.20-1}, we obtain that
\begin{align}\label{14.45-1}
&\f{d}{dt}\left(\int\f12\r^\v|\nabla\times\psi^\v|^2dx+\int\f12\r^\v\psi^\v B\psi^\v +\r^\v\psi^\v(Bu-n\times\omega)d\sigma  \right)
+\f12\mu\v\|\nabla\times(\nabla\times\psi^\v)\|^2\nonumber\\
&\leq
C\d\|\psi^\v_t\|^2+C\d\v\|\nabla^2\psi^\v\|^2+C_\d\left(\|(\phi^\v,\psi^\v)\|^2_{H^1}+\v^\f{1}{2}\right),~\mbox{general case},
\end{align}
and
\begin{align}\label{14.45-2}
&\f{d}{dt}\left(\int\f12\r^\v|\nabla\times\psi^\v|^2dx\right)
+\f12\mu\v\|\nabla\times(\nabla\times\psi^\v)\|^2\nonumber\\
&\leq
C\d\|\psi^\v_t\|^2+C\d\v\|\nabla^2\psi^\v\|^2+C_\d\left(\|(\phi^\v,\psi^\v)\|^2_{H^1}+\v^\f{3}{2}\right),~\mbox{flat case},
\end{align}
where we have used
\begin{eqnarray}\label{14.47}
\begin{cases}
|(\phi^\v,\psi^\v)|_{L^2}\leq \|(\phi^\v,\psi^\v)\|^{\f12}_{H^1}\cdot\|(\phi^\v,\psi^\v)\|^{\f12}
\leq  \d\|\nabla(\phi^\v,\psi^\v)\|^2+C_\d\v^{\f12},\\[2mm]
|(\phi^\v,\psi^\v)|^2_{L^2}\leq \|(\phi^\v,\psi^\v)\|_{H^1}\cdot\|(\phi^\v,\psi^\v)\|
\leq  \d\|\nabla(\phi^\v,\psi^\v)\|^2+C_\d\v^{\f32},
\end{cases}
\end{eqnarray}
which are consequences of the trace theorem and  \eqref{14.6}. It follows  from \eqref{3.1-1} that
\begin{align}\label{14.46}
&\|\nabla\times\psi^\v\|^2_{H^1}\leq C_1 \left(\|\nabla\times(\nabla\times\psi^\v)\|^2+\|\mbox{div}(\nabla\times\psi^\v)\|^2+\|\nabla\times\psi^\v\|^2+|n\times(\nabla\times\psi^\v)|^2_{H^{\f12}}\right)\nonumber\\
&\leq
\begin{cases}
C_1 (\|\nabla\times(\nabla\times\psi^\v)\|^2+\|\nabla\times\psi^\v\|^2+|B\psi^\v|^2_{H^{\f12}}+|(Bu)_\tau-n\times\omega|^2_{H^{\f12}}),~\mbox{general case}\nonumber\\
C_1 (\|\nabla\times(\nabla\times\psi^\v)\|^2+\|\nabla\times\psi^\v\|^2),~\mbox{flat case}
\end{cases}\\
&\leq
\begin{cases}
C_1 (\|\nabla\times(\nabla\times\psi^\v)\|^2+\|\psi^\v\|^2_{H^1}+C),~\mbox{general case},\\
C_1(\|\nabla\times(\nabla\times\psi^\v)\|^2+\|\psi^\v\|^2_{H^1}),~\mbox{flat case}.
\end{cases}
\end{align}
Substituting \eqref{14.46} into  \eqref{14.45-1} yields that
\begin{align}\label{14.45}
&\f{d}{dt}\left(\int\f12\r^\v|\nabla\times\psi^\v|^2dx+\int\f12\r^\v\psi^\v B\psi^\v +\r^\v\psi^\v(Bu-n\times\omega)d\sigma  \right)+c_0\v\|\nabla\times\psi^\v\|^2_{H^1}
\nonumber\\
&\leq
C\d\|\psi^\v_t\|^2+C\d\v\|\nabla^2\psi^\v\|^2+C_\d\left(\|(\phi^\v,\psi^\v)\|^2_{H^1}+\v^\f{1}{2}\right),~\mbox{general case},
\end{align}
and
\begin{align}\label{14.45-3}
&\f{d}{dt}\left(\int\f12\r^\v|\nabla\times\psi^\v|^2dx \right)+c_0\v\|\nabla\times\psi^\v\|^2_{H^1}
\nonumber\\
&\leq
C\d\|\psi^\v_t\|^2+C\d\v\|\nabla^2\psi^\v\|^2+C_\d\left(\|(\phi^\v,\psi^\v)\|^2_{H^1}+\v^\f{3}{2}\right),~\mbox{flat case}.
\end{align}
Integrating \eqref{14.45} and \eqref{14.45-3} over $[0,t]$ and using \eqref{14.47}, one  gets \eqref{14.22} and \eqref{14.22-1}, respectively.
Therefore, the proof of Lemma \ref{lem14.3} is completed. $\hfill\Box$

\

\noindent\textbf{Proof of Theorem \ref{thm1.2}:} It follows from  \eqref{3.1} that
\begin{eqnarray}\label{14.49}
\|\psi^\v\|^2_{H^1}&\leq& C\left(\|\nabla\times\psi^\v\|^2+\|\mbox{div}\psi^\v\|^2+\|\psi^\v\|^2+
|\psi^\v\cdot n|_{H^{\f12}}\right)\nonumber\\
&\leq&C\left(\|\nabla\times\psi^\v\|^2
+\|\mbox{div}\psi^\v\|^2+\|\psi^\v\|^2\right),
\end{eqnarray}
and
\begin{eqnarray}\label{14.50}
\|\psi^\v\|^2_{H^2}&\leq& C\left(\|\nabla\times\psi^\v\|^2_{H^1}+\|\mbox{div}\psi^\v\|^2_{H^1}+\|\psi^\v\|^2_{H^1}+
|\psi^\v\cdot n|_{H^{\f32}}\right)\nonumber\\
&\leq&C\left(\|\nabla\times\psi^\v\|^2_{H^1}+\|\mbox{div}\psi^\v\|^2_{H^1}
+\|\psi^\v\|^2_{H^1}\right).
\end{eqnarray}
While  $\eqref{14.2}_2$ implies that
\begin{equation}\label{14.5}
	\|\psi^\v_t\|_{L^2}^2\leq C\left(\|(\phi^\v,\psi^\v)\|_{L^2}^2
	+\v^2\|\nabla^2\psi^\v\|_{L^2}^2+\v^2 \right).
\end{equation}
Then, collecting \eqref{14.5}, \eqref{14.22}-\eqref{14.22-1}, \eqref{14.49}-\eqref{14.50}, \eqref{14.13}, \eqref{14.6} and choosing $\d$ suitably small, one obtains that
\begin{equation*}\label{14.52}
\|\nabla(\psi^\v,p^\v-p)\|^2+\v\int_0^t\|\psi^\v(\tau)\|^2_{H^2}d\tau\leq
\begin{cases}
 C\int_0^t \|\nabla(p^\v-p,\psi^\v)\|^2d\tau+C\v^\f{1}{2},~\mbox{general case},\\[2mm]
  C\int_0^t \|\nabla(p^\v-p,\psi^\v)\|^2d\tau+C\v^\f{3}{2},~\mbox{flat case},
\end{cases}
\end{equation*}
where one has used  $\|\phi^\v\|^2_{H^1}\lesssim C\|p^\v-p\|^2_{H^1}\lesssim\|\phi^\v\|^2_{H^1}$.
The Gronwall's inequality yields immediately that
\begin{eqnarray}\label{14.52-1}
\|\nabla(\psi^\v,p^\v-p)\|^2+\v\int_0^t\|\psi^\v(\tau)\|^2_{H^2}d\tau\leq
\begin{cases}
C\v^\f{1}{2},~\mbox{general case},\\
C\v^\f{3}{2},~\mbox{flat case}.
\end{cases}
\end{eqnarray}
Then, \eqref{14.6} and \eqref{14.52-1} imply \eqref{1.8-0}-\eqref{1.9-0} and \eqref{1.8-1}-\eqref{1.9-1}.  On the other hand, \eqref{1.10-2} and \eqref{1.10-3} are   immediately consequences of \eqref{1.8-0}, \eqref{1.8-1} and \eqref{15.1}-\eqref{15.3}.  Thus,  the proof of Theorem \ref{thm1.2} is completed.  $\hfill\Box$


\section{Appendix}

Let $S(t,\tau)$ be the $C^0$ evolution operator generated by the following equation
\begin{equation}\label{A.1}
[\partial_th+b_1(t,y)\partial_{y^1}h+b_2(t,y)\partial_{y^2}h+z b_3(t,y)\partial_zh]-\v d(t,y)\partial_{zz}h=0,~z>0,~t>\tau,
\end{equation}
with the boundary condition $h(t,y,0)=0$ and with the initial condition $h(\tau,y,z)=h_0(y,z)$. The coefficients are smooth and satisfies
\begin{equation}\label{A.1-1}
c_3\leq d(t,y)\leq \f1{c_3},~\mbox{and}~|b_i|\leq c_4,~i=1,2,3
\end{equation}
for some positive constant $c_3>0$ and $c_4>0$.

\

Then we have the following estimates which are generalizations of Lemma 15 in \cite{Masmoudi-R}.
\begin{lemma}\label{lemA.1}
It holds that, for $t\geq \tau\geq 0$
\begin{eqnarray}
\|S(t,\tau)h_0\|_{L^\infty}&&\leq \|h_0\|_{L^\infty},\label{A.2}\\
\|z\partial_zS(t,\tau)h_0\|_{L^\infty}&&\leq C(\|h_0\|_{L^\infty}+\|z\partial_zh_0\|_{L^\infty}),,\label{A.2-1}
\end{eqnarray}
where  $C>0$ is a uniform constant independent of the bounds for $d(t,y)$ and $b_j(t,y),~j=1,2,3$.
\end{lemma}

\noindent\textbf{Proof}. Set $h(t,y,z)=S(t,\tau)h_0$. Then $h$ solves the equation \eqref{A.1}. We first transform the half-plane problem into a problem in the whole space. Define $\tilde{h}$ as
\begin{equation}\label{A.3}
\tilde{h}(t,y,z)=h(t,y,z),~\mbox{for}~z> 0,~~~\tilde{h}(t,y,z)=-h(t,y,-z),~\mbox{for}~z<0.
\end{equation}
Then $\tilde{h}$ solves
\begin{equation*}\label{A.4}
[\partial_t\tilde{h}+b_1(t,y)\partial_{y^1}\tilde{h}+b_2(t,y)\partial_{y^2}\tilde{h}+z b_3(t,y)\partial_z\tilde{h}]-\v d(t,y)\partial_{zz}\tilde{h}=0,~z\in\mathbb{R},
\end{equation*}
with the initial condition $\tilde{h}(\tau,y,z)=\tilde{h}_0(y,z)$. Similar to \cite{Masmoudi-R}, we shall obtain the estimate by using an exact representation of the solution. Indeed, we can change the above equation into the generalized Fokker-Planck type equation, see \cite{Majda}.

Set
\begin{equation}\label{A.5}
v(t,y,z)=\tilde{h}(t,\Phi(t,\tau,y),z),
\end{equation}
where $\Phi$ is solution of
\begin{equation*}\label{A.6}
\partial_t\Phi=\left(\begin{array}{cccc} &b_1(t,\Phi)\\& b_2(t,\Phi)\end{array}\right),~~\Phi(\tau,\tau,y)=y.
\end{equation*}
Therefore, $v$ solves the equation
\begin{equation}\label{A.7}
\partial_tv+z \tilde{b}_3(t,y)\partial_zv-\v \tilde{d}(t,y)\partial_{zz}v=0,~z\in\mathbb{R},
\end{equation}
where
$$\tilde{d}(t,y)= d(t,\Phi(t,\tau,y)),~~\mbox{and}~~\tilde{b}_3(t,y)=b_3(t,\Phi(t,\tau,y)). $$
The equation \eqref{A.7} is just the one dimensional  generalized Fokker-Planck type equation with $y$ as parameter. We use the change of variables
\begin{equation}\nonumber
\tilde{z}=z e^{-\Gamma(t)},~~\tilde{t}=\v\int_{\tau}^te^{-2\Gamma(s)}\tilde{d}(s,y)ds,
\end{equation}
where $\Gamma(t)=\int_{\tau}^t\tilde{b}_3(s,y)ds$. Through this change of variables, the equation \eqref{A.7} reduces to the heat equation
\begin{equation}\nonumber
\begin{cases}
\partial_{\tilde{t}}v=\partial_{\tilde{z}\tilde{z}}v,\\
v(\tau,y,\tilde{z})=\tilde{h}_0(y,\tilde{z}).
\end{cases}
\end{equation}
Therefore, by using the standard heat kernel and transforming the variables $(\tilde{t},\tilde{z})$ into $(t,z)$, we obtain the explicit representation
\begin{equation*}\label{A.8}
v(t,y,z)=\int_{\mathbb{R}}k(t,\tau,y,z-z')\cdot \tilde{h}_0(y,z' e^{-\Gamma(t)})dz',
\end{equation*}
with
\begin{eqnarray}\label{A.9}
k(t,\tau,y,z-z')=\f1{\sqrt{4\pi\v\int_\tau^t\tilde{d}(s,y)e^{2(\Gamma(t)-\Gamma(s)) }ds}}
\exp\Big(-\f{|z-z'|^2}{4\v\int_\tau^t\tilde{d}(s,y)e^{2(\Gamma(t)-\Gamma(s)) }ds}\Big).
\end{eqnarray}
Since that the kernel $k$ is non-negative and that $\int_{\mathbb{R}}k(t,\tau,y,z)dz=1$, thus, it holds that
\begin{eqnarray}\label{A.10}
\|v\|_{L^\infty}\leq  \|\int_{\mathbb{R}}k(t,\tau,y,z')\cdot \sup_{z}|\tilde{h}_0(y,(z-z')e^{-\Gamma(t)})|dz'\|_{L^{\infty}_{t,z}}\leq \|\tilde{h}_0\|_{L^\infty},
\end{eqnarray}
which is the Maximum principle and proves \eqref{A.2}.

\

Next, we observe that
\begin{eqnarray}\nonumber
z\partial_zk(t,\tau,y,z-z')=(z-z') \partial_zk(t,\tau,y,z-z')-z'\partial_{z'}k(t,\tau,y,z-z'),
\end{eqnarray}
with
\begin{eqnarray}\nonumber
\int_{\mathbb{R}}|(z-z') \partial_zk(t,\tau,y,z-z')|dz' \lesssim 1.
\end{eqnarray}
Thus the integrating by parts, gives that
\begin{eqnarray}\label{A.10-1}
\|z\partial_zv\|_{L^\infty}&&\leq C\|\tilde{h}_0\|_{L^\infty}+  \|\int_{\mathbb{R}}k(t,\tau,y,z-z')\cdot e^{-\Gamma(t)}z'(\partial_z\tilde{h}_0)(y,z'e^{-\Gamma(t)}) dz'\|_{L^{\infty}}\nonumber\\
&&\leq C(\|\tilde{h}_0\|_{L^\infty}+\|z\partial_z\tilde{h}_0\|_{L^\infty}).
\end{eqnarray}
It follows from \eqref{A.3}, \eqref{A.5} and \eqref{A.10-1} that
\begin{eqnarray}\label{A.11}
\|z\partial_zh\|_{L^\infty}&&\leq C\|z\partial_z\tilde{h}\|_{L^\infty}\leq C\|z\partial_zv\|_{L^\infty}\leq C(\|\tilde{h}_0\|_{L^\infty}+\|z\partial_z\tilde{h}_0\|_{L^\infty})\nonumber\\
&&\leq C(\|h_0\|_{L^\infty}+\|z\partial_zh_0\|_{L^\infty}).
\end{eqnarray}
Thus the proof of Lemma \ref{lemA.1} is completed.
$\hfill\Box$

\begin{lemma}\label{lemA.2}
Let $h$ be a smooth solution to
\begin{equation}\label{A.13}
a(t,y)[\partial_th+b_1(t,y)\partial_{y^1}h+b_2(t,y)\partial_{y^2}h+z b_3(t,y)\partial_zh]-\v\partial_{zz}h=G,~z>0,~h(t,y,0)=0,
\end{equation}
for some smooth function $d(t,y)=\f1{a(t,y)}$ and vector fields $b=(b_1,b_2,b_3)^t(t,y)$  satisfying \eqref{A.1-1}. Assume that $h$ and $G$ are compactly supported in $z$. Then, it holds that:
\begin{eqnarray}\label{A.14}
\|h\|_{\mathcal{H}^{1,\infty}} &&\lesssim \|h_0\|_{\mathcal{H}^{1,\infty}}
+\int_0^t\|\f{1}{a}\|_{L^\infty} \|G\|_{\mathcal{H}^{1,\infty}}d\tau\nonumber\\
&&~~~~+\int_0^t(1+\|\f{1}{a}\|_{L^\infty})(1+\|b\|^2_{L^\infty}+\sum_{i=0}^2\|Z_i(a,b)\|^2_{L^\infty}) \|h\|_{\mathcal{H}^{1,\infty}}d\tau.
\end{eqnarray}

\end{lemma}

\noindent\textbf{Proof}. This will follow from Lemma \ref{A.1}. The estimates of $\|h\|_{L^{\infty}}$, $\|\partial_t h\|_{L^{\infty}}$ and $\|\partial_{y^i}h\|_{L^{\infty}}=\|Z_i h\|_{L^{\infty}},(i=1,2) $ follow easily from the maximum principle. Indeed, by Duhamel formula, one has that
\begin{equation}\label{A.17}
h(t)=S(t,0)h_0+\int_0^tS(t,\tau)\f{G(\tau)}{a(\tau)}d\tau.
\end{equation}
Consequently, \eqref{A.2} yields that
\begin{eqnarray}\label{A.18}
\|h(t)\|_{L^\infty}&&=\|S(t,0)h_0\|_{L^\infty}+\int_0^t\|S(t,\tau)\f{G(\tau)}{a(\tau)}\|_{L^\infty}d\tau\nonumber\\
&&\leq \|h_0\|_{L^\infty}+\int_0^t\|\f{G(\tau)}{a(\tau)}\|_{L^\infty}d\tau
\leq \|h_0\|_{L^\infty}+\int_0^t\|G\|_{L^\infty}\cdot\|\f1a\|_{L^\infty}d\tau.
\end{eqnarray}

Set $Z_0=\partial_t$. Then, applying $Z_i (i=0,1,2)$ to \eqref{A.13}, one has that
\begin{eqnarray}\label{A.19}
&&a(t,y)[\partial_tZ_ih+b_1(t,y)\partial_{y^1}Z_ih+b_2(t,y)\partial_{y^2}h+z b_3(t,y)\partial_zZ_ih]-\v\partial_{zz}Z_ih\nonumber\\
&&= Z_iG-Z_ia(t,y)\cdot[\partial_th+b_1(t,y)\partial_{y^1}h+b_2(t,y)\partial_{y^2}h+z b_3(t,y)\partial_zh]\nonumber\\
&&~~~+a(t,y)\cdot[\partial_th+Z_ib_1(t,y)\partial_{y^1}h+Z_ib_2(t,y)\partial_{y^2}h+z Z_ib_3(t,y)\partial_zh]\triangleq \mathcal{L}(\tau).
\end{eqnarray}
Consequently, \eqref{A.2} yields that
\begin{eqnarray}\label{A.20}
&&\|Z_ih(t)\|_{L^\infty}=\|S(t,0)Z_ih_0\|_{L^\infty}+\int_0^t\|S(t,\tau)\f{\mathcal{L}(\tau)}{a(\tau)}\|_{L^\infty}d\tau\nonumber\\
&& \lesssim\|Z_ih_0\|_{L^\infty}+\int_0^t\|\f{1}{a}\|_{L^\infty} \|Z_iG\|_{L^\infty}d\tau
+\int_0^t(1+\|Z_ib\|_{L^\infty})\|\MZ h\|_{L^\infty}d\tau\nonumber\\
&&~~+\int_0^t\|\f{1}{a}\|_{L^\infty}(1+\|b\|_{L^\infty})\|\MZ h\|_{L^\infty}\|Z_ia\|_{L^\infty}d\tau\nonumber\\
&& \lesssim \|Z_ih_0\|_{L^\infty}+\int_0^t\|\f{1}{a}\|_{L^\infty}\|Z_iG\|_{L^\infty}d\tau \nonumber\\
&&~~~~~~ +\int_0^t(1+\|\f{1}{a}\|_{L^\infty})(1+\|b\|^2_{L^\infty}+\|Z_i(a,b)\|^2_{L^\infty}) \|h\|_{\mathcal{H}^{1,\infty}}d\tau
\end{eqnarray}

It follows from \eqref{A.2-1} and \eqref{A.17} and the fact that  $h$ and $G$ are compactly supported in $z$ that
\begin{eqnarray}\label{A.21}
&&\|Z_3h(t)\|_{L^\infty} \lesssim \|h_0\|_{L^\infty}+\|z\partial_zh_0\|_{L^\infty}
+\int_0^t\|\f{G}{a}\|_{L^\infty}+\|\f{z\partial_zG}{a}\|_{L^\infty}d\tau\nonumber\\
&&\lesssim \|h_0\|_{\mathcal{H}^{1,\infty}}
+\int_0^t\|\f{1}{a}\|_{L^\infty}\cdot\|G\|_{\mathcal{H}^{1,\infty}}d\tau
\end{eqnarray}
Therefore, \eqref{A.18}, \eqref{A.20}, and \eqref{A.21} yield \eqref{A.14}. This completes the proof of Lemma \ref{lemA.2}.
$\hfill\Box$

\

\noindent {\bf Acknowledgments:}

Yong Wang is partially supported by National Natural
Sciences Foundation of China No. 11371064 and 11401565. Yan Yong is  supported by National Natural
Sciences Foundation of China No. 11201301. This research is also partially supported by Zheng Ge Ru Foundation, Hong Kong RGC Earmarked Research Grant CUHK4041/11P and CUHK4048/13P, NSFC/RGC Joint Research Scheme Grant N-CUHK443/14, a grant from Croucher Foundation, and Focus Area Grant at The Chinese University of Hong Kong.

\end{document}